\documentclass[12pt, a4paper]{amsart}
\usepackage[dvips]{epsfig}
\usepackage{amsmath}
\usepackage{amssymb}
\usepackage{amsfonts}
\usepackage{bbm}
\usepackage{amsthm}
\usepackage{amsbsy}
\usepackage{amsgen}
\usepackage{amscd}
\usepackage{amsopn}
\usepackage{amstext}
\usepackage{amsxtra}
\usepackage[utf8]{inputenc}
\usepackage{enumerate}
\usepackage{hyperref}
\usepackage{lineno}
\usepackage{marginnote}
\usepackage{verbatim}
\usepackage{calrsfs}
\usepackage{times, graphicx}
\usepackage{layout}
\usepackage{subcaption}
\usepackage{eso-pic}
\usepackage{lastpage}
\DeclareMathAlphabet{\pazocal}{OMS}{zplm}{m}{n}
\usepackage[foot]{amsaddr}

  \evensidemargin
\oddsidemargin
\numberwithin{equation}{section}

\makeatletter
\newcommand*\rel@kern[1]{\kern#1\dimexpr\macc@kerna}
\newcommand*\widebar[1]{%
  \begingroup
  \def\mathaccent##1##2{%
    \rel@kern{0.8}%
    \overline{\rel@kern{-0.8}\macc@nucleus\rel@kern{0.2}}%
    \rel@kern{-0.2}%
  }%
  \macc@depth\@ne
  \let\math@bgroup\@empty \let\math@egroup\macc@set@skewchar
  \mathsurround\z@ \frozen@everymath{\mathgroup\macc@group\relax}%
  \macc@set@skewchar\relax
  \let\mathaccentV\macc@nested@a
  \macc@nested@a\relax111{#1}%
  \endgroup
}
\makeatother

\newtheorem{theorem}{Theorem}[section]
\newtheorem{proposition}[theorem]{Proposition}
\newtheorem{lemma}[theorem]{Lemma}

\newtheorem{conjecture}[theorem]{Conjecture}

\theoremstyle{definition}
\newtheorem{definition}[theorem]{Definition}
\newtheorem{notation}[theorem]{Notation}
\newtheorem{remark}[theorem]{Remark}

\newcommand{\area}{\operatorname{Area}}

\newcommand{\Dc}{\mathcal{D}}

\newcommand{\Tc}{\mathcal{T}}
\newcommand{\Sc}{\pazocal{S}}

\newcommand{\Ec}{\pazocal{E}}

\newcommand{\Rc}{\mathcal{R}}

\newcommand{\Qc}{\pazocal{Q}}

\newcommand {\E} {\mathbb{E}}

\newcommand {\R} {\mathbb{R}}

\newcommand {\Z} {\mathbb{Z}}

\newcommand {\Cc} {\pazocal{C}}

\newcommand {\Pc} {\mathcal{P}}

\newcommand {\C} {\mathbb{C}}

\newcommand {\Ac} {\pazocal{A}}

\newcommand {\var} {\operatorname{Var}}

\newcommand {\cov} {\operatorname{Cov}}

\newcommand {\Vc} {\mathcal{V}}

\newcommand{\Mod}[1]{\ (\mathrm{mod}\ #1)}
\newcommand{\tmop}[1]{\ensuremath{\operatorname{#1}}}

 \allowbreak

\title[The Gauss circle problem]{Around the Gauss circle problem: Hardy's conjecture and the distribution of lattice points near circles}

\author{Stephen Lester}
\address{Department of Mathematics, King's College London, London WC2R 2LS, UK}
 \email{steve.lester@kcl.ac.uk}
\author{Igor Wigman}
\email{igor.wigman@kcl.ac.uk}

\begin{document}
\begin{abstract}

Hardy conjectured that the error term arising from approximating the number of lattice points lying in a radius-$R$ disc by its area is $O(R^{1/2+o(1)})$. One source of support for this conjecture is a folklore heuristic that uses i.i.d.
random variables to model the lattice points lying near the boundary and square-root cancellation of sums of these random variables.  We examine this heuristic by studying how these lattice points interact with one another and prove that their autocorrelation is determined in terms of a random model. Additionally, it is shown that lattice points near the boundary which are ``well separated" behave independently. We also formulate a conjecture concerning the distribution of pairs of these lattice points.


\end{abstract}
\date{\today}

\maketitle

\section{Introduction}
\label{sec:introduction}

Gauss initiated the study of the number of lattice points inside a radius-$R$ disc, $B(R)$, and showed that the number of such lattice points equals its area up to an error term of size $O(R)$. Numerous researchers have worked on improving the bound for the error term, including more recent works of Huxley ~\cite{Huxley} and Bourgain-Watt ~\cite{BoWa}, the latter of which shows that the error term is $O\left(R^{ \frac{1034}{1648}+o(1)}\right)$ (note that $\tfrac{1034}{1648}=0.6274\ldots$).  Hardy ~\cite{Hardy2} conjectured that
\begin{equation}\label{eq:Hardy conj alpha0=1/2}
N(R):=\#\{ \lambda \in \Z^2 : \lambda \in B(R) \}=\tmop{Area}(B(R))+O(R^{1/2+o(1)}).
\end{equation}
If true, this estimate would be nearly optimal, since Hardy proved the error term above is $$\neq O\left(R^{\frac12}(\log R)^{\frac14}\right),$$ which has been refined by Hafner \cite{hafner-1982} and Soundararajan \cite{sound-2003}.

To see why one might expect \eqref{eq:Hardy conj alpha0=1/2} should hold, recall Gauss' argument and consider squares centered at $\Z^2$-lattice points with unit area and sides parallel to the coordinate axes.
The area of the intersection of all the squares with $B(R)$ accounts for the main term $\tmop{Area}(B(R))$ and the error term arises
from the squares that are centered at lattice points lying in the annulus
\begin{equation}
\label{eq:Dc annulus sqrt(2) def}
\Dc_{R;1/\sqrt{2}}:=B(R+1/\sqrt{2})\setminus \overline{B(R-1/\sqrt{2})} = \{z\in\R^{2}:\: R-1/\sqrt{2} < \|z\|<R+1/\sqrt{2} \},
\end{equation}
of constant width $\sqrt{2}$. One might expect that each of these lattice points belongs ``randomly", with probability $\frac{1}{2}$ to the inner circle $B(R)$, independently of the other lattice points. Then, by
the standard square root cancellation laws (or, for example, the central limit theorem), a sequence of such growing sums of $\sim \tmop{Area}(\Dc_{R;1/\sqrt{2}})=\sqrt{8}\pi R$ i.i.d. Bernoulli random variables, is a.s. bounded by a function growing only slightly faster than $R^{1/2}$. Of course, the events that correspond to whether a lattice point lies in the inner circle
cannot be literally independent, since that would imply too far fetched conclusions, contradicting, in part, Hardy's lower bound; clearly, the relative positions of two neighboring lattice points are not independent.
However, a weaker statement, e.g. independence of the positions of sufficiently far apart lattice points, could still asymptotically hold.

In this article we explore this heuristic, and examine under what circumstances lattice points lying near the boundary of $B(R)$ interact independently from one another. Let
\begin{equation}
\label{eq:Gamma points annulus def}
\Gamma_{R}:=\Z^{2}\cap\Dc_{R;1/\sqrt{2}}
\end{equation}
be the set of lattice points lying in the annulus $\Dc_{R;1/\sqrt{2}}$. Also, let
$S=\left[-\tfrac12,\tfrac12\right]^2$ be the unit square, and for
$\lambda\in \Z^{2}$ denote $S(\lambda):=S+\lambda$ to be the unit square shifted by $\lambda$. Also, for $\lambda \in \Gamma_R$ let
\begin{equation*}
\Ac_{R}(\lambda):= \area \left(S(\lambda)\cap B(R)\right) = \left(\chi_{B(R)}*\chi_{S}\right)(\lambda)
\end{equation*}
where $\chi_{B(R)}$ and $\chi_{S}$ are the characteristic functions of the disc $B(R)$ and the square $S$ respectively, and $\chi_{B(R)}*\chi_{S}$ is their convolution. Namely, $\Ac_{R}(\lambda)$ is the area of the portion of $S(\lambda)$ lying inside $B(R)$.
The main focus of this work is the distribution, as $R\rightarrow\infty$, of the (deterministic) numbers
\begin{equation}
\label{eq:LambdaR def}
\Lambda_{R}:=\left\{\Ac_{R}(\lambda):\: \lambda\in \Gamma_{R}\right\},
\end{equation}
that appear naturally in the Gauss circle problem.


\subsection{Statement of the principal results}

\subsubsection{Distribution of the numbers $\Lambda_{R}$}
Given $R>0$, let
\begin{equation}
\label{eq:K(R) number points def}
K(R) = \# \Gamma_R.
\end{equation}
 Recall that $K(R)$ is asymptotic to
\begin{equation}
\label{eq:calKdef}
 \pazocal K(R):=\tmop{Area}(\Dc_{R;1/\sqrt{2}})=\sqrt{8}\pi R.
\end{equation}
Our first result concerns the expectation and the variance of the numbers $\Lambda_{R}$:

\begin{theorem}
\label{thm:C0 lim exist}

Let $\Gamma_{R}$ be the lattice points \eqref{eq:Gamma points annulus def}, and $K(R)$ their number \eqref{eq:K(R) number points def}.
Denote the expectation
\begin{equation*}
\Ec(R):= \frac{1}{K(R)}\sum\limits_{\lambda\in \Gamma_{R}} \Ac_{R}(\lambda)
\end{equation*}
of the numbers $\Lambda_{R}$, and
\begin{equation*}
\Vc(R):= \frac{1}{K(R)}\sum\limits_{\lambda\in \Gamma_{R}} \left(\Ac_{R}(\lambda)-\Ec(R)\right)^{2}.
\end{equation*}

\begin{enumerate}[i)]

\item One has
\begin{equation}
\label{eq:expectation individual 1/2}
\Ec(R) = \frac{1}{2} + O(R^{-1/3}).
\end{equation}

\item The variance of $\Lambda_{R}$ is given by
\begin{equation}
\label{eq:variance individual c0}
\Vc(R) = c_{0}+o(1), \qquad R \rightarrow \infty,
\end{equation}
with
\begin{equation}
\label{eq:c0 asymp var def}
c_{0}= \frac{1}{4}-\frac{1}{2\pi\sqrt{2}}\left( \frac{4}{15} +\frac{2\log(1+\sqrt{2})}{3}+\frac{2\sqrt{2}}{15}  \right) = 0.132642545\ldots >0.
\end{equation}

\end{enumerate}

\end{theorem}

While the expectation \eqref{eq:expectation individual 1/2} essentially follows from a well-known classical result due to independently Voronoi,
Sierpinski, and van der Corput, the variance result \eqref{eq:variance individual c0} is entirely new. Other than the variance, it is also possible to characterize, albeit indirectly, the {\em limit distribution} of the numbers $\Lambda_{R}$ (see Theorem \ref{thm:rand mod single} below). Next, we discuss the finer aspects of the numbers $\Lambda_{R}$, such as by-products of their {\em spacing distribution} (self-correlations) and their {\em pair correlation}.

\subsubsection{Self-correlations}

To define the self-correlations between the numbers $\Lambda_{R}$ (more generally, the joint distributions between the numbers $\Lambda_{R}$),
we assume that
\begin{equation}
\label{eq:GammaR ordered list}
\Gamma_{R} = \left\{\lambda_{j}:\: 1\le j\le K \right\}
\end{equation}
is ordered in nondecreasing order of the corresponding argument\footnote{We ignore the possibility that two or more angles are equal. This could only happen for lattice points lying on the coordinate axes, and they may be ordered arbitrarily.} in the polar representation of $\lambda_{j}$, thinking of the indices in the circular sense modulo $K$ (for example, $K+1:=1$). Of our particular interest is the correlations between the numbers $\Lambda_{R}$ at high distance, i.e. $\Ac_{R}(\lambda_{j})$ and $\Ac_{R}(\lambda_{j+k})$ as $k\rightarrow\infty$. Our next result is of fundamental nature: It asserts that the limit correlation of $\Ac_{R}(\lambda_{j}), \Ac_{R}(\lambda_{j+k})$, as $R\rightarrow\infty$ exists, for every $k\ge 1$.

\begin{theorem}
\label{thm:Ck existence}
For every $k\ge 1$, the limit
\begin{equation}
\label{eq:Ck existence}
\Cc_{k}:=\lim\limits_{R\rightarrow\infty}\frac{1}{K(R)}\sum\limits_{j=1}^{K(R)} \left(\Ac_{R}(\lambda_{j})\cdot \Ac_{R}(\lambda_{j+k}) - \frac{1}{4}\right)
\end{equation}
exists.
\end{theorem}

Since, by Theorem \ref{thm:C0 lim exist}(ii), the variance of $\Lambda_{R}$ is asymptotic to a constant, the numbers $\Cc_{k}$ in \eqref{eq:Ck existence} are
equivalent to the autocorrelations of $\Lambda_{R}$.
In Section \ref{sec:limit mod corr} below, the numbers $\Cc_{k}$ are characterized in terms of a certain random model, and, in principle, it is possible to numerically evaluate $\Cc_{k}$ for each given $k\ge 1$, though it fast becomes infeasible.
Of course, Theorem \ref{thm:C0 lim exist} could be included as case $k=0$ within Theorem \ref{thm:Ck existence}, with $$\Cc_{0}=c_{0},$$ though those are of different nature, as we will see below.
The following conjecture is a possible first rigorous manifestation of the aforementioned folklore heuristic argument  in support of Hardy's conjecture \eqref{eq:Hardy conj alpha0=1/2}, and hence, by itself, is supporting \eqref{eq:Hardy conj alpha0=1/2}.

\begin{conjecture}[Vanishing Correlation Conjecture]
\label{conj:Ck->0}
As $k\rightarrow\infty$, the numbers $\Cc_{k}$ vanish, that is,
\begin{equation*}
\lim\limits_{k\rightarrow\infty}\Cc_{k}=0.
\end{equation*}
\end{conjecture}

A stronger variant of Conjecture \ref{conj:Ck->0} is presented within Section \ref{sec:rand mod}, see Conjecture \ref{conj:independence}.
Though we do not settle the Vanishing Correlation Conjecture \ref{conj:Ck->0}, in what follows we state several results, including Theorem \ref{thm:Ck->0 average} and Theorem \ref{thm:Ck->0 constrained}, in support of that conjecture. Theorem \ref{thm:Ck->0 average} asserts that $\{\Cc_{k}\}$ vanish {\em on average}, weaker than the Vanishing Correlation Conjecture.

\begin{theorem}
\label{thm:Ck->0 average}
Let $\{\Cc_{k}\}_{k\ge 1}$ be the numbers \eqref{eq:Ck existence} defined in Theorem \ref{thm:Ck existence}. Then
\begin{equation*}
\lim\limits_{L\rightarrow\infty}\frac{1}{L}\sum\limits_{k=1}^{L}\Cc_{k} = 0.
\end{equation*}

\end{theorem}

\subsubsection{Pair correlation}
In Conjecture \ref{conj:Ck->0} the key assumption that $k\rightarrow \infty$ has the effect of separating $\lambda_j$ from $\lambda_{j+k}$. Instead of imposing the separation of the points through the ordering of their angles, we will now consider correlations of $\pazocal A_R(\lambda)$ at lattice points that are well separated in terms of the distance between their angles. Our main result establishes the analogue of Conjecture \ref{conj:Ck->0} in this setting.
To state this result, we introduce some more notation. Recall that $\Gamma_{R}$ is the collection
\eqref{eq:Gamma points annulus def} of all lattice points lying in the annulus $\Dc_{R;1/\sqrt{2}}$,
and their number $K=K(R)=\#\Gamma_{R}$ is asymptotic to $\pazocal K=\pazocal K(R)=\sqrt{8}\pi R$. For $\lambda\in\Gamma_{R}$ let $\theta_{\lambda}\in [0,2\pi)$ be its argument, i.e.
\begin{equation}
\label{eq:lambda polar unshifted}
\lambda=R_{\lambda}e^{i\theta_{\lambda}},
\end{equation}
with some $R_{\lambda}\in (R-1/\sqrt{2},R+1/\sqrt{2})$, thinking of $\lambda\in\R^{2}\cong\C$. Additionally, given a positive real number $L$, $[a,b] \subseteq \mathbb R$, and $\theta \in \mathbb R$, if there exists $j \in \mathbb Z$ such that
$\theta \in [a+jL, b+jL]$ then we write $\theta \Mod L \in [a,b]$.

\begin{theorem}
\label{thm:Ck->0 constrained}
One has
\begin{equation*}
\lim\limits_{\substack{k,R\rightarrow\infty \\ k\le R^{1/2-o(1)}}}\frac{1}{\pazocal K}\sum\limits_{\substack{\lambda,\mu\in\Gamma_{R} \\ \theta_{\lambda}-\theta_{\mu} \Mod{2\pi} \in\left[ \frac{k}{\pazocal K},\frac{k+2\pi}{\pazocal K}\right]}}
\left(\Ac_{R}(\lambda)-\frac{1}{2}\right)\cdot \left(\Ac_{R}(\mu)-\frac{1}{2}\right) = 0.
\end{equation*}
\end{theorem}

Some explanation is due.
As $\lambda$ varies over lattice points in $\Gamma_R$ one may expect that asymptotically as $R \rightarrow \infty$, there is on average one point $\mu\in\Gamma_{R}$ such that
\begin{equation}
\label{eq:k<tl-tm<k+1}
\theta_{\lambda}-\theta_{\mu} \Mod{2\pi} \in [\tfrac{k}{\pazocal K}, \tfrac{k+2\pi}{\pazocal K}]
\end{equation}
provided $k \rightarrow \infty$ sufficiently slowly with $R$  (see Theorem \ref{thm:jointdistribution}), since the region $$\left\{re^{i\theta} :R-\tfrac{1}{\sqrt{2}}<r<R+\tfrac{1}{\sqrt{2}} \,\, \& \,\, (\theta_0-\theta) \Mod{2\pi}\in [\tfrac{k}{\pazocal K},\tfrac{k+2\pi}{\pazocal K}]\right\}$$ has unit area for any $\theta_0 \in\mathbb R$, $k\ge 0$. Hence, for $\lambda=\lambda_{j}$ with some $j\le K$, and $\mu$ satisfying \eqref{eq:k<tl-tm<k+1}, one may expect that $\mu=\lambda_{j+k}$, upon ordering the lattice points as in \eqref{eq:GammaR ordered list}. If this intuition would be rigorously justifiable, that would imply the Vanishing Correlation Conjecture \ref{conj:Ck->0}. However, such a rigidity statement is not expected to hold, and, given $\lambda=\lambda_j\in \Gamma_{R}$, the number of lattice points $\mu$ satisfying \eqref{eq:k<tl-tm<k+1} fluctuates between $0$ and some absolute constant (see Lemma \ref{lem:lat pnt cons rect disc} below). Therefore, rather than the vanishing of $\Cc_{k}$ as $k\rightarrow\infty$, Theorem \ref{thm:Ck->0 constrained} demonstrates that some {\em mixture} of $\Cc_{k}$ vanishes.

\subsection{Outline of the paper}
In Section \ref{sec:random-results} we introduce random models that capture both the behavior of the lattice points lying near the boundary of the circle and their correlations; we also state theorems \ref{thm:rand mod single}, \ref{thm:avgconj}, \ref{thm:mixconj}, and \ref{thm:rand mod joint}, that contain theorems \ref{thm:C0 lim exist}, \ref{thm:Ck existence}, \ref{thm:Ck->0 average}, and  \ref{thm:Ck->0 constrained} as special cases. Further discussion of the results, random model, and proofs is given in Section \ref{sec:discussion}. Theorem \ref{thm:rand mod single} is proved in Section \ref{sec:random proof1} and Theorem \ref{thm:rand mod joint} is proved in Section \ref{sec:random2}. Several auxiliary results on the distribution of lattice points in sectors, theorems \ref{thm:jointdistribution} and \ref{thm:varbd}, are stated in Section \ref{sec:lattice-cor} and the proofs of Proposition \ref{eq:prop:(r,theta) equidist} as well as theorems \ref{thm:Ck->0 average},  \ref{thm:Ck->0 constrained},  \ref{thm:avgconj}, and \ref{thm:mixconj} are given in this section. The key technical estimates used to prove Theorem \ref{thm:rand mod joint} are given Section \ref{sec:random-proofs}. Section \ref{sec:lattice-intro} states the main lattice point estimates proved in the paper, propositions \ref{thm:paircorrelation} and \ref{thm:varasymp}, then uses these estimates to prove theorems \ref{thm:jointdistribution} and \ref{thm:varbd}. Sections \ref{sec:lattice-proof1} and \ref{sec:lattice-proof2} are devoted to proving propositions \ref{thm:paircorrelation} and \ref{thm:varasymp}.

\subsection*{Acknowledgement}

The authors of this manuscript are grateful to P. Kurlberg, Z. Rudnick and M. Sodin for their interest in our work, and the stimulating discussions around it. S.L. is partially supported by EPSRC Standard Grant EP/T028343/1.

\section{Randomness and lattice points} \label{sec:random-results}

\subsection{A random model for $\Ac_{R}$}
\label{sec:rand mod}

In this section, we introduce a random model for $\Lambda_{R}$ in \eqref{eq:LambdaR def}, that explains its asymptotic distribution, as $R\rightarrow\infty$. To motivate it, it is convenient to work with the (shifted) polar coordinates: For a lattice point $\lambda\in\Gamma_{R}$ let
\begin{equation}
\label{eq:lambda polar coord}
\lambda=(R+r_{\lambda})e^{i\theta_{\lambda}},
\end{equation}
where
$r_{\lambda}\in [ -\tfrac{1}{\sqrt{2}},\tfrac{1}{\sqrt{2}} ],$ and $\theta_{\lambda}\in [0,2\pi)$
(cf. \eqref{eq:lambda polar unshifted}). The following proposition asserts that, as $\lambda$ varies in $\Gamma_{R}$, the tuples $(r_{\lambda},\theta_{\lambda})$ equidistribute in
\begin{equation}
\label{eq:G rect def}
G:=[-\tfrac{1}{\sqrt{2}},\tfrac{1}{\sqrt{2}}]\times [0,2\pi).
\end{equation}

\begin{proposition}
\label{eq:prop:(r,theta) equidist}
For $\lambda\in\Gamma_{R}$ let $(r_{\lambda},\theta_{\lambda}) \in G$ be its polar coordinates \eqref{eq:lambda polar coord},
with $G$ the rectangle
\eqref{eq:G rect def}.
Then as $R\rightarrow\infty$ the random vectors $(r_{\lambda},\theta_{\lambda})\in G$ w.r.t. uniformly drawn $\lambda\in\Gamma_{R}$,
converge in distribution to the random uniform vector
$(r,\theta)\in G$, that is, the vector corresponding to the probability measure $\frac{dr d\theta}{\sqrt{8}\pi}$ on $G$.
\end{proposition}

Recall that $S=[-1/2,1/2]^{2}$ is the unit square, and $S(x)$ is its shift by $x\in\R^{2}$.
Let\ $\lambda = (R+r_{\lambda})\cdot e^{i\theta_{\lambda}}\in \Gamma_{R}$ be a lattice point and consider the square $S(\lambda)$ centered at $\lambda$. Then, thinking of the portion of the circle $\partial B(R)$ as ``approximately linear", we may rotate the emerging picture by $\frac{\pi}{2}-\theta$, so that the relevant sector of the ``straight line" $\partial B(R)$ will lie on the vertical axis, whereas $S$ is tilted by $(\pi/2-\theta)$ and shifted by $r$, see Figure \ref{fig:sect vs straight} left. Since $(r_{\lambda},\theta_{\lambda})$ asymptotically equidistribute on $G$ by Proposition \ref{eq:prop:(r,theta) equidist}, the following model is only natural for the distribution of $\Lambda_{R}$.

For $(r,\theta)\in G$ denote $S_{\theta}=T_{\theta}S$ the $\theta$-tilt of $S$ i.e. $T_{\theta} : \mathbb R^2 \rightarrow \mathbb R^2$ is the linear transformation that corresponds to $\begin{pmatrix}
 \cos \theta & -\sin \theta \\ \sin \theta & \cos \theta \end{pmatrix}$, and let
$$S_{r,\theta} := (r,0)+S_{\theta}$$ be the horizontal shift by $r$ of $S_{\theta}$.
Denote $\Ac_{\infty}$ to be the area of the portion of $S_{r,\theta}$ lying to the left of the $y$ axis, that is, $\Ac_{\infty}$ is the random variable
\begin{equation}
\label{eq:Ac inf def}
\Ac_{\infty}=\Ac_{\infty}(r,\theta) = \area(S_{r,\theta}\cap (\R_{\le 0}\times\R)) = \area\{ (x,y)\in S_{r,\theta}:\: x\le 0 \},
\end{equation}
with $(r,\theta)$ random, uniformly distributed in $G$. Here the $y$ axis models the circle $\partial B(R)$ separating the interior of $B(R)$,
modelled by the negative half-plane $\R_{<0}\times \R$, from its exterior, modelled by the positive half-plane $\R_{>0}\times \R$.

\begin{figure}[h]
  \centering
\includegraphics[height=7cm]{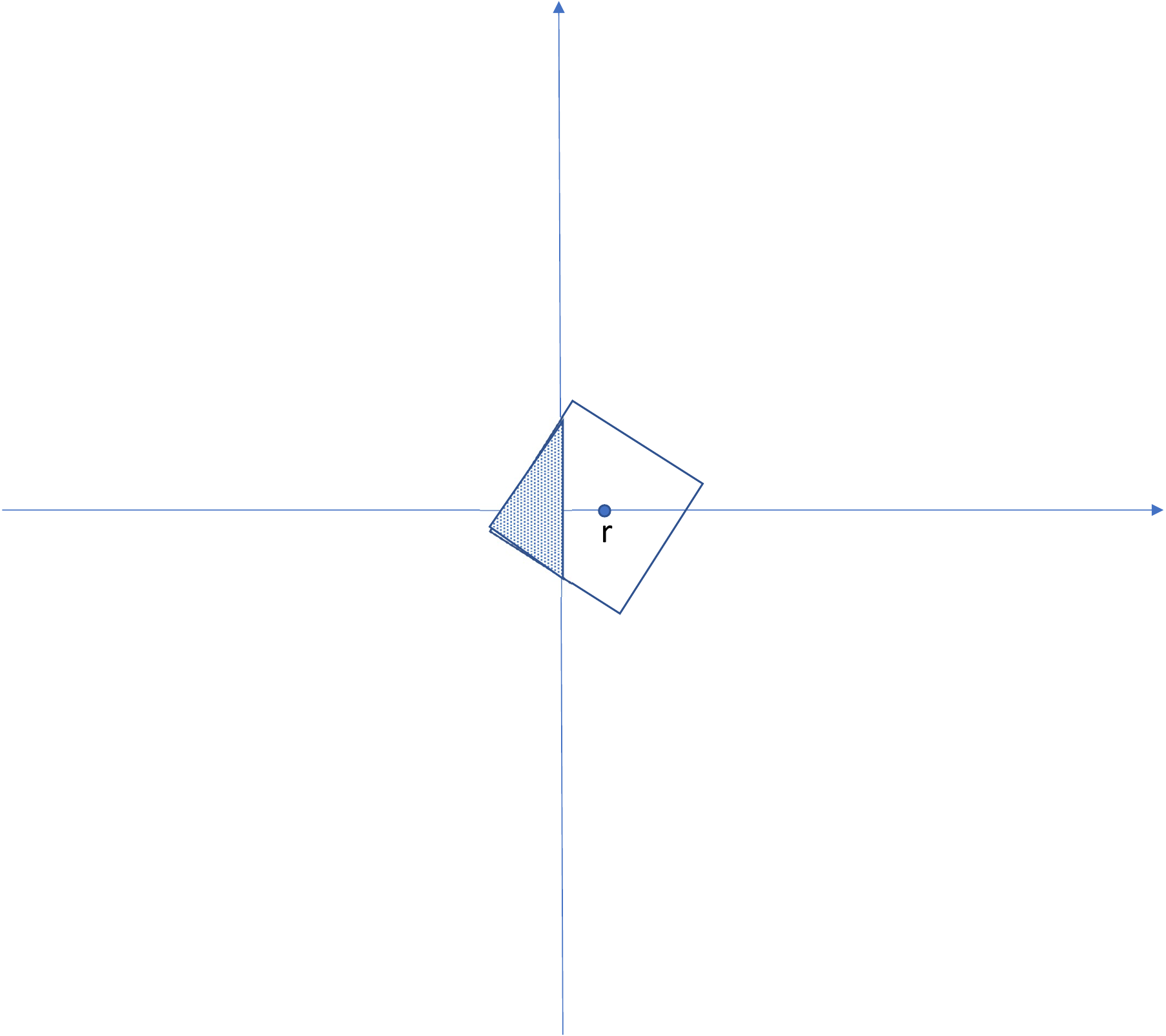}
\includegraphics[height=7cm]{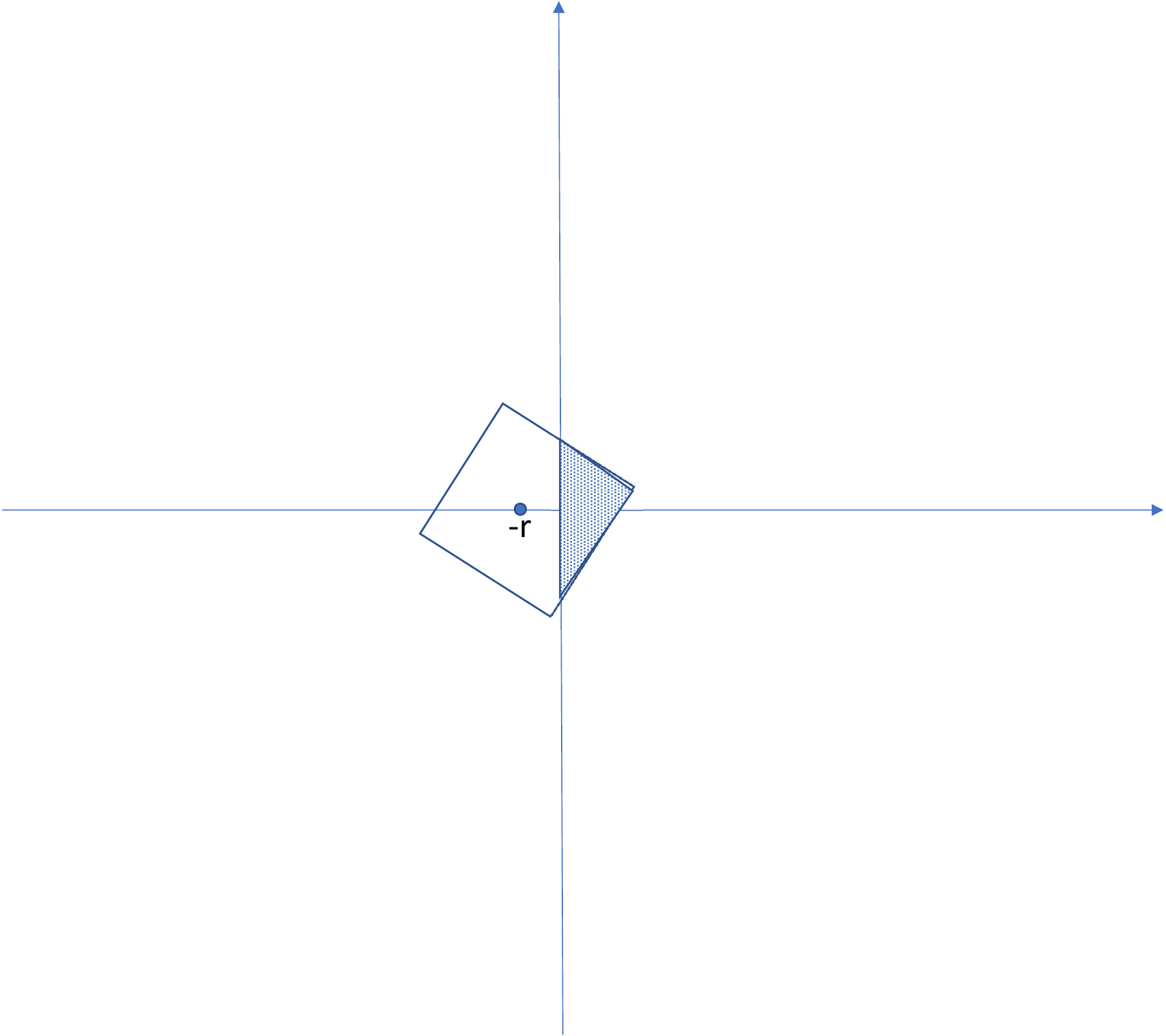}
\caption{$\Ac_{\infty}(r,\theta)$ is the area of the portion of the square to the left of $y$ axis: $\Ac_{\infty}(r,\theta)$ (left), $\Ac_{\infty}(-r,\theta)$ (right). The blue shaded areas are equal.}
  \label{fig:A(-r,t)=1-A(r,t)}
\end{figure}

Since the distribution of $\Ac_{\infty}$ is invariant w.r.t. $t\mapsto 1-t$ (see Figure \ref{fig:A(-r,t)=1-A(r,t)}
and Lemma \ref{lem:Ainf explicit}(ii)), it follows, in particular, that
\begin{equation}
\label{eq:exp inf = 1/2}
\E\left[  \Ac_{\infty} \right] = \frac{1}{2} ,
\end{equation}
consistent with \eqref{eq:expectation individual 1/2}.
The following theorem asserts that $\Lambda_{R}$ (considered as a random variable that draws a random uniform number in $\Lambda_{R}$), converges in distribution to $\Ac_{\infty}$.

\begin{theorem}
\label{thm:rand mod single}

\begin{enumerate}[i)]

\item The numbers $\Lambda_{R}$,  converge, in distribution,
to the random variable $\Ac_{\infty}$.

\item One has
\begin{equation*}
\var\left( \Ac_{\infty} \right) = c_{0},
\end{equation*}
with $c_{0}$ given by \eqref{eq:c0 asymp var def}.

\end{enumerate}

\end{theorem}

\begin{proof}[Deducing Theorem \ref{thm:C0 lim exist}(ii) from Theorem \ref{thm:rand mod single}]

Since $|\Ac_{R}(\cdot)|$ is bounded (by $1$), it follows ~\cite[Theorem 3.5]{Billingsley}
that {\em all} the moments of $\Ac_{R}$ are asymptotic to the moments
of $\Ac_{\infty}$, by Theorem \ref{thm:rand mod single}(i). Together with Theorem \ref{thm:rand mod single}(ii),
and bearing in mind \eqref{eq:exp inf = 1/2} and Theorem \ref{thm:C0 lim exist}(i), Theorem \ref{thm:C0 lim exist}(ii) follows.
\end{proof}

Theorem \ref{thm:rand mod single} allows for stating the following strengthened version of the Vanishing Correlation Conjecture \ref{conj:Ck->0}.
Mind that, by the virtue of Theorem \ref{thm:rand mod single}(i), for every $k\ge 1$, both the marginal distributions of the vector $$\{(\Ac_{R}(\lambda_{j}),\Ac_{R}(\lambda_{j+k}))\}_{1\le j\le K}$$ converge to $\Ac_{\infty}$.

\begin{conjecture}
\label{conj:independence}
The random vectors $$\{(\Ac_{R}(\lambda_{j}),\Ac_{R}(\lambda_{j+k}))\}_{1\le j\le K}$$ converge, in distribution 
in the double limit $R\rightarrow\infty$ then $k\rightarrow\infty$, to the random vector
\begin{equation}
\label{eq:Y1 Y2 theta diag}
(Y_{1},Y_{2}):=(\Ac_{\infty}(r,\theta),\Ac_{\infty}(r',\theta))
\end{equation}
where $(r,\theta,r')$ is random uniform in
$$G\times \left[-\frac{1}{\sqrt{2}},\frac{1}{\sqrt{2}} \right] = \left[-\frac{1}{\sqrt{2}},\frac{1}{\sqrt{2}} \right] \times [0,2\pi)\times \left[-\frac{1}{\sqrt{2}},\frac{1}{\sqrt{2}} \right].$$

\end{conjecture}

The random variables $Y_{1},Y_{2}$ in Conjecture \ref{conj:independence} are uncorrelated, but {\em not} independent. Further, conditioned on $\theta$, $Y_{1},Y_{2}$ are independent, see the discussion in section \ref{sec:indep vs uncorr}. In support of this conjecture we have the following results, Theorem \ref{thm:avgconj} and Theorem \ref{thm:mixconj}.

\begin{theorem} \label{thm:avgconj}
Let $M$ tend to infinity with $R$ in such a way
so that $M \le R^{1/10-o(1)}$. Also, let $L>0$ satisfy $M^{1/2+o(1)} \le L \le M$.
We have that \[\left\{\left(\pazocal A_R(\lambda_j),
\pazocal A_R(\lambda_{j+k})\right)\right\}_{1 \le j \le K;\, M \le k \le M+L}\] converges, in distribution as $R \rightarrow \infty$, to the random vector
$(Y_{1},Y_{2})$ of \eqref{eq:Y1 Y2 theta diag}.
\end{theorem}

\begin{theorem} \label{thm:mixconj}
Let $k$ tend to infinity with $R$ in such a way
so that $k \le R^{1/2-o(1)}$.
As $(\lambda,\mu)$ varies over pairs of lattice points
in $\Gamma_R\times \Gamma_R$ with $$\theta_{\lambda}-\theta_{\mu} \Mod{2\pi} \in\left[ \frac{k}{\pazocal K},\frac{k+2\pi}{\pazocal K}\right],$$ we have that
$(\pazocal A_R(\lambda),
\pazocal A_R(\mu))$ converges, in distribution as $R \rightarrow \infty$, to the random vector
$(Y_{1},Y_{2})$ of \eqref{eq:Y1 Y2 theta diag}.
\end{theorem}

\subsection{A random model for the correlations $\Cc_{k}$}

\label{sec:limit mod corr}

For $(r,\theta)\in G$ let $\widetilde{\Rc}$ be the (semi-infinite) rectangle
\begin{equation}
\label{eq:R semi-inf rectangle}
\widetilde{\Rc}:=[0,+\infty)\times [-\tfrac{1}{\sqrt{2}},\tfrac{1}{\sqrt{2}}]
\end{equation}
and $$\widetilde{\Rc}_{\theta}= T_{\theta}\widetilde{\Rc}$$ its tilt
by $\theta$.
Now let
\begin{equation}
\label{eq:Rc tilted shifted}
\Rc=\Rc_{r,\theta}=\widetilde{\Rc}_{\theta}+ r\cdot e^{i(\theta+\pi/2)}
\end{equation}
be the same (tilted) rectangle shifted by $r$ along the short side of $\widetilde{\Rc}_{\theta}$.
Denote the resulting linear map $$Q=Q_{r,\theta}:x\mapsto T_{\theta}x+  r\cdot e^{i(\theta+\pi/2)},$$ so that $Q\widetilde{\Rc}=\Rc$.
Equivalently, we may obtain $\Rc_{r,\theta}$ from $\widetilde{\Rc}$ by first shifting by $r$ (positive or negative) in the vertical direction so that to put the origin at $(0,r)\in \R^{2}$, and then tilting by $\theta$ via the coordinate axes, also shifted by $(0,r)$.

\begin{figure}[h]
  \centering
\includegraphics[height=6cm]{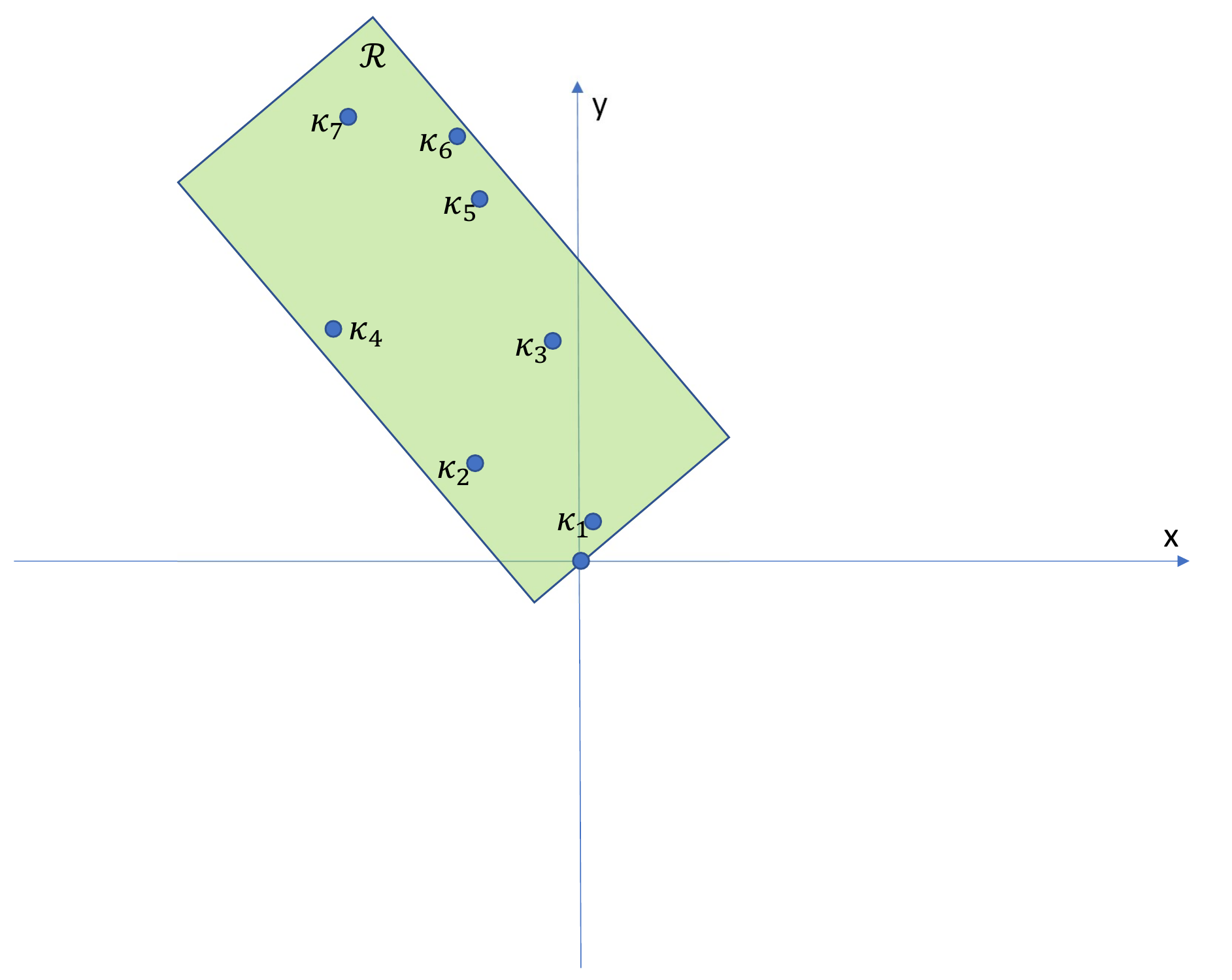}
\includegraphics[height=6cm]{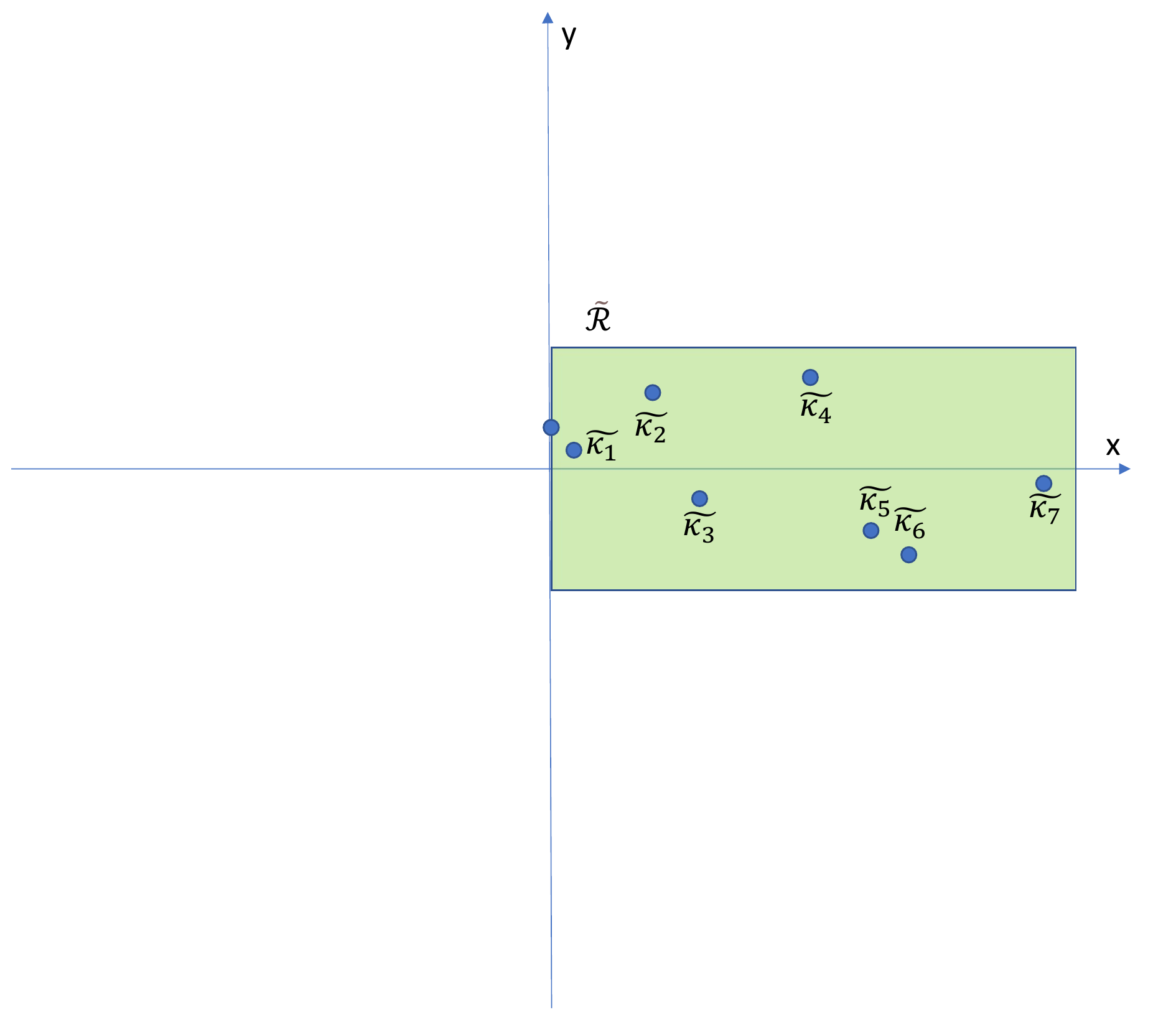}
\caption{Left: A finite part of the semi-infinite rectangle $\Rc$ is depicted, containing $7$ lattice points. Right: its image onto $\widetilde{\Rc}$. }
\label{fig:pointsR}
\end{figure}

Since, by the above procedure, $$\Rc=T_{\theta}\widetilde{\Rc}+r\cdot e^{i(\theta+\pi/2)},$$ we may invert to recover $\widetilde{\Rc}$: $$\widetilde{\Rc}=T_{-\theta}(\Rc-r\cdot e^{i(\theta+\pi/2)}).$$ Denote the resulting map $P=P_{r,\theta}=Q^{-1}$, $$x\mapsto Px=T_{-\theta}(x-r\cdot e^{i(\theta+\pi/2)}),$$ and let
\begin{equation}
\label{eq:lambdal in R}
\kappa_{1}=\kappa_{1}(r,\theta),\ldots, \kappa_{\ell}=\kappa_{\ell}(r,\theta),\ldots\in \Rc
\end{equation}
be the sequence of lattice points of $\Rc$, not including the origin, ordered by nondecreasing
projection onto the long side of $\Rc$.
We denote the sequence of the images of the $\kappa_{\ell}$ in $\widetilde{\Rc}$:
\begin{equation}
\label{eq:lambdaj=Ptild}
\{\widetilde{\kappa}_{1}=P\kappa_{1},\ldots \widetilde{\kappa}_{\ell}=P\kappa_{\ell},\ldots\}\subseteq \widetilde{\Rc},
\end{equation}
illustrated in Figure \ref{fig:pointsR}, left.
For $\ell\ge 1$ we write
\begin{equation}
\label{eq:lambdaj=(tj,rj)}
\widetilde{\kappa}_{\ell}=(t_{\ell},\rho_{\ell})\in \widetilde{\Rc}
\end{equation}
in {\em Cartesian coordinates} with
\begin{equation}
\label{eq:rho proj def}
\rho_{\ell}= \rho_{\ell}(r,\theta)\in [-\tfrac{1}{\sqrt{2}},\tfrac{1}{\sqrt{2}}]
\end{equation}
the projection of $\kappa_{\ell}$ onto the short side of $\Rc$, and $t_{\ell}\in (0,+\infty]$ is the projection of $\kappa_{\ell}$ on the long side of $\Rc$, illustrated in Figure \ref{fig:pointsR}, right.
(Equivalently, $\rho_{\ell}$ and $t_{\ell}$ are the projections of $\widetilde{\kappa}_{\ell}$ onto the short or the long side of $\widetilde{\Rc}$ respectively.) Since the $\kappa_{\ell}$ were ordered according to nondecreasing projection onto the long side of $\Rc$, the corresponding $\widetilde{\kappa}_{\ell}$ are {\em sorted} by nondecreasing order of $t_{\ell}>0$.

Given $(r,\theta)$ corresponding to $\lambda_{j}\in\Gamma_{R}$, and a correlation distance $k\ge 1$, the quantity which we associate to $\Ac_{R}(\lambda_{j+k})$ in our model is $\rho_{k}=\rho_{k}(r,\theta)$, that is used to define its by-product $\Ac_{\infty}(\rho_{k},\theta)$
with $\Ac_{\infty}(\cdot,\cdot)$ as in \eqref{eq:Ac inf def}.
Definition \ref{def:Ac inf def} will explicate this and will also introduce the {\em model correlations}.
We decompose
$$\widetilde{\Rc}=([0,\infty)\times [-\tfrac{1}{\sqrt{2}},0])\cup  ([0,\infty)\times [0,\tfrac{1}{\sqrt{2}}]) =:\widetilde{\Rc}_{1}\cup  \widetilde{\Rc}_{2} , $$
and their images
\begin{equation*}
\Rc = Q\widetilde{\Rc_{1}}\cup Q\widetilde{\Rc_{2}} =:\Rc_{1}\cup \Rc_{2},
\end{equation*}
where $\Rc_{1}=\Rc_{1;\theta}$ and $\Rc_{2}=\Rc_{2;\theta}$.

\begin{definition}
\label{def:Ac inf def}

\begin{enumerate}[i)]

\item
For $k\ge 1$, and $(r,\theta)\in G$ let
\begin{equation}
\label{eq:Ac inf k def}
\Ac_{\infty,k}(r,\theta) = \area\left( (S + Q_{r,\theta}(t,\rho_{k}))\cap \Rc_{1}  \right),
\end{equation}
where $(t,\rho_{k})\in \widetilde{\Rc}$, and $t>0$ is chosen arbitrary sufficiently large so that $S + Q(t,\rho_{k})$ does not intersect the short side of $\Rc$.

\item
For a number $k\ge 1$ we denote the covariance, w.r.t. $(r,\theta)\in G$ random uniform, of $\Ac_{\infty}(r,\theta)$ and $\Ac_{\infty}(\rho_{k},\theta)=\Ac_{\infty,k}(r,\theta)$, i.e.
\begin{equation}
\label{eq:Cktild def}
\widetilde{\Cc}_{k}:= \E_{r,\theta}\left[\Ac_{\infty}(r,\theta)\cdot \Ac_{\infty,k}(r,\theta)\right] - \frac{1}{4}=
\cov_{r,\theta}\left(\Ac_{\infty}(r,\theta),\Ac_{\infty,k}(r,\theta)\right),
\end{equation}
with $\Ac_{\infty,k}(r,\theta)$ as in \eqref{eq:Ac inf k def}.

\end{enumerate}
\end{definition}

\vspace{2mm}

\begin{figure}[h]
  \centering
\includegraphics[height=10cm]{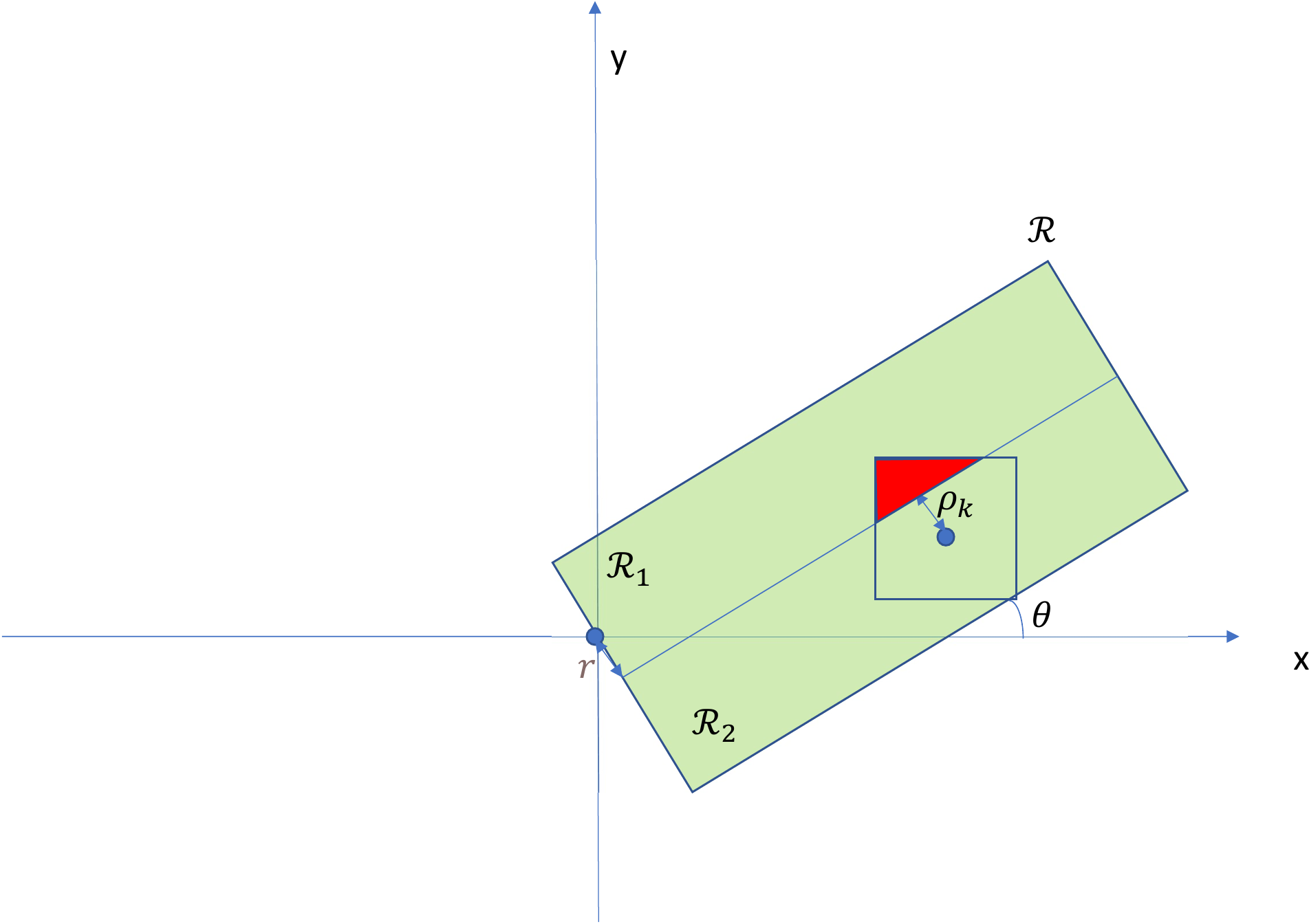}
\caption{Definition \ref{def:Ac inf def} illustrated. Here $\Ac_{\infty,k}(r,\theta)$ is designated by the color red. }
\label{fig:corr illustration}
\end{figure}

Definition \ref{def:Ac inf def}(i), with \eqref{eq:Ac inf k def} is clearly equivalent to
\begin{equation}
\label{eq:Ainfk,r=Ainfrk}
\Ac_{\infty,k}(r,\theta)=\Ac_{\infty}(\rho_{k},\theta),
\end{equation}
and is illustrated in Figure \ref{fig:corr illustration}.
Mind that $S + Q(t,\rho)$ not intersecting the short side of $\Rc$ in Definition \ref{def:Ac inf def}(i) is equivalent to $PS+(t,\rho)$ not
intersecting the short side $\{0\}\times [-1/\sqrt{2},1/\sqrt{2}]$ of $\widetilde{\Rc}$.
One may naturally extend the definition of $\Ac_{\infty,k}(r,\theta)$ for $k=0$ by
replacing $\rho_{k}$ with $r$, in which case this notion coincides with $\Ac_{\infty}$ as in \eqref{eq:Ac inf def}, and the equality $$\Ac_{\infty,k}(r,\theta)=\Ac_{\infty}(\rho_{k},\theta)   $$ holds true for $k\ge 0$.
The principal result of the section is the following theorem, evidently implying Theorem \ref{thm:Ck existence}:

\begin{theorem}
\label{thm:rand mod joint}
For every $k\ge 1$, the limit
\begin{equation*}
\Cc_{k}:=\lim\limits_{R\rightarrow\infty}\frac{1}{K(R)}\sum\limits_{j=1}^{K(R)} \left(\Ac_{R}(\lambda_{j})\cdot \Ac_{R}(\lambda_{j+k}) - \frac{1}{4}\right),
\end{equation*}
exists and coincides with $\widetilde{\Cc}_{k}$ as in \eqref{eq:Cktild def}.
\end{theorem}

\section{Discussion} \label{sec:discussion}

\subsection{On the proofs of the main results}

Given $\lambda=(R+r_{\lambda})e^{i\theta_{\lambda}}$ as in \eqref{eq:lambda polar coord}, it will be shown that as $R\rightarrow\infty$ the value of $\pazocal A_R(\lambda)$ depends, up to lower order terms, only on $(r_{\lambda},\theta_{\lambda})$ and not $R$, see \eqref{eq:AR=Ainf+O(1/R)}. Since the leading term in the r.h.s. of \eqref{eq:AR=Ainf+O(1/R)} is a continuous functional $\Ac_{\infty}(r,\theta)$ as in \eqref{eq:Ac inf def},
it is possible to infer the distribution of $\{\Ac_{\infty}(r_{\lambda},\theta_{\lambda})\}_{\lambda\in\Gamma_{R}}$ from the equidistribution of Proposition \ref{eq:prop:(r,theta) equidist} of the shifted polar coordinates $$\{(r_{\lambda},\theta_{\lambda})\}_{\lambda\in\Gamma_{R}}\subseteq G$$ of lattice points, via an application of the continuous mapping theorem. This approach will yield the asymptotic distribution part of Theorem \ref{thm:rand mod single}. The variance part of Theorem \ref{thm:rand mod single} is done by a direct and explicit, alas somewhat long and tedious computation of $\Ac_{\infty}(r,\theta)$, and its variance.

\vspace{1mm}

To prove Theorem \ref{thm:rand mod joint} one aims to rigorously justify the limit model of Section \ref{sec:limit mod corr}.
In analogy to the above, the correlations $$\sum\limits_{j=1}^K\pazocal A_R(\lambda_j)\pazocal A_R(\lambda_{j+k})$$ may also be understood in terms of the joint distribution of the tuples $(r_{j},\theta_{j})$ and $(r_{j+k},\theta_{j+k})$ with $j \in [1,\ldots, K]$. Rather, given $k\ge 1$ and $(r_{j},\theta_{j})$, one has $\theta_{j}\approx \theta_{j+k}$, and the random model suggests that $r_{j+k}$ is almost fully determined by $(r_{j},\theta_{j})$ essentially via
\begin{equation}
\label{eq:rj+k=rk(rj)}
r_{j+k}=\rho_{k}(r_{j},\theta_{j})+O_k\bigg( \frac{1}{R}\bigg),
\end{equation}
which holds for a subset of $\{ j : 1\le j \le K\}$ of full asymptotic density.
The function $\rho_{k}(\cdot,\cdot)$ is somewhat implicit, and difficult to use for practical purposes, though perhaps it could be employed for numerical computations of the values of $\pazocal C_{k}$. However, for the {\em existence} of the limit of Theorem \ref{thm:rand mod joint}, only the continuity of $\rho_{k}(\cdot,\cdot)$ outside a measure zero sets of values of $(r,\theta)\in G$ will be needed in order to apply the continuous mapping theorem again.

A number of obstacles exist to justify the above. First, the lattice points in the annulus $\Dc_{R;1/\sqrt{2}}$ next to a given $\lambda_{j}$ might, in general, fail to be the same as those lying in the approximating rectangle. Second, even if these two sets coincide, the ordering, w.r.t. the polar angle in the annulus might differ from the ordering and w.r.t. the projection on the long side of the approximating rectangle. If either of these two problematic scenarios occur, then one may no longer infer the important approximate identity \eqref{eq:rj+k=rk(rj)}. Proposition \ref{prop:rect pert arc lat pnts} and its proof shows that either of these two obstacles is rare: Discrepancies between the said two sets would imply a lattice point in a narrow sector, whereas the ordering violation would imply a very restrictive condition on the angle of $\lambda_{j}$ or one of its neighbors.

\vspace{2mm}

The strategy used to prove theorems \ref{thm:avgconj} and \ref{thm:mixconj} is similar, and
we will only describe the ideas used in the proof of Theorem \ref{thm:mixconj}.
Let $$\pazocal B=\left\{ (\lambda, \mu) \in \Gamma_R \times \Gamma_R : \theta_{\lambda}-\theta_{\mu} \Mod{2\pi} \in \left[ \frac{k}{\pazocal K},\frac{k+2\pi}{\pazocal K}\right]\right\},$$
where $\pazocal K= \pazocal K(R)=\sqrt{8} \pi R$.
By our earlier observation, to prove Theorem \ref{thm:mixconj} it suffices to compute the joint distribution of
  $(r_{\lambda}, \theta_{\lambda}, r_{\mu})_{(\lambda, \mu)\in \pazocal B}$, that is, we wish to asymptotically estimate
\begin{equation} \label{eq:dist-sum}
\frac{1}{\pazocal K}\sum_{\substack{(\lambda, \mu) \in \pazocal B  \\ (r_{\lambda}, \theta_{\lambda}, r_{\mu}) \in (I_1,J,I_2) }} 1
\end{equation}
as $R \rightarrow \infty$ where $I_1,I_2 \subseteq[-\tfrac{1}{\sqrt{2}},\tfrac{1}{\sqrt{2}}]$ and $J \subseteq[0,2\pi)$ are intervals.
To analyze the above sum, we expand the condition $\theta_{\lambda}-\theta_{\mu} \Mod{2\pi} \in [ \frac{k}{\pazocal K},\frac{k+2\pi}{\pazocal K}]$ as, roughly, a trigonometric polynomial of degree $\pazocal K$, i.e.
\[
\frac{1}{\pazocal K}\sum_{\substack{(\lambda, \mu) \in \pazocal B \\ (r_{\lambda}, \theta_{\lambda}, r_{\mu}) \in (I_1,J,I_2) } } 1 \approx  \frac{1}{\pazocal K^2}\sum_{|\ell| \le \pazocal K}  e^{-2\pi i \ell k} \sum_{\substack{\mu, \lambda \in  \Gamma_R \\ (r_{\lambda}, \theta_{\lambda}, r_{\mu}) \in (I_1,J,I_2)}}  e^{i\ell (\theta_{\lambda}-\theta_{\mu})} .
\]
We then transform the inner sum on the r.h.s. using Poisson summation (to make this rigorous we use a smooth function that approximates the condition  $(r_{\lambda}, \theta_{\lambda}, r_{\mu}) \in (I_1,J,I_2)$) and collect a main term from $\ell=0$, which corresponds to the equidistribution of $(r_{\lambda}, \theta_{\lambda}, r_{\mu})$. In our range of interest the dual sum after applying Poisson summation will be shorter than the original sum, providing a gain at this step. However if one bounds the $\ell \neq 0$  terms individually this gives a worse than trivial estimate, since our sum over $\ell$ is long and even the optimal conjectural bounds for the individual sums would not suffice here.

We proceed to average the dual sums over $ |\ell| \le \pazocal K$ and this procedure yields a secondary main term. For fixed $k$ the aforementioned secondary term is the same size as the main term, which reflects the rigidity of neighboring lattice points, whereas for $k \rightarrow \infty$ the analysis is more subtle. We express the secondary main term as an oscillatory integral with a phase function that includes the parameter $k$ as a linear factor. Applying a stationary phase estimate, we conclude that the secondary main term decays relative to our main term as $k \rightarrow \infty$. This establishes the equidistribution of $(r_{\lambda}, \theta_{\lambda}, r_{\mu})_{(\lambda, \mu)\in \pazocal B}$ from which Theorem \ref{thm:mixconj} follows.


\subsection{Second order results for the Gauss circle problem}
Gauss' classical argument for estimating lattice points inside a circle shows that a bound for the remainder term in the asymptotic
expression \eqref{eq:expectation individual 1/2} for the expectation
\begin{equation*}
\frac{1}{K(R)}\sum\limits_{\lambda\in \Gamma_{R}} \Ac_{R}(\lambda)
\end{equation*}
and their number \eqref{eq:K(R) number points def} give a bound for the error term $\Delta(R)$ in the Gauss circle problem.
Therefore, the finer aspects of the distribution of the numbers $\{\Ac_{R}(\lambda)\}$ are not expected to
improve the known bounds for $\Delta(R)$, and the assertions presented within this manuscript are {\em second-order} results giving an important insight on the Gauss circle problem.

\subsection{Equidistribution of lattice points in Diophantine rectangles}

Though not explicitly stated, the essence of the vanishing {\em on average} of the $\Cc_{k}$ in Theorem \ref{thm:Ck->0 average} has to do with the following equidistribution fact, encapsulated within the proof of Theorem \ref{thm:avgconj}. Let $\theta\in [0,2\pi)$ be an angle, $r\in [-\tfrac{1}{\sqrt{2}},\tfrac{1}{\sqrt{2}} ]$ and consider the tilted and shifted semi-infinite rectangle $\Rc=\Rc_{r,\theta}$ as in \eqref{eq:Rc tilted shifted}. Then, for {\em ``generic"} $\theta$, the projection of $\{\kappa_{\ell}\}$, the lattice points \eqref{eq:lambdal in R} lying in $\Z^{2}\cap \Rc$, ordered w.r.t. the projection onto the long side of $\Rc$, onto the short side of $\Rc$, are equidistributed in $I:=[-\tfrac{1}{\sqrt{2}},\tfrac{1}{\sqrt{2}} ]$, {\em on average} w.r.t. $\ell$. That is, the set of the projections of $\{\kappa_{\ell}\}_{1\le \ell\le L}$ becomes equidistributed in $I$ as $L\rightarrow\infty$, for a set of $\theta$ of almost full measure (depending on $L$) inside $[0,2\pi)$.

The generic $\theta$ above are those so that $\tan\theta$ is far from a rational number with small denominator, or, what is equivalent, the direction $(\cos{\theta},\sin{\theta})$ is far from rational direction. Otherwise, if $(\cos\theta,\sin\theta)$ is rational, then $\Z^{2}\cap \Rc$ is periodic, so that their projections only the short side are far from equidistributed; in some extreme cases, such as, for example $\theta=0$ or $\theta=\pi/4$, the distribution of the projections is dominated by a single atom at $-r$. The major challenge in proving Conjecture \ref{conj:Ck->0} is showing the same equidistribution, without averaging w.r.t. $\ell$, and exploiting the averaging w.r.t. $r$ (which is not used in the aforementioned argument).

\subsection{Random tilts vs. random lattices}

Taking the fixed lattice $\Z^{2}$ and intersecting it with a randomly tilted semi-infinite rectangle is equivalent to intersecting the fixed semi-infinite rectangle with nonstandard, randomly tilted lattices. This construction is reminiscent of the random lattice problems that could be addressed using techniques involving ergodicity and the equidistibution of the given lattices w.r.t. the Haar measure defined on the space of lattices, see e.g. the recent survey ~\cite{MarklofRandLat}. Unfortunately, the (random) transformations applied on the lattice $\Z^{2}$, namely the random tilts (and the random shifts) are unlikely to equidistribute in the space of lattices w.r.t. the Haar measure, since the stretch transformation is not present within the allowed repertoire. Hence these ergodic methods are unlikely to be directly applicable in our problem, and the new techniques developed within this manuscript are necessary.

\subsection{Independence and uncorrelatedness at large distance}
\label{sec:indep vs uncorr}

Since the $\Cc_{k}$, the output of Theorem \ref{thm:Ck existence}, is an outcome of a limit, as $R\rightarrow\infty$, of the correlations between $\Ac_{R}(\lambda_{j})$ and $\Ac_{R}(\lambda_{j+k})$ with $k$ {\em fixed}, that implies that the angles $\theta_{j}$ and $\theta_{j+k}$ (polar representation of $\lambda_{j}$ and $\lambda_{j+k}$ respectively) are asymptotic, and cannot possibly be independent. However, as a by-product of our analysis presented within the proofs of Theorem \ref{thm:avgconj} and Theorem \ref{thm:mixconj}, it follows that the $(r_{j},r_{j+k})$ components of the shifted polar representation \eqref{eq:lambda polar coord} of $\lambda_{j}$ and $\lambda_{j+k}$ respectively, are still asymptotically independent, on average w.r.t. $k$, leading to Conjecture \ref{conj:independence}, with $(Y_{1},Y_{2})$ as in
\eqref{eq:Y1 Y2 theta diag}, rather than two {\em independent} copies of $\Ac_{\infty}(r,\theta)$.

Using the same reasoning as above, it makes sense to expect that as $R,k\rightarrow\infty$ {\em simultaneously}, with $k=t\cdot R$ for some $t\in [0,\sqrt{8}\cdot \pi]$ fixed, the vector $(\Ac_{R}(\lambda_{j}),\Ac_{R}(\lambda_{j+k}))$ converges in distribution to the limit random vector
$$\left(\Ac_{\infty}(r,\theta),\Ac_{\infty}\left(r,\theta+t/(\sqrt{2}\pi)\right)\right).$$ On the other hand, since the distribution of $\Ac_{\infty}(r,\theta)$, conditioned on $\theta$, genuinely depends on the value of $\theta$, the random variables $(Y_{1},Y_{2})$ are not independent, and their being uncorrelated has to do with the symmetric distribution conditioned on $\theta$, for {\em every} $\theta$.

\vspace{2mm}

The setting of Conjecture \ref{conj:independence} (or Conjecture \ref{conj:Ck->0}), where no averaging w.r.t. $k$ takes place, is very different from the averaging regime of Theorem \ref{thm:Ck->0 average}. Here, given $\theta\in [0,2\pi)$, the $\rho_{k}(r,\theta)$ is {\em fully dependent} on $r$, and, a forteriori, so does $\Ac_{\infty,k}(r,\theta)$ fully depend on $\Ac_{\infty}(r,\theta)$. For some small values of $k$ it should be possible to analytically evaluate $\rho_{k}(r,\theta)$ though it quickly becomes unfeasible. Since, as $k$ grows, the complexity of the map $r\mapsto \rho_{k}(r,\theta)$ seems
to grow, with a growing number of different ranges, it seems plausible that the map $\rho_{k}$ exhibits statistical pseudorandomness, giving some confidence in Conjecture \ref{conj:independence}.

\section{Limit distribution: Proof of Theorem \ref{thm:rand mod single}} \label{sec:random proof1}

First, we will require the following lemma, where, in particular, an explicit expression for $\Ac_{\infty}(\cdot, \cdot)$ is derived.

\begin{lemma}
\label{lem:Ainf explicit}
Let $\Ac_{\infty}$ be as in \eqref{eq:Ac inf def}.

\begin{enumerate}[i)]

\item
For every $\theta\in [0,2\pi)$, one has the equality
\begin{equation}
\label{eq:cond exp Ainf=1/2}
\frac{1}{\sqrt{2}}\int\limits_{-1/\sqrt{2}}^{-1/\sqrt{2}} \Ac_{\infty}(r,\theta)dr = \frac{1}{2}.
\end{equation}

\item For $\theta\in [0,\pi/2]$ and $r\in [-\tfrac{1}{\sqrt{2}},\tfrac{1}{\sqrt{2}}   ]$ one has the equality
\begin{equation}
\label{eq:Ainf(r,theta)=Ainf(-r,pi/2-theta)}
\Ac_{\infty}(r,\theta)=\Ac_{\infty}(r,\pi/2-\theta)= 1-\Ac_{\infty}(-r,\theta).
\end{equation}

\item
For $\theta\in (0,\pi/4)$ and $r\in [-\tfrac{1}{\sqrt{2}},\tfrac{1}{\sqrt{2}}   ]$, one has
\begin{equation}
\label{eq:Ainf(r,theta) explicit}
\Ac_{\infty}(r,\theta) = \begin{cases}
1 & -\frac{1}{\sqrt{2}}\le r<-\frac{\sin\left(\theta+\frac{\pi}{4}\right)}{\sqrt{2}}, \\
1-  \frac{\left(r+\frac{\sin\left(\frac{\pi}{4}+\theta  \right)}{\sqrt{2}}   \right)^{2}}{\sin(2\theta)}
&   -\frac{\sin\left(\theta+\frac{\pi}{4}\right)}{\sqrt{2}} \le r< -\frac{\cos\left(\theta+\frac{\pi}{4}\right)}{\sqrt{2}}, \\
\frac{\tan{\theta}}{2}-\frac{r- \frac{\cos\left( \frac{\pi}{4}+\theta  \right)}{\sqrt{2}}}{\cos{\theta}}
& -\frac{\cos\left(\theta+\frac{\pi}{4}\right)}{\sqrt{2}} \le r < \frac{\cos\left(\theta+\frac{\pi}{4}\right)}{\sqrt{2}} ,\\
\frac{1}{\sin(2\theta)} \left( \frac{\sin\left( \frac{\pi}{4}+\theta \right)}{\sqrt{2}}-r \right)^{2}
& \frac{\cos\left(\theta+\frac{\pi}{4}\right)}{\sqrt{2}} \le r < \frac{\sin\left(\theta+\frac{\pi}{4}\right)}{\sqrt{2}}, \\
0   &\frac{\sin\left(\theta+\frac{\pi}{4}\right)}{\sqrt{2}} \le r \le \frac{1}{\sqrt{2}}.
\end{cases}
\end{equation}

\end{enumerate}

\end{lemma}

A formal proof of Lemma \ref{lem:Ainf explicit} is omitted in this manuscript. With a little thought, it is not hard to see why \eqref{eq:Ainf(r,theta)=Ainf(-r,pi/2-theta)} is true (cf. Figure \ref{fig:A(-r,t)=1-A(r,t)}) and \eqref{eq:cond exp Ainf=1/2} follows from \eqref{eq:Ainf(r,theta)=Ainf(-r,pi/2-theta)}. The proof of \eqref{eq:Ainf(r,theta) explicit} is elementary, though somewhat long and tedious angle and length
chasing, standard in Euclidean geometry.
The l.h.s. of \eqref{eq:cond exp Ainf=1/2} is interpreted as the conditional expectation
\begin{equation*}
\E_{r}\left[\Ac_{\infty}   \big| \theta\right] := \frac{1}{\sqrt{2}}\int\limits_{-1/\sqrt{2}}^{1/\sqrt{2}} \Ac_{\infty}(r,\theta)dr,
\end{equation*}
that will be found useful below, as will Lemma \ref{lem:Ainf explicit}(i).

\begin{proof}[Proof of Theorem \ref{thm:rand mod single}]

First we prove Theorem \ref{thm:rand mod single}(i). Fix a lattice point $\lambda = (R+r_0)e^{i\theta_0}\in\Gamma_{R}$ in its polar representation \eqref{eq:lambda polar coord} with some $r_0\in [-\tfrac{1}{\sqrt{2}},\tfrac{1}{\sqrt{2}}]$, $\theta_0\in [0,2\pi)$. We claim that, in this case,
\begin{equation}
\label{eq:AR=Ainf+O(1/R)}
\Ac_{R}(\lambda) = \Ac_{\infty}(r_0,\pi/2-\theta_0)+O\left( \frac{1}{R} \right) = \Ac_{\infty}(r_0,\theta_0)+O\left( \frac{1}{R} \right),
\end{equation}
with the constant involved in the $``O"$-notation absolute,
Once we show the first equality of \eqref{eq:AR=Ainf+O(1/R)}, the $2$nd one is immediate from the symmetry \eqref{eq:Ainf(r,theta)=Ainf(-r,pi/2-theta)}, and
the statement of Theorem \ref{thm:rand mod single}(i) will be easily deduced from \eqref{eq:AR=Ainf+O(1/R)} thereafter.

To prove \eqref{eq:AR=Ainf+O(1/R)} we observe the picture, emerging in Figure \ref{fig:sect vs straight}. The number $\Ac_{R}(\lambda)$ is the portion of the square $S(\lambda)$ lying inside the disc $B(R)$, shown to the left. It is approximated by the portion of $S(\lambda)$ bounded by the tangent of $B(R)$ at $Re^{i\theta}$, intersecting the $x$ axis at angle $\pi/2+\theta$, which we tilt to the vertical axis in the right picture of Figure \ref{fig:sect vs straight}. This way, what was $S(\lambda)$ is now tilted by $\pi/2-\theta$, hence the portion of what is now $S(\lambda)$ (namely, the unit square, tilted by $\pi/2-\theta$, centered at $(r,0)$) intersecting the left half-plane is $\Ac_{\infty}(r,\pi/2-\theta)$. Finally, since (in the left picture), the angle, in the polar representation, of all points of $S(\lambda)$ is $\theta+O(1/R)$, the discrepancy area, bounded between the circle $\partial B(R)$ and the tangent line, intersecting $S(\lambda)$, is $O(1/R)$, hence the estimate \eqref{eq:AR=Ainf+O(1/R)}.

\begin{figure}[h]
  \centering
\includegraphics[height=7cm]{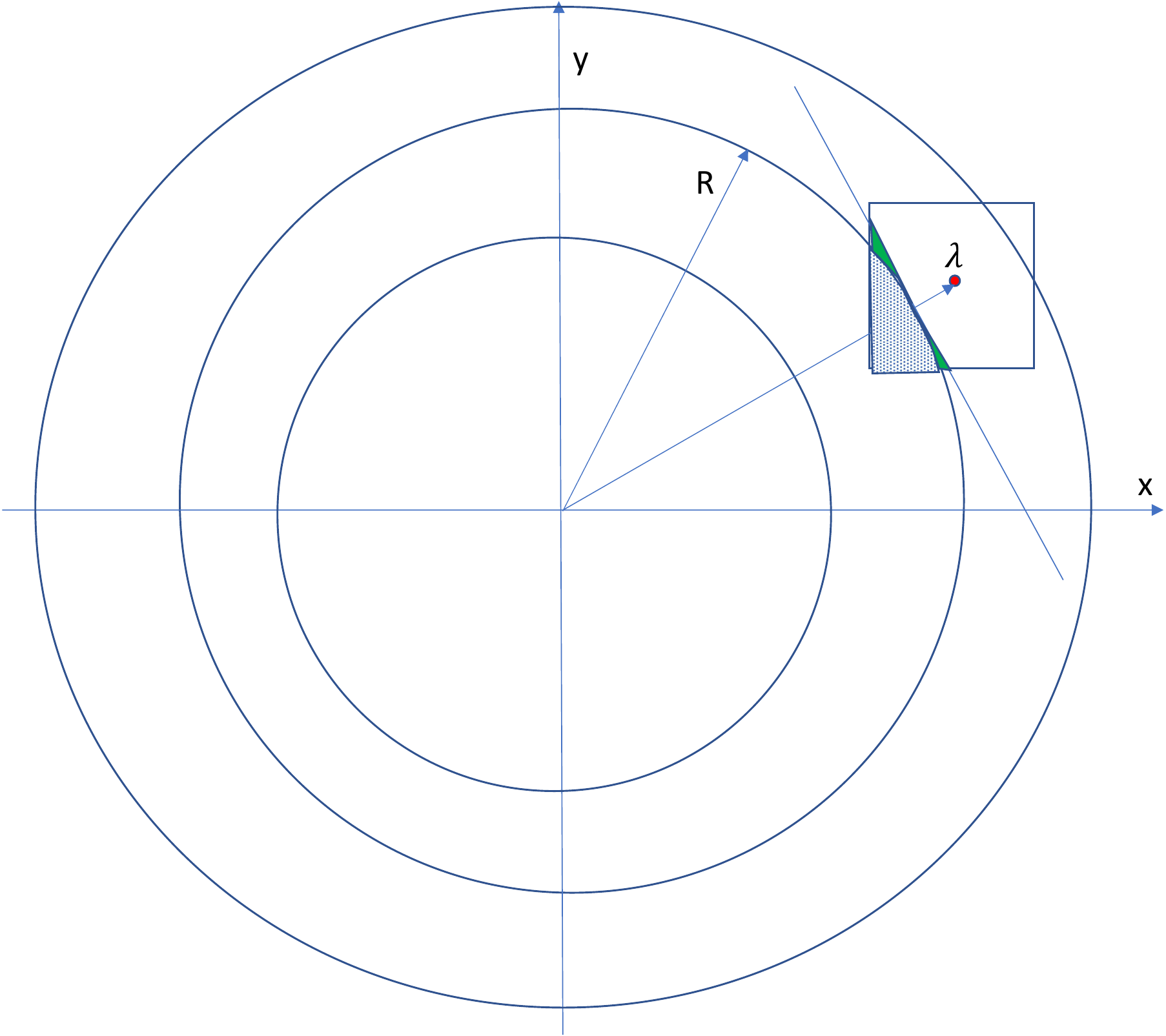}
\includegraphics[height=7cm]{sect_vs_straight_tilted.pdf}
\caption{Approximation of $\Ac_{R}$ with $\Ac_{\infty}$. Left: $\Ac_{R}(\lambda)$ is the portion of the rectangle bounded by the radius-$R$ circle. The discrepancy between these two is designated by the color green. Right: $\Ac_{\infty}(r,\pi/2-\theta)$ is the portion of the tilted square to the left of the $y$ axis, equalling the portion of the square in the left picture bounded by the straight line.}
  \label{fig:sect vs straight}
\end{figure}

Now we deduce the statement of Theorem \ref{thm:rand mod single}(i) from \eqref{eq:AR=Ainf+O(1/R)}. It follows directly from
Proposition \ref{eq:prop:(r,theta) equidist}, via the continuous mapping theorem applied on the functional $$(r,\theta)\mapsto \Ac_{\infty}(r,\theta),$$ continuous on $G$, that the numbers $$\{\Ac_{\infty}(r_{\lambda},\theta_{\lambda})\}_{\lambda\in\Gamma_{R}}$$ converge, in distribution, to $\Ac_{\infty}(r,\theta)$, with $(r,\theta)\in G$ random uniform. The estimate \eqref{eq:AR=Ainf+O(1/R)} implies that so do the numbers
$$\Lambda_{R} = \{\Ac_{R}(\lambda)\}_{\lambda\in\Gamma_{R}},$$ that is the statement of Theorem \ref{thm:rand mod single}(i).

\vspace{2mm}

To prove Theorem \ref{thm:rand mod single}(ii) we observe that, in light of the symmetry \eqref{eq:Ainf(r,theta)=Ainf(-r,pi/2-theta)}, and the symmetry of the square w.r.t. rotation by $\pi/2$, we may express the second moment of $\Ac_{\infty}$ in terms of an integral with a restricted range of $\theta$:
\begin{equation}
\label{eq:var Ainf single int}
\E_{r,\theta}\left[\Ac_{\infty}^2  \right] =  \int\limits_{0}^{\pi/4}\int\limits_{-1/\sqrt{2}}^{1/\sqrt{2}}\Ac_{\infty}(r,\theta)^{2}\frac{dr}{\sqrt{2}}\frac{d\theta}{\pi/4}.
\end{equation}
The statement of Theorem \ref{thm:rand mod single}(ii) follows upon substituting the explicit expression \eqref{eq:Ainf(r,theta) explicit} into
\eqref{eq:var Ainf single int}, taking into account
\begin{equation}
\label{eq:exp r,theta Acinf = 1/2}
\E_{r,\theta}\left[ \Ac_{\infty}\right] = \frac{1}{2},
\end{equation}
which obviously follows from Lemma \ref{lem:Ainf explicit}(i).

\end{proof}

\section{Existence of limit correlations: Proof of Theorem \ref{thm:rand mod joint}} \label{sec:random2}

\subsection{Preparatory results towards the proof of Theorem \ref{thm:rand mod joint}}
\label{sec:prep result}

\begin{notation}
Given $R>0$ and angles $\theta_{1},\theta_{2}\in \Sc^{1}$, where $\Sc^1=\{ x \in \mathbb R^2 : \| x\|=1\}$, we denote $$\Gamma_{R}(\theta_{1},\theta_{2}) = \Gamma_{R}\cap \left\{x=\| x\|e^{i\theta}:\: \theta\in [\theta_{1},\theta_{2}]\right\},$$ to be the set of lattice points lying in the sector inside the annulus $\Dc_{R;1/\sqrt{2}}$ bounded by the rays $\theta=\theta_{j}$, $j=1,2$, in polar coordinates.
\end{notation}

\begin{lemma}
\label{lem:lat pnt cons rect disc}
Let \eqref{eq:GammaR ordered list} be the enumeration of the lattice points $\Gamma_{R}$ lying in the annulus $\Dc_{R;1/\sqrt{2}}$, and $$\lambda_{j}=(R+r_{j})e^{i\theta_{j}}$$ be the polar representation \eqref{eq:lambda polar coord} of $\lambda_{j}$,
ordered as in \eqref{eq:GammaR ordered list}. Then:

\begin{enumerate}[i)]

\item There exists an absolute number $C>0$ sufficiently large so that for all $\theta_{1},\theta_{2}\in \Sc^{1}$ satisfying
\begin{equation}
\label{eq:t2-t1>=1/R assump}
|\theta_{1}-\theta_{2}|> \frac{1}{R}
\end{equation}
(l.h.s. understood as the distance on the circle, i.e. in the $2\pi$-periodic sense), one has $$\#\Gamma_{R}(\theta_{1},\theta_{2}) \le C\cdot |\theta_{1}-\theta_{2}|\cdot R.   $$

\item There exist absolute numbers $c>0$ sufficiently small and $C>0$ sufficiently large, so that for all $\theta_{1},\theta_{2}\in \Sc^{1}$ satisfying $$|\theta_{1}-\theta_{2}|> \frac{C}{R},$$
one has $$\#\Gamma_{R}(\theta_{1},\theta_{2}) \ge c\cdot |\theta_{1}-\theta_{2}|\cdot R.   $$

\item There exists an absolute constant $C>0$ so that $$|\theta_{j}-\theta_{j+k}| \le C\cdot \frac{k}{R},$$ for every $1\le j\le K(R)$, $R>0$, and $k \ge 1$.

\item There exist absolute constants $R_{0}>0$, $k_{0}\ge 1$ and $c>0$ so that $$|\theta_{j+k}-\theta_{j}| \ge c\cdot \frac{k}{R},$$
for every $1\le j\le K(R)$, $R>R_{0}$, $k\ge k_{0}$.

\item The magnitudes of the differences $\|\lambda_{j+1}-\lambda_{j}\|$ are bounded by an absolute constant for every $1 \le j \le K(R)$, $R>0$.

\end{enumerate}
\end{lemma}

Part (iv) of Lemma \ref{lem:lat pnt cons rect disc} is added for completeness of the treatise of the subject and will not be used in the rest of the paper.

\begin{lemma}
\label{lem:tj+1-tj>>1/R^2}
Let \eqref{eq:GammaR ordered list} be the ordering of $\Gamma_{R}$. Then, for a density one set of $\lambda_{j}$, the angle between consecutive lattice points grows relatively to $\frac{1}{R^{2}}$ as $R \rightarrow \infty$. Explicitly, for every $C>0$, $$\frac{\#\left\{j\le K:\: |\theta_{j+1}-\theta_{j}| \le \frac{C}{R^{2}}  \right\}}{K} \rightarrow 0,$$   where $\theta_{j+1}-\theta_{j}$ is understood modulo $2\pi$.
\end{lemma}

We introduce the following notation.

\begin{notation} \label{not:qc}
Let $R>0$ and $\theta\in [0,2\pi)$ be an angle, $k\ge 1$ a positive integer, and $C'>2C$
with $C>0$ the sufficiently large constant given in Lemma \ref{lem:lat pnt cons rect disc}(iii).
Denote the rectangle $$\widetilde{\Qc}_{k}(\theta) = [0,C'k]\times [-1/\sqrt{2},1/\sqrt{2}],$$ and
$$\Qc_{R,k}(\theta):=T_{\theta+\pi/2}\widetilde{\Qc}_{k}(\theta) + R\cdot e^{i\theta}$$ its $(\theta+\pi/2)$-tilt, shifted by $R\cdot e^{i\theta}$.
\end{notation}

The following proposition asserts that, at the scale commensurable to $1/R$, in the vicinity of generic lattice points one may approximate the nearby sector with a rectangle of the same direction, leaving the set of corresponding lattice points unimpaired.

\begin{proposition}
\label{prop:rect pert arc lat pnts}
For every positive integer $k\ge 1$ and $R>0$ there exists a set of lattice points $\widetilde{\Gamma}_{R}\subseteq \Gamma_{R}$ satisfying:
\begin{enumerate}[i)]

\item The set $\widetilde{\Gamma}_{R}$ is of full asymptotic density inside $\Gamma_{R}$. Equivalently,
$$\#\left(\Gamma_{R}\setminus \widetilde{\Gamma}_{R}\right)= o(R), \qquad R \rightarrow \infty.$$

\item
For a lattice point $$\lambda_{j}=(R+r_{j})e^{i\theta_{j}}\in \widetilde{\Gamma}_{R}$$ let $\Qc$ be the (tilted) rectangle $\Qc=\Qc_{R,k}(\theta_{j})$, and
$$k'=k'(k;R,\theta_{j}):= \#\Qc\cap\Z^{2}-1.$$
Then one has
\begin{equation}
\label{eq:k<=k'<<k}
k\le k'\ll k,
\end{equation}
and
$$\Gamma_{R}(\theta_{j},\theta_{j+k'}) = \Qc\cap\Z^{2},$$ i.e. the $(k'+1)$ lattice points lying in the sector $\Gamma_{R}(\theta_{j},\theta_{j+k'})$
(inclusive of $\lambda_{j}$) are precisely those lying in the rectangle $\Qc$.

\item
Further, for $\lambda_{j}\in \widetilde{\Gamma}_{R}$ let $\Qc$ and $k'$ be as in part (ii). Then the ordering $\lambda_{j},\lambda_{j+1},\ldots,\lambda_{j+k'}$
in $\Gamma_{R}(\theta_{j},\theta_{j+k'})$ (i.e. nondecreasing w.r.t. $\theta_{j}$)
coincides with the internal ordering of $\Qc\cap\Z^{2}$ (i.e. w.r.t. the projection of $\lambda\in \Qc\cap\Z^{2}$ onto the long side of $\Qc$).

\item For every $\lambda_{j}\in \widetilde{\Gamma}_{R}$ and $0\le k''\le k'$, one has
$$\Ac_{R}(\lambda_{j+k''}) =  \Ac_{\infty,k''}(r_{j},\theta_{j})+ O\left(\frac{k}{R} \right),$$
where $\Ac_{\infty,k''})$ is as in \eqref{eq:Ac inf k def}, and the constant involved in the $``O"$-notation is absolute.

\end{enumerate}

\end{proposition}

\subsection{Proof of Theorem \ref{thm:rand mod joint} assuming lemmas \ref{lem:lat pnt cons rect disc}-\ref{lem:tj+1-tj>>1/R^2} and
Proposition \ref{prop:rect pert arc lat pnts}}

\begin{proof}[Proof of Theorem \ref{thm:rand mod joint}]

We will prove directly that, as $R\rightarrow\infty$,
\begin{equation}
\label{eq:1/K(R)sum->Cktild}
\frac{1}{K(R)}\sum\limits_{j=1}^{K(R)} \left(\Ac_{R}(\lambda_{j})\cdot \Ac_{R}(\lambda_{j+k}) - \frac{1}{4}\right) \rightarrow \widetilde{\Cc}_{k},
\end{equation}
which contains both assertions of Theorem \ref{thm:rand mod joint}.
Let $\widetilde{\Gamma}_{R}\subseteq\Gamma_{R}$ be the subset of $\Gamma_{R}$ described
in Proposition \ref{prop:rect pert arc lat pnts}. Since, as $R\rightarrow\infty$, the sets $$\{(\theta_{j},r_{j}):\: \lambda_{j}\in\Gamma_{R}\}$$ equidistribute in $G$ by Proposition \ref{eq:prop:(r,theta) equidist}, so do the sets $$\{(\theta_{j},r_{j}):\: j\le K(R),\, \lambda_{j}\in\widetilde{\Gamma}_{R}\},$$
by the property of $\widetilde{\Gamma}_{R}$ asserted by Proposition \ref{prop:rect pert arc lat pnts}(i).
Then as $R \rightarrow \infty$
\begin{equation}
\label{eq:sumAA func}
\begin{split}
&\sum\limits_{j=1}^{K(R)} \left(\Ac_{R}(\lambda_{j})\cdot \Ac_{R}(\lambda_{j+k}) - \frac{1}{4}\right) =
\sum\limits_{j:\: \lambda_{j}\in\widetilde{\Gamma}_{R}} \left(\Ac_{R}(\lambda_{j})\cdot \Ac_{R}(\lambda_{j+k}) - \frac{1}{4}\right) + o(R)
\\&=
\sum\limits_{j:\: \lambda_{j}\in\widetilde{\Gamma}_{R}} \left(\Ac_{\infty}(r_{j},\theta_{j})\cdot \Ac_{\infty,k}(r_{j},\theta_{j}) - \frac{1}{4}\right) + o(R),
\end{split}
\end{equation}
by Proposition \ref{prop:rect pert arc lat pnts}(iv), and \eqref{eq:K(R) *R}.

Now recall that $G$ is the rectangle \eqref{eq:G rect def}, and $$\rho_{k}(\cdot,\cdot):G\rightarrow [-1/\sqrt{2},1/\sqrt{2}]$$ is the map \eqref{eq:rho proj def}, so that \eqref{eq:Ainfk,r=Ainfrk} is satisfied. We claim that, for every $k\ge 1$, the map $\rho_{k}$ is continuous outside of a measure $0$ set in $G$. First, it is clear that $\rho_{k}$ is continuous outside the union
$$\Vc:=\bigcup\limits_{k'\le k}\rho_{k'}^{-1}\left(\left\{\pm \frac{1}{\sqrt{2}}\right\}\right)\subseteq G,$$ since by perturbing $(r,\theta)$,
there are no new lattice points that emerge to the relevant rectangle $\Rc$, nor lattice points that disappear from it. Moreover, given $\theta$,
the sides of $\Rc$ contain a lattice point for a discrete set of $r$ only. That implies that the measure of $\Vc$ is $0$ w.r.t. the measure $dU:=\frac{drd\theta}{\sqrt{2}\pi}$ on $G$. Since $\Ac_{\infty}(\cdot,\cdot)$ is continuous everywhere (being piecewise $C^{1}$), the map
$$(r,\theta)\mapsto \Ac_{\infty,k}(r,\theta)=\Ac_{\infty}(\rho_{k}(r,\theta),\theta)$$ is continuous outside
a set of probability $0$ (again, w.r.t. $U$), and so is the map
$$(r,\theta)\mapsto \Ac_{\infty}(r,\theta)\cdot \Ac_{\infty,k}(r,\theta)-\frac{1}{4}.$$

Since, as it was mentioned above, the random vectors $$\{(\theta_{j},r_{j}):\: j\le K(R),\, \lambda_{j}\in\widetilde{\Gamma}_{R}\}_{R>0}$$  equidistribute in $G$,
we may invoke the continuous mapping theorem to yield
the convergence, in distribution, of the {\em bounded} random variables
$\left\{X_{R}:=\Ac_{\infty}(r_{j},\theta_{j})\cdot \Ac_{\infty,k}(r_j,\theta_j) - \frac{1}{4}\right\}_{\lambda_{j}\in\widetilde{\Gamma}_{R}} $
to the random variable
$$\left\{X_{\infty}:=\Ac_{\infty}(r,\theta)\cdot \Ac_{\infty,k}(r,\theta) - \frac{1}{4}\right\}$$ w.r.t. the uniform measure $U$ on $G$.
Since $X_{R}$ (and $X_{\infty}$) are bounded, convergence in distribution implies ~\cite[Theorem 3.5]{Billingsley} the convergence of expectations, i.e.
\begin{equation}
\label{eq:EXR->EXinf}
\E[X_{R}] \rightarrow \E[X_{\infty}].
\end{equation}
We identify the r.h.s. of \eqref{eq:EXR->EXinf} as
$$\widetilde{\Cc}_{k}:= \E_{r,\theta}\left[\Ac_{\infty}(r,\theta)\cdot \Ac_{\infty,k}(r,\theta)\right]-\frac{1}{4}$$ in \eqref{eq:Cktild def},
and the l.h.s. of \eqref{eq:EXR->EXinf} as
\begin{equation*}
\begin{split}
&\frac{1}{\#\widetilde{\Gamma}_{R}} \sum\limits_{j:\: \lambda_{j}\in\widetilde{\Gamma}_{R}} \left(\Ac_{\infty}(r_{j},\theta_{j})\cdot \Ac_{\infty,k}(r_{j},\theta_{j}) - \frac{1}{4}\right)
\\&= \frac{1}{K(R)}\sum\limits_{j=1}^{K(R)} \left(\Ac_{R}(\lambda_{j})\cdot \Ac_{R}(\lambda_{j+k}) - \frac{1}{4}\right) +o(1),
\end{split}
\end{equation*}
as $R\rightarrow \infty$,
by \eqref{eq:sumAA func} and Proposition \ref{prop:rect pert arc lat pnts}(i). That concludes the proof of \eqref{eq:1/K(R)sum->Cktild}, and, therefore, of Theorem \ref{lem:lat pnt cons rect disc}.
\end{proof}

\section{Lattice points in sectors and the vanishing of correlations: proofs of Proposition \ref{eq:prop:(r,theta) equidist}, theorems \ref{thm:Ck->0 average},  \ref{thm:Ck->0 constrained},  \ref{thm:avgconj}, and \ref{thm:mixconj}} \label{sec:lattice-cor}

\subsection{Lattice points in sectors} Recall the notation $\lambda=(R+r_{\lambda})e^{i\theta_{\lambda}}$
given in \eqref{eq:lambda polar coord}.
While our ultimate goal is to understand the limiting behavior of $\pazocal A_R(\lambda)$ and its correlations, as we have seen, this problem is closely related to the joint distribution of $(r_{\lambda},\theta_{\lambda})$ as $\lambda$ varies over lattice points in $\Gamma_R$, and their correlations. In this section we will first record several estimates describing the distribution of lattice points within $\Gamma_R$ then show how to use these results to prove Proposition \ref{eq:prop:(r,theta) equidist} as well as theorems \ref{thm:Ck->0 average},  \ref{thm:Ck->0 constrained},  \ref{thm:avgconj}, and \ref{thm:mixconj}.

Our first result provides an estimate for the number of lattice points in $\Gamma_R$ which also lie in a narrow sector. In particular, the result gives an asymptotic formula for the number of such lattice points provided the length of the sector is $\ge R^{1/2+o(1)}$ and establishes that the distribution of $(r_{\lambda},\theta_{\lambda})_{\lambda \in \Gamma_R}$ as $R \rightarrow \infty$ tends to $(r,\theta)$ where $(r,\theta)$ is random uniform on $[-\tfrac{1}{\sqrt{2}}, \tfrac{1}{\sqrt{2}}]\times \mathbb R/(2\pi \mathbb Z)$.
\begin{proposition} \label{prop:equid}
Let $\varepsilon>0$.
Let $c,d \in \mathbb R$ with
$c<d$ and $I \subseteq [-\frac{1}{\sqrt{2}}, \frac{1}{\sqrt{2}}]$ be a closed interval. Also let $\theta \in \mathbb R$. Suppose that $R^{-1/2} \le |d-c| \le 2\pi $. Then
uniformly in $c,d,I$ we have that
\[
\# \bigg\{ \lambda \in \Gamma_R :  \theta_{\lambda}\Mod {2\pi} \in [ c,d ] \, \& \, r_{\lambda} \in I \bigg\}= 2\pi \sqrt{2} R \cdot \frac{|d-c|}{2\pi} \frac{|I|}{\sqrt{2}}+O\bigg(|d-c|^{1/3}R^{2/3+\varepsilon}\bigg).
\]
\end{proposition}
Taking $I=[-\frac{1}{\sqrt{2}}, \frac{1}{\sqrt{2}}]$ and $[c,d]=[0,2\pi]$ this result recovers the classical estimate due to independently Voronoi,
Sierpinski, and van der Corput. We have made little effort to optimize the range of $|d-c|$ for which the asymptotic holds and it would be interesting to increase this range.
Using Proposition \ref{prop:equid} we immediately
get Proposition \ref{eq:prop:(r,theta) equidist}.
\begin{proof}[Proof of Proposition \ref{eq:prop:(r,theta) equidist}]
Recall that $K=2\sqrt{2}\pi R+O(R^{2/3})$. Apply Proposition \ref{prop:equid} to get that
\begin{equation}\label{eq:lambda theta N=1}
\frac1K\# \bigg\{ \lambda \in \Gamma_R :  \theta_{\lambda} \Mod {2\pi} \in [ c,d ] \, \& \, r_{\lambda} \in I \bigg\} =  \frac{|d-c|}{2\pi} \frac{|I|}{\sqrt{2}}+O\bigg(|d-c|^{1/3}R^{-1/3+\varepsilon}\bigg).
\end{equation}
The main term on the r.h.s. of \eqref{eq:lambda theta N=1} equals $\mathbb P((r,\theta) \in I \times [c,d])$ where $r,\theta$ are independent random variables with $r$ uniformly distributed on $[-\tfrac{1}{\sqrt{2}},\tfrac{1}{\sqrt{2}}]$ and $\theta$ uniformly distributed on $\mathbb R/(2\pi \mathbb Z)$.
\end{proof}

\vspace{2mm}

The next result concerns the number of expected lattice points $\mu \in \Gamma_R$ such that $\theta_{\lambda}-\theta_{\mu} \Mod{2\pi}$ lies in $[\tfrac{k}{\pazocal K}, \tfrac{k+2\pi}{\pazocal K}]$, when averaged w.r.t. $\lambda \in \Gamma_R$, provided $k$ tends to infinity with $R$. In contrast with
Proposition \ref{prop:equid} (and Theorem \ref{thm:varbd} below), for each $\lambda\in \Gamma_R$ the lattice points $\mu$ which we are counting here lie within a region with \textit{area equal to one} and subsequently contains a bounded number of lattice points (the diameter of the region is bounded). Nevertheless we obtain an asymptotic formula when $k \rightarrow \infty$.

\begin{theorem} \label{thm:jointdistribution}
Let $\varepsilon>0$. Also, let $0 < \eta < 1/3$.
Suppose $1 \le k \le R^{\frac{1}{2+5\eta)}-\varepsilon}$.
Then uniformly for intervals $I_1,I_2 \subseteq [-\tfrac{1}{\sqrt{2}},\tfrac{1}{\sqrt{2}}]$ and
$J=[c,d] \subseteq \mathbb R$ with $|J| \le 2\pi$ we have that
\begin{equation} \label{eq:count}
\begin{split}
&\frac{1}{\pazocal K}\# \bigg\{ \lambda,\mu \in \Gamma_R :   \theta_{\lambda}-\theta_{\mu} \Mod {2\pi}\in \bigg[ \frac{k}{\pazocal K}, \frac{k+2\pi}{\pazocal K}\bigg] \, \& \,
(r_{\lambda},r_{\mu}, \theta_{\lambda} \Mod {2\pi}) \in I_1\times  I_2 \times J \bigg\}\\
&\qquad \qquad \qquad =\frac{|I_1|}{\sqrt{2}} \frac{|I_2|}{\sqrt{2}} \frac{|J|}{2\pi} +O\bigg(\frac{1}{k^{\eta-\varepsilon}}+\frac{1}{k^{\frac12-\frac{3\eta}{2} -\varepsilon}}+\frac{k^{2+5\eta}}{R^{1-\varepsilon}}+\frac{k^{7\eta/2}}{R^{1/2-\varepsilon}} \bigg).
\end{split}
\end{equation}
\end{theorem}

For the choice $\eta=1/5$, and assuming that $k \le R^{5/16}$, the error term in \eqref{eq:count} is $\ll k^{-1/5+\varepsilon}$.

\begin{remark} \label{rem:correlation}
Taking $I_1=I_2=[-\tfrac{1}{\sqrt{2}},\tfrac{1}{\sqrt{2}}]$ and $J=[0,2\pi)$ shows that the mean value, averaging over $\lambda \in \Gamma_R$, of the number of lattice points $\mu \in \Gamma_R$ with $  \theta_{\lambda}-\theta_{\mu} \Mod {2\pi} \in [ \tfrac{k}{\pazocal K}, \tfrac{k+2\pi}{\pazocal K}]$ is asymptotically equal to one as $k \rightarrow \infty$ for $k \le R^{\frac{1}{2(1+2\eta)}-\varepsilon}$.
\end{remark}

In the case where $N$ is large and $N=o(R)$ we also study the distribution of lattice points in $\Gamma_R$
whose angles lie within an arc of length $1/N$ centered at $\theta \Mod{2\pi}$. The expected number of lattice points in such a region is asymptotically equal to its area, where we average over the position of our arc $\theta \Mod{2\pi}$. Our main result in this setting provides an {\em asymptotic} formula for the variance.

\begin{notation} Given $\xi \in \mathbb R$ and $f \in L^1(\mathbb R)$ let  $e(\xi)= e^{2\pi i\xi}$ and
\[
\widehat f(\xi)=\int_{\mathbb R} f(t) e(-t\xi) \, dt.
\]
\end{notation}

\begin{theorem} \label{thm:varbd}
Let $\delta_2>0$ be fixed. Then for $I \subseteq [-\tfrac{1}{\sqrt{2}},\tfrac{1}{\sqrt{2}}]$, $c<d$, $R^{-1} \le |d-c| \le R^{-9/10-\delta_2}$, we have that
\begin{equation} \label{eq:varbd}
\begin{split}
\int_0^{\pi/2} \bigg| \sum_{\substack{ \lambda \in \Gamma_R, r_{\lambda} \in I \\
(\theta_{\lambda}-\ \theta) \Mod{2\pi} \in [c, d ]}} 1- R  |d-c| |I|  \bigg|^2 \frac{d\theta}{\pi/2} =  R|d-c| \, \pazocal D(I) +O( ( R|d-c|)^{1-9\delta_2/10} ),
\end{split}
\end{equation}
where
\[
\pazocal D(I)=  \frac{1}{2\pi^2}\sum_{ \lambda \in \mathbb Z^2 \setminus \{0\}} \frac{|\widehat \chi_I(\| \lambda \|)|^2}{\| \lambda \|}.
\]
\end{theorem}

Theorem \ref{thm:varbd} implies that, for generic $\theta$, the size of the error term in the lattice point counting problem is the square root of the main term.

\subsection{Proof of theorems \ref{thm:Ck->0 constrained} and \ref{thm:mixconj}}
\begin{proof}[Proof of Theorem \ref{thm:mixconj} assuming Theorem \ref{thm:jointdistribution}]
Let $\epsilon>0$. We have by \eqref{eq:AR=Ainf+O(1/R)} that for $\lambda \in \Gamma_R$
\[
\pazocal{A}_R(\lambda)=\pazocal{A}_{\infty}(r_{\lambda},\pi/2-\theta_{\lambda})+O(R^{-1}).
\]
Let $H: [0,1]\times[0,1]\rightarrow \mathbb R$ be a continuous function.
Hence, we see that for $R > R_0(\epsilon)$
\begin{equation}
\label{eq:passtoinfty}
\begin{split}
&\frac{1}{\pazocal K} \sum_{\substack{\lambda, \mu \in \Gamma_R \\ \theta_{\lambda}- \theta_{\mu} \Mod {2\pi} \in [\frac{k}{\pazocal K},\frac{k+2\pi}{\pazocal K}]}} H(\pazocal{A}_R(\lambda), \pazocal{A}_R(\mu))\\
& \qquad \qquad \qquad = \frac{1}{\pazocal K} \sum_{\substack{\lambda, \mu \in \Gamma_R \\ \theta_{\lambda}- \theta_{\mu} \Mod {2\pi} \in [\frac{k}{\pazocal K},\frac{k+2\pi}{\pazocal K}]}}  H(\pazocal{A}_{\infty}(r_{\lambda},\pi/2-\theta_{\lambda}), \pazocal{A}_{\infty}(r_{\mu},\pi/2-\theta_{\mu}))+O(\epsilon)
\end{split}
\end{equation}
where we have used Lemma \ref{lem:lat pnt cons rect disc}(i) to handle the error term. Since $\pazocal{A}_{\infty}$ is uniformly continuous and bounded by $1$ we have that, for $R$ sufficiently large in terms of $\epsilon$, that
\begin{equation} \label{eq:precedingformula1}
H(\pazocal{A}_{\infty}(r_{\lambda},\pi/2-\theta_{\lambda}),\pazocal{A}_{\infty}(r_{\mu},\pi/2-\theta_{\mu}))=H(\pazocal{A}_{\infty}(r_{\lambda},\pi/2-\theta_{\lambda}),\pazocal{A}_{\infty}(r_{\mu},\pi/2-\theta_{\lambda}))+O(\epsilon)
\end{equation}
for $\lambda,\mu$ with $\theta_{\lambda}- \theta_{\mu} \Mod {2\pi} \in \left[\frac{k}{\pazocal K},\frac{k+2\pi}{\pazocal K}\right]$.
Note that Theorem \ref{thm:jointdistribution} implies that as $(\lambda,\mu)$ vary over pairs of points in $\Gamma_R\times \Gamma_R$ with $(\theta_{\lambda}-\theta_{\mu}) \Mod{2\pi} \in [\tfrac{k+2\pi}{\pazocal K},\tfrac{k+2\pi}{\pazocal K}]$ that the joint distribution of $(r_{\lambda},\theta_{\lambda},r_{\mu})$ tends to $(r,\theta,r')$ where $(r,\theta,r')$ is uniform in $\left[-\tfrac{1}{\sqrt{2}},\tfrac{1}{\sqrt{2}}\right] \times \mathbb R/(2\pi) \times \left[-\tfrac{1}{\sqrt{2}},\tfrac{1}{\sqrt{2}}\right]$.
Combining this observation with \eqref{eq:passtoinfty} and \eqref{eq:precedingformula1}
we see that for $R$ sufficiently large that
\begin{equation} \label{eq:dist}
\begin{split}
&\frac{1}{\pazocal K} \sum_{\substack{\lambda, \mu \in \Gamma_R \\  \theta_{\lambda}- \theta_{\mu} \Mod{2\pi} \in [\frac{k}{\pazocal K},\frac{k+2\pi}{\pazocal K}]}} H(\pazocal{A}_R(\lambda), \pazocal{A}_R(\mu))\\
&\qquad \qquad =\int_{-\frac{1}{\sqrt{2}}}^{\frac{1}{\sqrt{2}}} \int_{-\frac{1}{\sqrt{2}}}^{\frac{1}{\sqrt{2}}} \int_0^{2\pi} H(\pazocal{A}_{\infty}(r_1, \theta), \pazocal{A}_{\infty}(r_2,\theta)) \, \frac{dr_1}{\sqrt{2}}\frac{dr_2}{\sqrt{2}} \frac{d\theta}{2\pi}+O(\epsilon),
\end{split}
\end{equation}
where we transformed variables in the integral w.r.t. $\theta$. This completes the proof of Theorem \ref{thm:mixconj}. \end{proof}

\begin{proof}[Proof of Theorem \ref{thm:Ck->0 constrained}]
Using Lemma \ref{lem:Ainf explicit}(i) we can evaluate the integral on the r.h.s. of \eqref{eq:dist} for the choice $H(x,y)=xy$ of the test function, by changing the order of integration and evaluating directly the integrals over $r_1,r_2$. This gives that
\begin{equation} \label{eq:uncorr}
\int_{-\frac{1}{\sqrt{2}}}^{\frac{1}{\sqrt{2}}} \int_{-\frac{1}{\sqrt{2}}}^{\frac{1}{\sqrt{2}}} \int_0^{2\pi} \pazocal{A}_{\infty}(r_1, \theta) \pazocal{A}_{\infty}(r_2,\theta) \, \frac{dr_1}{\sqrt{2}}\frac{dr_2}{\sqrt{2}} \frac{d\theta}{2\pi}=\frac14.
\end{equation}
Hence, we conclude from \eqref{eq:dist} that
\begin{equation} \label{eq:corrform}
\frac{1}{\pazocal K} \sum_{\substack{\lambda, \mu \in \Gamma_R \\ \theta_{\lambda}- \theta_{\mu} \Mod {2\pi} \in [\frac{k}{\pazocal K},\frac{k+2\pi}{\pazocal K}]}} \pazocal{A}_R(\lambda) \pazocal{A}_R(\mu)=\frac14+O(\epsilon).
\end{equation}
Repeating the argument above with straightforward modifications gives
\begin{equation}
\label{eq:exp AR approx}
\frac{1}{\pazocal K} \sum_{\substack{\lambda, \mu \in \Gamma_R \\ \theta_{\lambda}- \theta_{\mu} \Mod {2\pi} \in [\frac{k}{\pazocal K},\frac{k+2\pi}{\pazocal K}]}} \pazocal{A}_R(\lambda) =\frac12+O(\epsilon), \quad \frac{1}{\pazocal K} \sum_{\substack{\lambda, \mu \in \Gamma_R \\ \theta_{\lambda}- \theta_{\mu} \Mod {2\pi} \in [\frac{k}{\pazocal K},\frac{k+2\pi}{\pazocal K}]}} \pazocal{A}_R(\mu) =\frac12+O(\epsilon).
\end{equation}
Using the estimates \eqref{eq:uncorr}, \eqref{eq:corrform} and \eqref{eq:exp AR approx}, and Theorem \ref{thm:jointdistribution} with $I_1=I_2=[-\tfrac{1}{\sqrt{2}}, \tfrac{1}{\sqrt{2}}]$ and $J=[0,2\pi)$
completes the proof of Theorem \ref{thm:Ck->0 constrained} (see Remark \ref{rem:correlation}).
\end{proof}

\subsection{Proof of theorems \ref{thm:Ck->0 average} and \ref{thm:avgconj}}
\begin{proof}[Proof of Theorem \ref{thm:avgconj} assuming Theorem \ref{thm:varbd}]
Let $\Lambda, M$ tend to infinity with $R$ with $M \le R^{1/10}$ and suppose $\Lambda\le M^{1/6}$.
Consider\footnote{Note that the area of the region $\Dc_{R,1/\sqrt{2}} \cap \{ x=\| x \|e^{ i\theta} : \theta \in (0,\frac{2\pi M}{\pazocal K})\}$ is equal to $ \tfrac12 \frac{M}{\sqrt{2}R} ((R+\frac{1}{\sqrt{2}})^2-(R-\frac{1}{\sqrt{2}})^2)=M$. }
\[
S_{\Lambda,M}=\bigg\{ 0 \le \theta \le \pi/2
: \bigg| \sum_{\substack{\lambda \in \Gamma_R \\  \theta_{\lambda}- \theta \Mod {2\pi} \in [0,\frac{ 2\pi M}{\pazocal K}]}} 1-M \bigg| \le \Lambda \sqrt{M}
\bigg\}
\]
and $S_{\Lambda, M}^c =[0,\pi/2] \setminus S_{\Lambda,M}$ .
Recall the notation $\theta_j=\theta_{\lambda_j}$.
Also let
\[\pazocal S_{\Lambda,M}=\bigg\{ 1 \le j \le K :  \theta_j \Mod  {2\pi}  \in S_{\Lambda,M} \bigg\}.
\]
Using
Theorem \ref{thm:varbd} we have that
\begin{equation} \label{eq:measbd}
\tmop{meas}(S_{\Lambda,M}^c) \le \frac{1}{\Lambda^2 M}
\int_{0}^{\pi/2} \bigg| \sum_{\substack{\lambda \in \Gamma_R \\  \theta_{\lambda}- \theta \Mod {2\pi} \in [0,\frac{2\pi M}{\pazocal K}]}} 1-M \bigg|^2 \, \frac{d\theta}{\pi/2} \ll \frac{1}{\Lambda^2}=o(1).
\end{equation}
Hence, using this along with Proposition \ref{prop:equid} we conclude that
\begin{equation} \label{eq:cardinality}
\# \pazocal S_{\Lambda,M}=K(1+o(1)).
\end{equation}

By construction and Lemma \ref{lem:lat pnt cons rect disc}(i)-(ii), for $j \in \pazocal S_{\Lambda,M}$ we have
\begin{equation} \label{eq:located}
 \theta_{j+M} - \theta_j \Mod {2\pi} \in \bigg[ \frac{2\pi M-\Lambda^2 \sqrt{M}}{\pazocal K}, \frac{2\pi M+\Lambda^2 \sqrt{M}}{\pazocal K}\bigg].
\end{equation}
Let $M_1=M+L$, where $L=\Lambda^3 \sqrt{M}$. For $j \in \pazocal S_{\Lambda,M_1} \cap \pazocal S_{\Lambda,M}$
using \eqref{eq:located} along with an analogue of this formula where $M$ is replaced with $M+L$
we have that
\[
\begin{split}
&
\# \bigg(\bigg\{ \lambda_{j+k} : M \le k \le M+L  \bigg\}
\triangle \bigg\{ \lambda \in \Gamma_R :  \theta_j-\theta_{\lambda} \Mod {2\pi} \in \bigg[\frac{2\pi M}{\pazocal K}, \frac{2\pi(M+L)}{\pazocal K} \bigg] \bigg\} \bigg)
\\
 & \qquad  \le  \sum_{\substack{\lambda \in \Gamma_R \\ \theta_{\lambda}- \theta_j \Mod {2\pi} \in [\frac{2\pi M-\Lambda^2 \sqrt{M}}{\pazocal K}, \frac{2\pi M}{\pazocal K}]}} 1+\sum_{\substack{\lambda \in \Gamma_R \\ \theta_{\lambda}- \theta_j \Mod {2\pi} \in [\frac{2\pi M_1 }{\pazocal K}, \frac{2\pi M_1+\Lambda^2\sqrt{M}}{\pazocal K}]}} 1 \ll \Lambda^2 \sqrt{M},
\end{split}
\]
where $ A \triangle B$ denotes the symmetric difference of the sets $A,B$ and the last inequality follows from Lemma \ref{lem:lat pnt cons rect disc}(i).

Let $H:[0,1]^2 \rightarrow \mathbb R$ be a continuous function.
Hence, using the above estimate along with \eqref{eq:cardinality} along with an analogue of this formula with $M$ replaced by $M_1$ we have that
\begin{equation} \label{eq:passoff}
\begin{split}
&\sum_{M \le k \le M+L} \sum_{j=1}^K H( \pazocal A_R(\lambda_j), \pazocal A_R(\lambda_{j+k}))\\
& \qquad \qquad =\sum_{M \le k \le M+L} \sum_{j \in \pazocal S_{\Lambda, M} \cap \pazocal S_{\Lambda, M+L}} H( \pazocal A_R(\lambda_j), \pazocal A_{R} (\lambda_{j+k}))+o\bigg(
\Lambda^3 \sqrt{M} R \bigg) \\
&\qquad \qquad = \sum_{\substack{ \lambda,\mu \in \Gamma_R \\ \theta_{\lambda}- \theta_{\mu} \Mod {2\pi} \in [\frac{2\pi M}{\pazocal K}, \frac{2\pi(M+L)}{\pazocal K}]}} H( \pazocal A_R(\lambda), \pazocal A_R(\mu))+o\bigg(
\Lambda^3 \sqrt{M} R \bigg)+O(\Lambda^2 R \sqrt{M}).
\end{split}
\end{equation}
In the r.h.s. of \eqref{eq:passoff}, using \eqref{eq:dist} with $k=2\pi M,2\pi (M+1),\ldots,2\pi( M+\lfloor L\rfloor)$, where $\lfloor L \rfloor$ is the largest integer smaller than or equal to $L$, together with Lemma \ref{lem:lat pnt cons rect disc} (i), and recalling \eqref{eq:corrform} completes the proof of Theorem \ref{thm:avgconj}.
\end{proof}

\begin{proof}[Proof of Theorem \ref{thm:Ck->0 average}]
Using Theorem \ref{thm:avgconj} with $L=M$ and noting \eqref{eq:exp r,theta Acinf = 1/2} and \eqref{eq:uncorr} we get that $\frac1M \sum_{k=M}^{2M} \pazocal C_k \rightarrow 0$ as $M \rightarrow \infty$.
Since $|\pazocal C_k| \le 1/4$ we can conclude $\sum_{k=1}^M \pazocal C_k =o(M)$ by splitting the sum over intervals of the form $[M/2^j, 2M/2^j]$.
\end{proof}

\section{Proof of the preparatory results of section \ref{sec:prep result}} \label{sec:random-proofs}

We will require an auxiliary notation for annuli and their sectors.

\begin{notation}
\label{not:annuli, sec}

Let $0<R_{1}<R_{2}<+\infty$, and $\theta_{1}<\theta_{2}$ be two angles on the circle, i.e. inequality modulo $2\pi$.

\begin{enumerate}

\item We denote the annulus $$A_{R_{1},R_{2}}=\{x\in\R^{2}:\: R_{1}\le\|x\|\le R_{2}\} .$$

\item We denote the annular sector of $A_{R_{1},R_{2}}$ by
$$A_{R_{1},R_{2}}(\theta_{1},\theta_{2})=\{x\in\R^{2}:\: x=\|x\|e^{i\theta}, \, \theta_{1}\le \theta\le \theta_{2}\}\cap A_{R_{1},R_{2}} .$$

\end{enumerate}

\end{notation}

\begin{proof}[Proof of Lemma \ref{lem:lat pnt cons rect disc}]

First, we prove (i). Let the polygonal domain
$$\Pc_{R}(\theta_{1},\theta_{2})=\bigcup\limits_{\lambda\in \Gamma_{R}(\theta_{1},\theta_{2})}\left( \lambda+\left[-\frac{1}{2},\frac{1}{2} \right]^{2} \right) $$ be the union of all unit squares centered at some $\lambda\in \Gamma_{R}(\theta_{1},\theta_{2})$. Then $\Pc_{R}(\theta_{1},\theta_{2})$ is contained in the $\sqrt{2}$-neighborhood of the annular sector
$A_{R_{1},R_{2}}(\theta_{1},\theta_{2})$ with $R_{1}=R-1/\sqrt{2}$, $R_{2}=R+1/\sqrt{2}$, that itself is contained in a slightly bigger annular sector:
\begin{equation*}
\Pc_{R}(\theta_{1},\theta_{2})\subseteq A_{R_{1},R_{2}}(\theta_{1},\theta_{2})_{1/\sqrt{2}}:= \left\{x\in\R^{2}:\: d(x,A_{R_{1},R_{2}}(\theta_{1},\theta_{2}))\le 1/\sqrt{2}\right\}
\subseteq A_{R_{3},R_{4}}(\theta_{3},\theta_{4}),
\end{equation*}
where $R_{3}=R_{1}-1/\sqrt{2}=R-\sqrt{2}$, $R_{4}=R_{2}+1/\sqrt{2}=R+\sqrt{2}$, $\theta_{3}=\theta_{1}-C/R$, $\theta_{4}=\theta_{2}+C/R$ with sufficiently large absolute constant $C>0$, and $d(\cdot,\cdot)$ designates the Euclidean distance between a point and a set.
Then
\begin{equation*}
\begin{split}
&\# \Gamma_{R}(\theta_{1},\theta_{2})=\area\left(\Pc_{R}(\theta_{1},\theta_{2})\right) \le \area(A_{R_{3},R_{4}}(\theta_{3},\theta_{4})) \ll R\cdot (\theta_{4}-\theta_{3})
\\&\le R\cdot (\theta_{2}-\theta_{1}+2C/R) \ll R(\theta_{2}-\theta_{1}),
\end{split}
\end{equation*}
by the assumption \eqref{eq:t2-t1>=1/R assump} on $\theta_{1},\theta_{2}$.

\vspace{2mm}

Next, we prove (ii). It is sufficient to prove that there exists a sufficiently large constant $C>0$ so that, for every $\theta\in [0,2\pi)$,
an annular sector $A_{R_{1},R_{2}}(\theta,\theta+C/R)$ contains at least one lattice point. We will assume w.l.o.g. that $\theta\in [0,\pi/2)$, and, further,
that $\theta \in [\pi/4,\pi/2)$ (otherwise flip the coordinate axes). Let $$x=(x_{1},x_{2})=R\cdot e^{i\theta},$$ and take
$m:=\lfloor x_{1}\rfloor$. We denote the unit square $T:=[m-1,m]\times [ 0,1]$, and consider its vertical shifts by an
integer: $$\Tc:=\{(0,n)+T:\:n\in\Z_{\ge 0}\}.$$ If we can show that one of the elements of $\Tc$ is fully contained in $A_{R_{1},R_{2}}$, then it means that
we found a lattice point $\lambda\in A_{R_{1},R_{2}}$, whose $x$ coordinate differs by at most $2$ from $x_{1}$. Thus, bearing in mind $\theta>\frac{\pi}{4}$,
that implies that $\lambda\in A_{R_{1},R_{2}}(\theta,\theta+C/R)$ with $C>0$ sufficiently large absolute, and we are done.

\begin{figure}[h]
  \centering
\includegraphics[height=13cm]{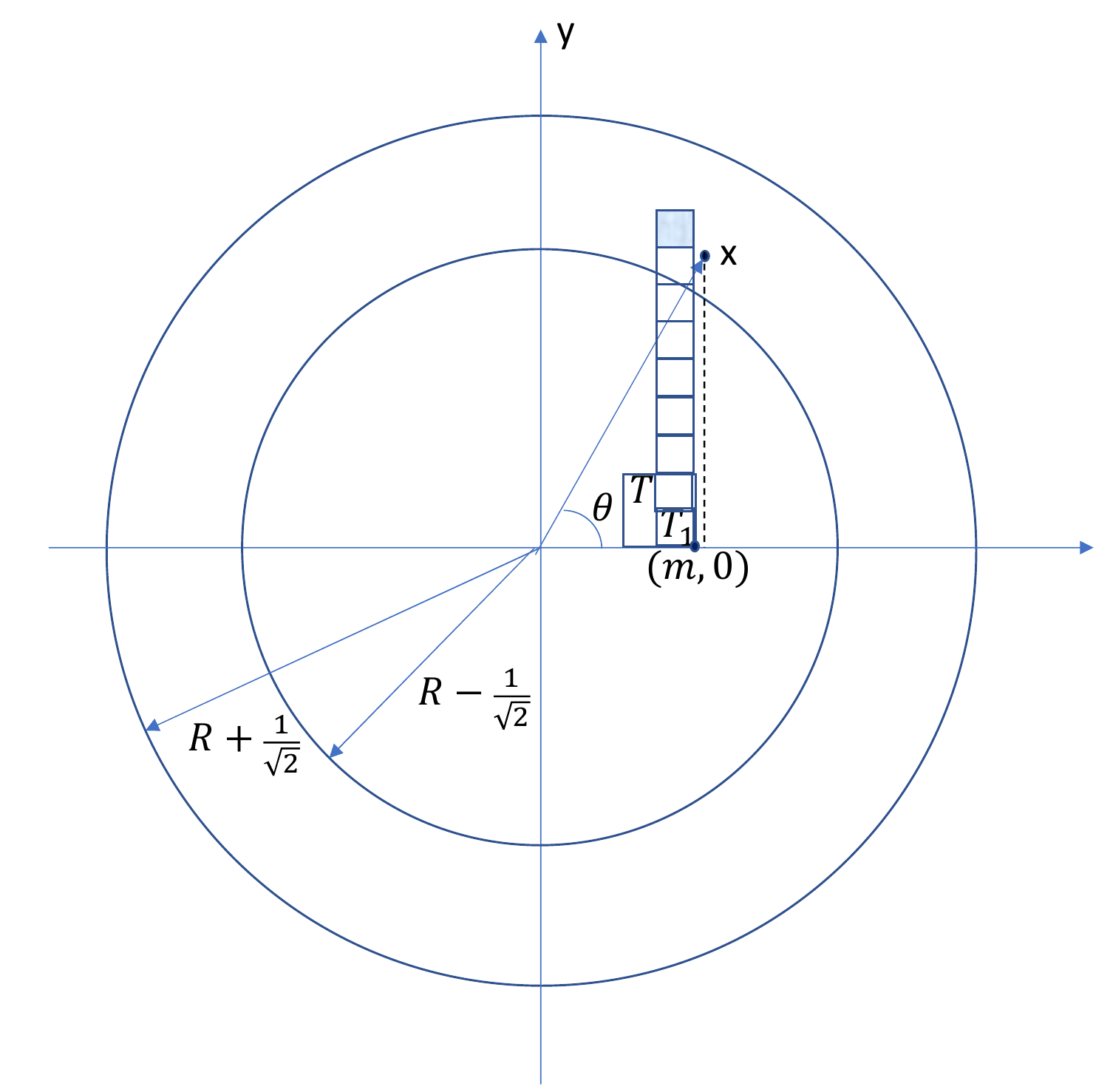}
\caption{Finding a lattice point near a given one. One of the corners of the solid square is guaranteed to be a lattice point.}
  \label{fig:int approx}
\end{figure}

However, the problem is that the diagonal length of $T$ is $\sqrt{2}$, and we require this square to fit inside an annulus of the same width, by
integral translation. The upshot is that we may repeat the same argument with {\em quarters} of $T$, see the illustration in
Figure \ref{fig:int approx}.
Namely, let $T_{1}=[m-1/2,m]\times [0,1/2]$, and consider the family
$$\Tc_{1}:=\{B_{n}=(0,n/2)+T_{1}:\:n\in\Z_{\ge 0}\}$$ of half-integer vertical shifts. Since each element of $\Tc_{1}$ contains at least one integral corner, if for some $n$,
$B_n\in\Tc_{1}$, we have that $B_n\subseteq A_{R_{1},R_{2}}$, then that yields a lattice point $\lambda\in B_n$, and we are done, by the above logic. To this end,
we argue that the diagonal (and the diameter) of the rectangle $B_{n}\cup B_{n+1}$ for each $n\in\Z_{\ge 0}$ is $\frac{\sqrt{5}}{2}<\sqrt{2}$; hence there exists an $n$ for which $B_{n}$ intersects the
internal boundary $\partial B(R-1/\sqrt{2})$ of $A_{R_{1},R_{2}}$ and $B_{n+1}$ does not intersect the boundary, in which case $B_{n+1}$ is fully contained in $A_{R_{1},R_{2}}$, concluding the proof of part (ii).

Part (iii) is a direct consequence of part (ii). Part (iv) is a direct consequence of part (i). Part (v) is a direct consequence of part (iii) and the triangle inequality:
\begin{equation*}
\begin{split}
&\left\|\lambda_{j+1}-\lambda_{j}\right\| \le \left\| \|\lambda_{j+1}\|e^{i\theta_{j}}-\lambda_{j} \right\| + \left\|\lambda_{j+1}- \|\lambda_{j+1}\|e^{i\theta_{j}}\right\|
\\&\le \big|\|\lambda_{j+1}\|-\|\lambda_{j}\|\big|+(R+1/\sqrt{2})\cdot \left|\theta_{j+1}-\theta_{j}\right| \ll 1.
\end{split}
\end{equation*}

\end{proof}

\begin{figure}[h]
  \centering
\includegraphics[height=13cm]{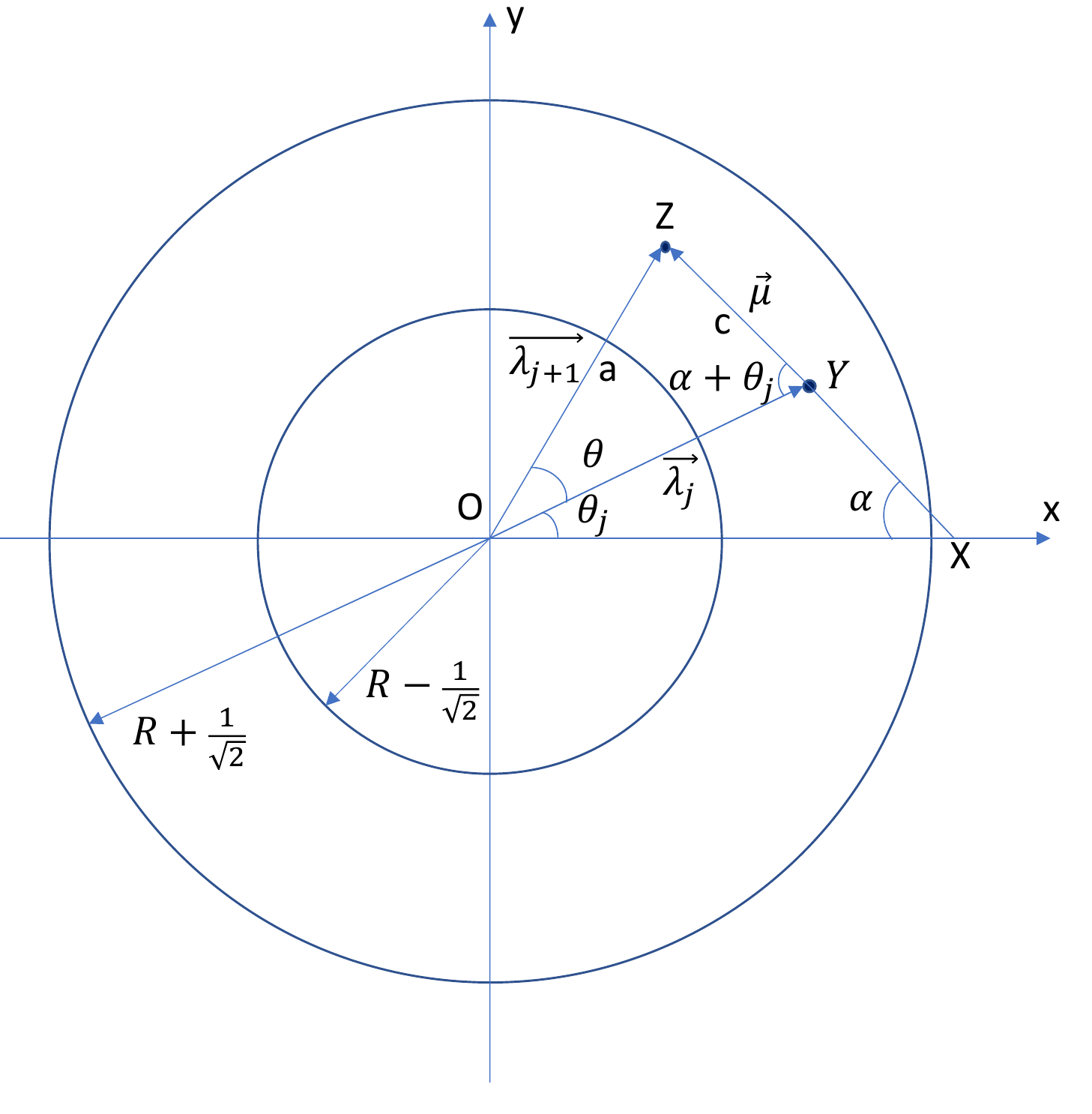}
\caption{Small angle $\theta_{j+1}-\theta_{j}$ forces $\theta_{j}$ to be almost parallel to $\mu=\lambda_{j+1}-\lambda_{j}$.}
  \label{fig:angle chasing}
\end{figure}

\begin{proof}[Proof of Lemma \ref{lem:tj+1-tj>>1/R^2}]

Let $C>0$ be given, and define $\Sigma_{R,C}$ by
\begin{equation*}
\Sigma_{R,C}=\left\{\lambda_{j}\in\Gamma_{R}:\: |\theta_{j+1}-\theta_{j}| \le \frac{C}{R^{2}}  \right\}.
\end{equation*}
Lemma \ref{lem:tj+1-tj>>1/R^2} asserts that, for every $C>0$, the relative density of $\Sigma_{R,C}\subseteq \Gamma_{R}$ vanishes as $R\rightarrow\infty$.

Now take $\gamma>0$ to be a sufficiently large constant so that, as guaranteed by Lemma \ref{lem:lat pnt cons rect disc}(v), for every $j\le K$, one has $\|\lambda_{j+1}-\lambda_{j}\|< \gamma$,
or, put otherwise, for every
$j\le K$, all the vectors $$\lambda_{j+1}-\lambda_{j}\in B(\gamma)$$ belong to the radius-$\gamma$ centered disc. Then we may decompose
\begin{equation*}
\Gamma_{R}= \bigcup\limits_{\mu\in B(\gamma)\cap\Z^{2}}\Delta_{R}(\mu),
\end{equation*}
with
\begin{equation*}
\Delta_{R}(\mu):= \{\lambda_{j}\in\Gamma_{R}:\: \lambda_{j+1}-\lambda_{j}=\mu\},
\end{equation*}
and, accordingly,
\begin{equation*}
\Sigma_{R,C}=\bigcup\limits_{\mu\in B(\gamma)\cap\Z^{2}}\Delta_{R}(\mu)\cap \Sigma_{R,C},
\end{equation*}
a (fixed) finite union of subsets $\Delta_{R}(\mu)\cap \Sigma_{R,C}$ of $\Gamma_{R}$. It is therefore sufficient to prove that, for every $\mu\in B(\gamma)$ the subset
$$\Delta_{R}(\mu)\cap \Sigma_{R,C}\subseteq\Gamma_{R}$$ is of asymptotically vanishing relative density.

\vspace{2mm}

Let
\begin{equation}
\label{eq:lambda in Delta and Sigma}
\lambda_{j}\in \Delta_{R}(\mu)\cap \Sigma_{R,C},
\end{equation}
and assume w.l.o.g., that the corresponding angle $\theta_{j}\in [0,\pi/2)$, i.e. $\lambda_{j}$ belongs to the $1$st quadrant in
$\R^{2}$. We think of $\mu$ as the vector between $\lambda_{j}$ and $\lambda_{j+1}$, and claim that the condition \eqref{eq:lambda in Delta and Sigma} on $\lambda_{j}$ induces a constraint
on $\theta_{j}$ that is only satisfied by
non-generic $\lambda_{j}$ that are {\em almost parallel} to $\mu$, in a quantitative sense to be made precise. Let $O$ be the origin, and $X$ be the intersection of
the continuation of $\mu$ with the $x$ axis, $Y=\lambda_{j}$, $Z=\lambda_{j+1}$, $\alpha$ the angle $\alpha=\measuredangle OXY$, $a=\|\lambda_{j+1}\|=R+r_{j+1}$,
and $c=\|\mu\|\ge 1$. Then, by the usual angle chasing, $\measuredangle OYZ=\theta_{j}+\alpha$, and denote $$\theta:=\theta_{j+1}-\theta_{j}=\measuredangle ZOY.$$
Figure \ref{fig:angle chasing} illustrates the emerging picture.

We then
apply the sine theorem in the triangle $OYZ$ to yield that
\begin{equation*}
\frac{a}{\sin(\theta_{j}+\alpha)} = \frac{c}{\sin\theta},
\end{equation*}
which, in turn, implies
\begin{equation*}
\sin(\theta_{j}+\alpha)=a\sin(\theta)\cdot \frac{1}{c}.
\end{equation*}
Since $c$ is constant and $a\le R+1/\sqrt{2}$, the above, together with the assumption $\lambda_{j}\in\Sigma_{R,C}$ (equivalent to $\theta\le \frac{C}{R^{2}}$) gives
\begin{equation*}
\sin(\theta_{j}+\alpha) \ll_{\mu,C} \frac{1}{R},
\end{equation*}
hence $$\theta_{j}=-\alpha+O\left(\frac{1}{R}\right),$$ meaning that there exists a sufficiently large $C_0=C_0(\mu,C)$ such that $$\lambda_{j}\in \Gamma_{R}(-\alpha-C_0/R,-\alpha+C_0/R).$$ That \footnote{The number of exceptions is $O_{C}(1)$, much stronger than required.} $\Delta_{R}(\mu)\cap \Sigma_{R,C}$ is of vanishing relative asymptotic density in $\Gamma_{R}$
follows from a straightforward application of Lemma \ref{lem:lat pnt cons rect disc}(i), which, as it was mentioned above, in turn, also yields the assertion of Lemma \ref{lem:tj+1-tj>>1/R^2}.

\end{proof}

\begin{proof}[Proof of Proposition \ref{prop:rect pert arc lat pnts}]

Let $k\ge 1$ be given. In what follows we assume w.l.o.g. that $\lambda_{j}=(R+r_{j})e^{i\theta_{j}}$ with $\theta_j\in [0,\pi/4]$, so that, in particular, $\lambda_{j}$ lies in the $1$st quadrant of $\R^{2}$, and let $\Qc=\Qc_{R,k}(\theta_{j})$ be as in Notation \ref{not:qc}. First, regardless of whatever $\widetilde{\Gamma}_{R}$ is, the number $k'$ of Proposition \ref{prop:rect pert arc lat pnts}(ii) satisfies
$k'\ge k$ for every $j\le K$, by construction of $\Qc_{R,k}(\theta)$, and Lemma \ref{lem:lat pnt cons rect disc}(iii). The bound
$k'\ll k$ with constant involved in the $``\ll"$-notation absolute follows, again, by construction and
Lemma \ref{lem:lat pnt cons rect disc}(iv). Hence \eqref{eq:k<=k'<<k}
will follow automatically once $\widetilde{\Gamma}_{R}$ satisfying all the other properties of Proposition \ref{prop:rect pert arc lat pnts} will be constructed.
Further,
one has that
\begin{equation}
\label{eq:theta'-theta<<1/R}
0\le \theta_{j+k''}-\theta_{j}\le C\frac{k'}{R} \ll \frac{k}{R},
\end{equation}
due to Lemma \ref{lem:lat pnt cons rect disc}(i) and Lemma \ref{lem:lat pnt cons rect disc}(iii).

If $$\lambda \in \Gamma_{R}(\theta_{j},\theta_{j+k'}),$$ then necessarily $\lambda=\lambda_{j+k''}$ with some $1\le k''\le k'$, and if
$\lambda \notin \Qc\cap \Z^{2}$ (i.e. $\lambda \notin \Qc$). Let $\Dc_{R;1/\sqrt{2}}(\theta_{j},\theta_{j+k'})$ be the sector
$$\Dc_{R;1/\sqrt{2}}(\theta_{j},\theta_{j+k'}) = \Dc_{R;1/\sqrt{2}}\cap \left\{x=\|x\|e^{i\theta}:\: \theta \in [\theta_{j},\theta_{j+k'}]\right\},$$
so that $$\Gamma_{R}(\theta_{j},\theta_{j+k'}) = \Dc_{R;1/\sqrt{2}}(\theta_{j},\theta_{j+k'})\cap \Z^{2}.$$

\begin{figure}[h]
  \centering
\includegraphics[height=13cm]{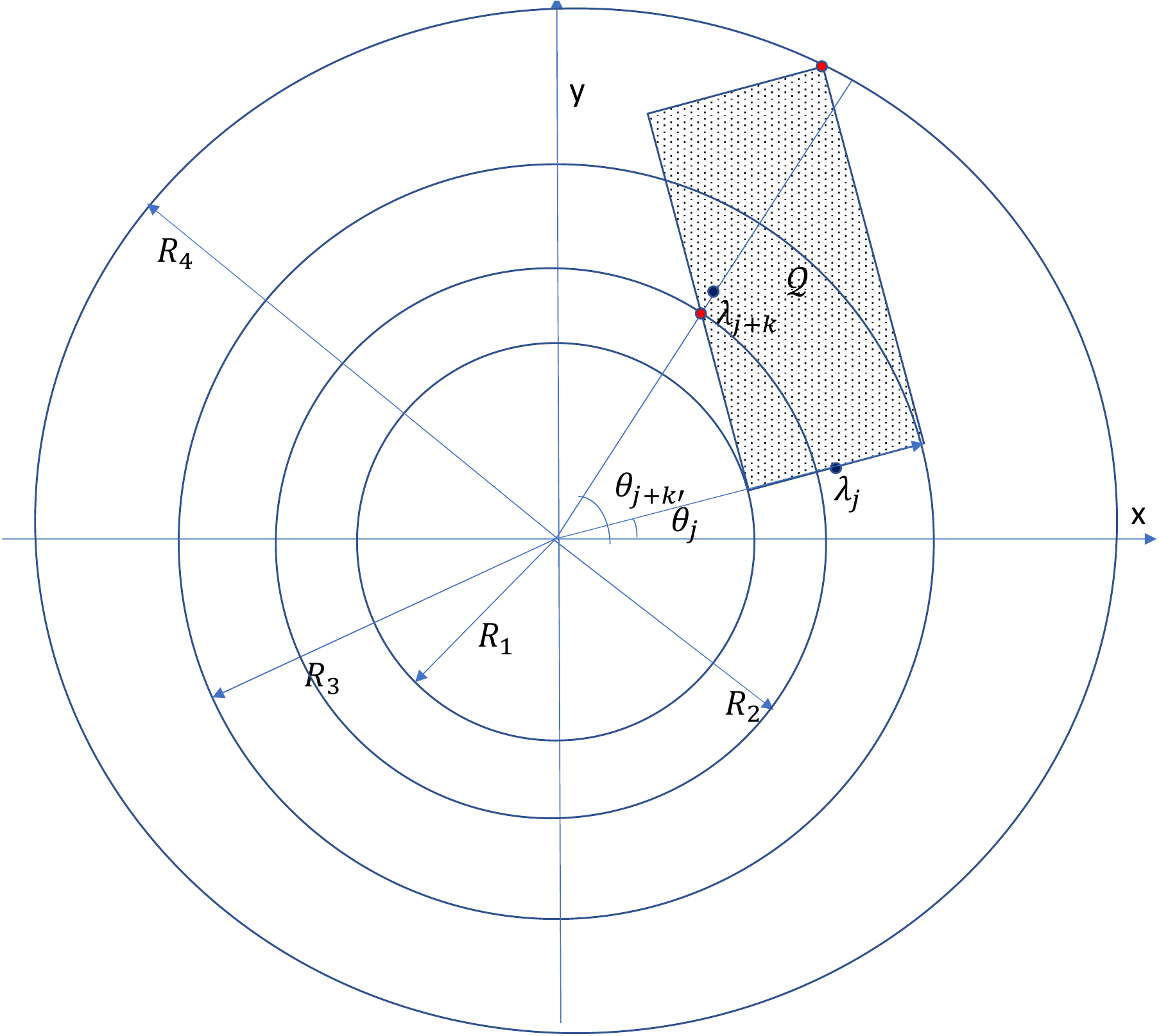}
\caption{The symmetric difference $\Dc_{R;1/\sqrt{2}}(\theta_{j},\theta_{j+k'})\,\triangle\,\Qc$ is contained within a union of two annuli. The red intersections determine the radii $R_{2}$ and $R_{4}$.}
  \label{fig:sec vs rect}
\end{figure}

The emerging picture is illustrated within Figure \ref{fig:sec vs rect}.
Since $$\Dc_{R;1/\sqrt{2}}(\theta_{j},\theta_{j+k'})\setminus \Qc$$ is contained
within the annulus $A_{R_{1},R_{2}}$ (recall Notation \ref{not:annuli, sec}),
with $R_{1}:=R-1/\sqrt{2}$ and
\begin{equation}
\label{eq:R2=R1+O(1/R)}
R_{2}:=\frac{R-1/\sqrt{2}}{\cos(\theta_{j+k'}-\theta_{j})} = R_{1}+O_{k}(1/R),
\end{equation}
by invoking \eqref{eq:theta'-theta<<1/R},
\begin{equation*}
\lambda\in \Gamma_{R}(\theta_{j},\theta_{j+k'})\setminus \Qc \subseteq \Dc_{R;1/\sqrt{2}}(\theta_{j},\theta_{j+k'})\setminus \Qc
\end{equation*}
forces that
$\lambda \in A_{R_{1},R_{2}}$.
For the (somewhat) symmetric case $$\lambda\in (\Qc\cap\Z^{2})\setminus \Gamma_{R}(\theta_{j},\theta_{j+k'}),$$ we have $$\Qc\setminus \Dc_{R;1/\sqrt{2}}(\theta_{j},\theta_{j+k'}) \subseteq A_{R_{3},R_{4}},$$
where $R_{3}:=R+1/\sqrt{2}$, and
\begin{equation}
\label{eq:R4=R3+O(1/R)}
R_{4}:= \sqrt{R_{3}^{2}+(C'k)^{2}}  = R_{3}+O_{k}(1/R),
\end{equation}
by Pythagoras, where $C'k$ is the long side length of the rectangle $\Qc$. It follows that, if $\lambda\in (\Qc\cap\Z^{2})\setminus \Gamma_{R}(\theta_{j},\theta_{j+k'})$, then $\lambda \in A_{R_{3},R_{4}}$.

\vspace{2mm}

Now let $$\Gamma^{0}_{R}:= \{\lambda_{j}\in\Gamma_{R}: \lambda_{j}\in A_{R_{1},R_{2}}\cup A_{R_{3},R_{4}} \} =
\Gamma_{R}\cap (A_{R_{1},R_{2}}\cup A_{R_{3},R_{4}})\subseteq \Gamma_{R} $$ be the set of lattice points of $\Gamma_{R}$ lying inside either of the small annuli $A_{R_{1},R_{2}}$ or
$A_{R_{3},R_{4}}$, with the same indices as $\Gamma_{R}$. Using elementary methods, armed with \eqref{eq:R2=R1+O(1/R)}
or \eqref{eq:R4=R3+O(1/R)},
we may bound the size of $\Gamma^{0}_{R}$ as $$\#\Gamma^{0}_{R} = O\left(R^{2/3}\right) .$$ We define the $k'$-thickening of $\Gamma_R^0$ by
\begin{equation}
\label{eq:Gamma1=Gamma0 k'-thickening}
\begin{split}
\Gamma^{1}_{R}&=\Gamma^{0}_{R}(k'):= \bigcup\limits_{k''=0}^{k'}\left\{\lambda_{j}\in\Gamma_{R}:\: \lambda_{j+k''}\in\Gamma_{R}^{0}\right\} =
 \bigcup\limits_{k''=0}^{k'}\left\{\lambda_{j-k''}\in\Gamma_{R}:\: \lambda_{j}\in\Gamma_{R}^{0}\right\} \\&=
\left\{\lambda_{j}\in\Gamma_{R}:\: \exists 0\le k''\le k' \, \text{ with } \lambda_{j+k''}\in\Gamma_{R}^{0}\right\},
\end{split}
\end{equation}
also of size
\begin{equation}
\label{eq:Gamma1 comp size bnd}
\#\Gamma^{1}_{R} \le (k'+1)\#\Gamma^{0}_{R}\ll k\#\Gamma^{0}_{R} \ll_{k} R^{2/3}=o(R).
\end{equation}

We claim that the complement set
\begin{equation}
\label{eq:Gamma hat excise def}
\widehat{\Gamma}_{R} := \Gamma_{R}\setminus \Gamma^{1}_{R}
\end{equation}
satisfies properties (i)-(ii) of Proposition \ref{prop:rect pert arc lat pnts}. Indeed, first,
(i) of Proposition \ref{prop:rect pert arc lat pnts} is a direct consequence of the size estimate \eqref{eq:Gamma1 comp size bnd}. Next, we check that
(ii) of Proposition \ref{prop:rect pert arc lat pnts} is satisfied. Let $\lambda_{j}\in \Gamma_{R}$ be given, and assume that for some $0\le k''\le k'$, one has that
$$\lambda_{j+k''}\in (\Qc\cap\Z^{2})\, \triangle\, \Gamma_{R}(\theta_{j},\theta_{j+k'}),$$ the symmetric difference of $\Qc\cap\Z^{2}$ and $\Gamma_{R}(\theta_{j},\theta_{j+k'})$. Then, by the given
argument, it would force $$\lambda_{j+k''}\in  A_{R_{1},R_{2}}\cup A_{R_{3},R_{4}},$$ hence $\lambda_{j+k''}\in \Gamma_{R}^{0}$. That would, in turn, imply, by the construction of the excised set
$\Gamma^{1}_{R}$, that $\lambda_{j}\in \Gamma^{1}_{R}$. Hence for $\lambda_{j}\in \widehat{\Gamma}_{R}$,
$$ \Qc\cap\Z^{2} = \Gamma_{R}(\theta_{j},\theta_{j+k'}) = \left\{\lambda_{j},\lambda_{j+1},\,\ldots,\lambda_{j+k'}\right\},$$ as claimed.

\vspace{2mm}

\begin{figure}[h]
  \centering
\includegraphics[height=12cm]{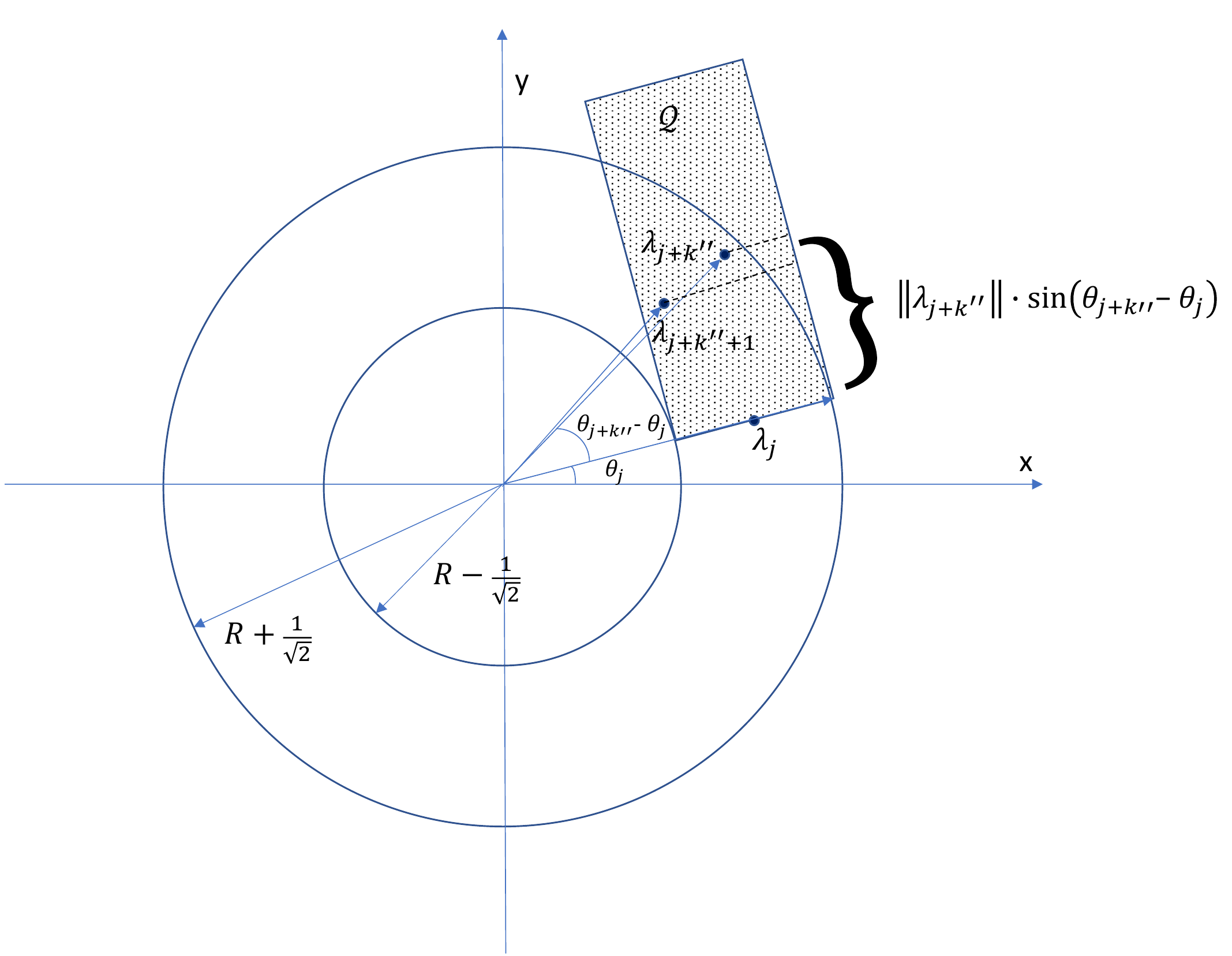}
\caption{An example when orderings w.r.t. the angle and w.r.t. the projection onto the long side of $\Qc$ are different. The projection onto the long side
of $\lambda_{j+k''}$ is $\|\lambda_{j+k''}\|\cdot \sin(\theta_{j+k''}-\theta_{j})$.   }
  \label{fig:order violate projection}
\end{figure}

It then remains to excise a further thin set from $\widehat{\Gamma}_{R}$ to also make it satisfy the property (iii) of Proposition \ref{prop:rect pert arc lat pnts}. The postulated property (iii) of $\widetilde{\Gamma}_{R}$ in Proposition \ref{prop:rect pert arc lat pnts} would be violated, if (and only if) for some $\lambda_{j}\in \widetilde{\Gamma}_{R}$ there exists $0\le k''\le k'-1$ so that the projection of $\lambda_{j+k''}$ onto the long side of $\Qc=\Qc_{R,k}(\theta_{j})$ is greater than that of
$\lambda_{j+k''+1}$, despite the opposite inequality between the angles: $\theta_{j+k''}<\theta_{j+k''+1}$, an example illustrated within Figure \ref{fig:order violate projection}. Since, given $\lambda_{j+k''}$, its projection onto the long side of
$\Qc$ is $$\|\lambda_{j+k''}\|\sin(\theta_{j+k''}-\theta_{j})=(R+r_{j+k''})\cdot \sin(\theta_{j+k''}-\theta_{j}),$$ upon using the shorthand
\begin{equation}
\label{eq:theta def shorthand}
\theta:= \theta_{j+k''}-\theta_{j}=O_{k}(1/R),
\end{equation}
and $$0<\triangle\theta:= \theta_{j+k''+1} - \theta_{j+k''}= (\theta_{j+k''+1}-\theta_{j}) - (\theta_{j+k''}-\theta_{j})=O(1/R),$$
so that $\theta_{j+k''+1}-\theta_{j}=\theta+\triangle\theta$,
this amounts to the inequality
\begin{equation*}
(R+r_{j+k''})\cdot \sin(\theta) > (R+r_{j+k''+1})\cdot \sin(\theta+\triangle\theta).
\end{equation*}
Further, substituting, as we may thanks to \eqref{eq:theta def shorthand}, $$\sin(\theta+\triangle\theta)=\sin{\theta}+\triangle\theta+o_{R\rightarrow\infty}(\triangle\theta),$$
as a result of standard manipulations, we obtain
\begin{equation*}
1+\frac{r_{j+k''}-r_{j+k''+1}}{R+r_{j+k''+1}}=\frac{R+r_{j+k''}}{R+r_{j+k''+1}} > \frac{\sin(\theta+\triangle\theta)}{\sin{\theta}} = 1 + \frac{\triangle\theta}{\sin{\theta}}(1+o(1)),
\end{equation*}
so that
\begin{equation*}
\frac{r_{j+k''}-r_{j+k''+1}}{R+r_{j+k''+1}} >  \frac{\triangle\theta}{\sin{\theta}}(1+o(1)).
\end{equation*}

We may then deduce from the above that
\begin{equation}
\label{eq:r-r>>R^2dt}
r_{j+k''}-r_{j+k''+1} \gg R\cdot \frac{\triangle\theta}{|\sin{\theta}|} \gg_{k} R^{2}\cdot \triangle\theta ,
\end{equation}
thanks to \eqref{eq:theta def shorthand}.
The upshot is that \eqref{eq:r-r>>R^2dt} forces that
\begin{equation}
\triangle\theta = \theta_{j+k''+1} - \theta_{j+k''} \ll_{k} (r_{j+k''}-r_{j+k''+1}) \frac{1}{R^{2}} \ll \frac{1}{R^{2}},
\end{equation}
since $$0<r_{j+k''}-r_{j+k''+1}<2\sqrt{2}$$ is bounded by an absolute constant. Hence
that, in turn, forces that $\theta_{j+k''}$ belongs to the thin set $\Gamma_{R}^{2}\subseteq \Gamma_{R}$
as in Lemma \ref{lem:tj+1-tj>>1/R^2}. Let $$\Gamma_{R}^{3}:= \Gamma_{R}^{2}(k')$$ be its $k'$-thickening as in \eqref{eq:Gamma1=Gamma0 k'-thickening},
that is also thin, and recall that $\widehat{\Gamma}_{R}$ as in \eqref{eq:Gamma hat excise def} is a subset of $\Gamma_{R}$ satisfying
the properties (i)-(ii) of Proposition \ref{prop:rect pert arc lat pnts}.

Then, defining
$$\widetilde{\Gamma}_{R} := \widehat{\Gamma}_{R}\setminus \Gamma_{R}^{3},$$ by further excising $\Gamma_{R}^{3}$ satisfies properties (i)-(iii) of Proposition
\ref{prop:rect pert arc lat pnts}, and we claim that it also satisfies (iv) of Proposition
\ref{prop:rect pert arc lat pnts}. We assume that $\lambda_{j}\in \widetilde{\Gamma}_{R}$, and let $0\le k''\le k'$. The crucial observation is that, since $\Z^{2}$ is invariant w.r.t. shifts by lattice points, one has that
\begin{equation*}
\left(\Qc\cap \Z^{2}\right) -\lambda_{j} = (\Qc-\lambda_{j})\cap \Z^{2} \subseteq \Rc_{r_{j},\theta_{j}+\pi/2}\cap \Z^{2},
\end{equation*}
with $\Rc_{\cdot,\cdot}$ as in \eqref{eq:R semi-inf rectangle}.
Therefore, since by properties (ii) and (iii) of Proposition \ref{prop:rect pert arc lat pnts}, satisfied by the $\widetilde{\Gamma}_{R}$ constructed immediately above, the lattice points of $\Qc$ are those lying inside $\Gamma_{R}$ in the vicinity of $\lambda_{j}$.
We deduce that $$\lambda_{j+k''}=\lambda_{j}+\kappa_{k''},$$ with $\kappa_{k''}=\kappa_{k''}(r_{j},\theta_{j})$ as in \eqref{eq:lambdal in R}.

To bound $$|\Ac_{\infty,k''}(r_{j},\theta_{j})- \Ac_{R}(\lambda_{j+k''})|$$ we rotate the circle so that to put $\lambda_{j}$ on the horizontal axis, so may assume with no loss of generality that $\theta_{j}=0$ in the first place. Then $\Ac_{\infty,k''}(r_{j},\theta_{j})$ is the area of
$$S(\lambda_{j+k''})\cap \{(x_{1},x_{2})\in\R^{2}:\: x_{1}\le R\},$$ whereas $\Ac_{R}(\lambda_{j+k''})$ is the area of
$$S(\lambda_{j+k''})\cap \{(x_{1},x_{2})\in\R^{2}:\: x_{1}^{2}+x_{2}^{2}\le R^{2}\}.$$
Therefore, by the emerging geometric picture (cf. Figure \ref{fig:sec vs rect}
and Figure \ref{fig:order violate projection}),
and taking into account the convexity of the disc,
\begin{equation*}
\begin{split}
0 &< \Ac_{\infty,k''}(r_{j},\theta_{j})- \Ac_{R}(\lambda_{j+k''}) \\&< \area\left(\left\{x=(x_{1},x_{2})\in\R^{2}: R\cos(C\cdot k'/R)\le x_{1}\le R \right\} \cap
\left(\lambda_{j+k''}+ [-1/2,1/2]^{2}\right)\right)\\& = O_{k}\left( \frac{1}{R}\right).
\end{split}
\end{equation*}
Proposition \ref{prop:rect pert arc lat pnts} is proved. \end{proof}

\section{Lattice points in $\Gamma_R$ lying in narrow sectors: Main results and preliminaries}
\label{sec:lattice-intro}

\subsection{Setup}
Let
\begin{equation} \label{eq:PkNdef}
2 \le P \le R^{1-\varepsilon}, 1 \le N \le 100R \quad \text{ and }  \quad 1 \le k \le R^{1-\varepsilon}.
\end{equation}
To streamline our analysis we will use various smooth functions.
\begin{notation} \label{not:smoothfns}
For each $j=1,2$, we take $f_j,g,W_j,W_{R,j},$ and $F_{N,j}$ as follows.
\begin{itemize}
\item $W_j$ is a Schwartz function that is compactly supported on $[-2,2]$ such that $W_j^{(l)}(x) \ll P^l$ for each nonnegative integer $l$.
\item $W_{R,j}(x)=W_j((x-1)R)$.
\item $g,f_2$
are Schwartz functions with $g^{(l)}(x),f_2^{(l)}(x) \ll \frac{P^l}{1+|x|^A}$ for each $A >0$, and nonnegative integer $l$.
\item $f_1(x)=g(x-k)$.
\item $F_{N,j}$ is the $\pi/2$-periodic function
\[
F_{N,j}(\theta)=\sum_{l \in \mathbb Z} f_j\bigg(\frac{N}{\pi/2} \bigg(\theta+l \cdot \frac{\pi}{2} \bigg) \bigg).
\]
\end{itemize}
\end{notation}
The parameter $N$ is our scaling parameter that determines the sector width we consider and $P$ keeps track of the rate of change of our function which will be used to pass from smooth approximations to sharp cutoffs e.g. we will later choose $W_{R,1}(\cdot/R)$ so that it approximates $\chi_{\Dc_{R;1/\sqrt{2}}}$ in $L^1(\mathbb R)$ up to an error of size $O(R/P)$. Also, the parameter $k$ shifts the position of the sector we consider.
Before proceeding let us record the following estimates,
\begin{equation}\label{eq:f1bd}
|\widehat f_1(x)| \ll \bigg(\frac{P}{1+|x|}\bigg)^A, \quad  \int_{\mathbb R} |t^l f_1(t)| \, dt \ll k^l+1,
\end{equation}
for all nonnegative integers $A,l$.
Also for each $j=1,2$ we have by Poisson summation that
\begin{equation} \label{eq:fourierseries}
F_{N,j}(\theta)=\frac{1}{N} \sum_{\ell \in \mathbb Z} \widehat f_j\bigg( \frac{\ell}{N} \bigg) e^{4i \theta \ell}.
\end{equation}

\subsection{Main results}
Theorems \ref{thm:jointdistribution} and \ref{thm:varbd}
will follow from the following propositions. To state our main result let us introduce some further notation. Given $b \in \mathbb Z$, $t \in \mathbb R$, nonzero $\lambda \in \mathbb Z^2$ let
\begin{equation} \label{eq:hdef}
h(t,\lambda,b)=2\pi(\tfrac{t}{4}+\tfrac{b}{4})+\theta_{\lambda}.
\end{equation}

\begin{proposition} \label{thm:paircorrelation}
Let $\varepsilon>0$.
Suppose that $P^{5/2} k \le R^{1/2-\delta_0}$ where $\delta_0>0$ is fixed. Then
\begin{equation} \label{eq:angles}
\begin{split}
&\frac{1}{\pazocal K}\sum_{\lambda,\mu \in \mathbb Z^2}  F_{\pazocal K,1}(\theta_{\lambda}-\theta_{\mu}) F_{1,2}(\theta_{\lambda}) W_{R,1}\bigg( \frac{\|\lambda\|^2}{R^2}\bigg)W_{R,2}\bigg( \frac{\|\mu\|^2}{R^2}\bigg)\\
&\qquad  =\frac{1}{8} \widehat W_1(0) \widehat W_2(0) \widehat f_1(0) \widehat f_2(0)+\mathcal M(W_1,W_2,f_1,f_2)+O\bigg(\frac{P^{7/2}}{ R^{1/2-\varepsilon}}+\frac{k^2 P^5}{R^{1-\varepsilon}} \bigg)
\end{split}
\end{equation}
where
\begin{equation} \label{eq:Mdef}
\begin{split}
&\mathcal M(W_1,W_2,f_1,f_2) =\frac{1}{32} \sum_{0 \le b < 4} \sum_{\substack{\lambda \in \mathbb Z^2 \setminus\{0\} }} \int_{\mathbb R} \widehat W_1\bigg( -\frac{\|\lambda\|}{2} \sin(h(t,\lambda,b))\bigg)\\
&\qquad \qquad \times\widehat W_2\bigg( \frac{\|\lambda\|}{2} \sin(h(t,\lambda,b))\bigg)  \widehat f_1\bigg(\frac{\pi }{2\sqrt{2}} \|\lambda\|  \cos(h(t,\lambda,b)) \bigg) f_2(-t) \, dt.
\end{split}
\end{equation}
\end{proposition}

Roughly, the sum on the l.h.s. of \eqref{eq:angles} counts the number of pairs of lattice points $\lambda,\mu$
such that the difference between their corresponding angles lies in an interval in $\mathbb R/\mathbb Z$ of length $\asymp 1/\pazocal K$. For $\lambda \in \Gamma_R$ recall $\lambda=(R+r_{\lambda})e^{i\theta_{\lambda}}$.
The first term on the r.h.s. of \eqref{eq:angles} corresponds to the joint equidistribution of $(r_{\lambda},r_{\mu},\frac{2}{\pi}\theta_{\lambda})$ within $[\tfrac{-1}{\sqrt{2}},\tfrac{1}{\sqrt{2}}] \times [\tfrac{-1}{\sqrt{2}},\tfrac{1}{\sqrt{2}}] \times \mathbb R/\mathbb Z$, over such pairs of lattice points $\lambda,\mu$.  The secondary main term $\mathcal M$ does not depend on $R$ provided $P,k$ do not depend on $R$ and quantifies how $r_{\mu}$ and $r_{\lambda}$ depend on one another.

\vspace{2mm}

\begin{notation} \label{not:tildenotation} Given $\theta \in \mathbb R$ we let $\widetilde \theta=\tfrac{2}{\pi} \theta$.
\end{notation}
The following formula for the variance of the number of lattice points lying in a narrow sector is established next:

\begin{proposition} \label{thm:varasymp}
Let $\varepsilon>0$.
Let $c,d \in \mathbb R$ with $c<d$. Also, let $N=|d-c|^{-1}$. Suppose that $ P^{5/2} R^{1/2+\delta_1} \le N \le 100 R$ for some fixed $\delta_1>0$. Then uniformly for $I \subseteq [-\tfrac{1}{\sqrt{2}}, \tfrac{1}{\sqrt{2}}]$ we have that
\begin{equation} \label{eq:varform}
\begin{split}
&\int_0^{\pi/2} \bigg| \sum_{\substack{ \lambda \in \Gamma_R, \, r_{\lambda} \in I \\
(\widetilde \theta_{\lambda}-\widetilde \theta) \Mod 1 \in [c, d ]}} 1-2\pi \sqrt{2} R   |d-c|
\frac{|I|}{ \sqrt{2}}  \bigg|^2 \frac{d\theta}{\pi/2}\\
&\qquad  \qquad =4 \pi^2 R^2 |d-c|^2\sum_{1 \le \| \lambda \| \le P^{1+\varepsilon}} \int_{-1/2}^{1/2} \bigg| \widehat \chi_{I}( \|\lambda\| \sin(2\pi t) )
\widehat \chi_{[cN,dN]}( \tfrac{\pi^2 R}{N} \|\lambda\| \cos(2\pi t) ) \bigg|^2 \, dt\\
&\qquad \qquad\qquad \qquad \qquad+ O\bigg( \frac{R^2}{N^2 P^{1-\varepsilon}}+\frac{R^{2+\varepsilon}P^4}{N^2}  \bigg( \frac{RP}{N^2}+\frac{1}{R^{1/2}}\bigg)\bigg).
\end{split}
\end{equation}
\end{proposition}
In the range $N \asymp R$ we can take  $P=R^{1/10}$ and the
error term is $\ll R^{-1/10+\varepsilon}$.

\subsection{Smooth approximations}
In the proofs of theorems \ref{thm:varbd} and \ref{thm:jointdistribution} and Proposition \ref{prop:equid}
we will approximate the indicator functions of intervals in terms of Schwartz functions. In this section we will describe the constructions of the smooth approximations that will be employed in the proofs of these results.

For $j=1,2$ and $I_j \subseteq [-\tfrac{1}{\sqrt{2}}, \tfrac{1}{\sqrt{2}}]$ with $I_j=[a_j,b_j]$ we take $W_j^{+}$ to be a Schwartz function satisfying
\[
W_j^{+}(x)=
\begin{cases}
1  & \text{ if } 2a_j-\frac{1}{P} \le x \le 2b_j+\frac{1}{P}, \\
0 & \text{ if } x \le 2a_j-\frac{2}{P}  \text{ or } x \ge 2b_j+\frac{2}{P},
\end{cases}
\]
$0 \le W_j^{+}(x) \le 1$
and $(W_j^{+})^{(l)}(x) \ll P^l$ uniformly with respect to $I_j$.
Similarly, we take $W_j^{-}$ to be a Schwartz function satisfying
\[
W_j^{-}(x)=
\begin{cases}
1  & \text{ if } 2a_j+\frac{2}{P} \le x \le 2b_j-\frac{2}{P}, \\
0 & \text{ if } x \le 2a_j+\frac{1}{P}  \text{ or } x \ge 2b_j-\frac{1}{P},
\end{cases}
\]
 $0 \le W_j^{-}(x) \le 1$ and $(W_j^{-})^{(l)}(x) \ll P^l$ uniformly with respect to $I_j$.
We also take $W_{R,j}^{\pm}(x)=W_j^{\pm}((x-1)R)$. Hence, for $\lambda \in \mathbb Z^2$ writing $\|\lambda\|=R+r_{\lambda}$
we have that
\begin{equation} \label{eq:Wbd}
0 \le W_{R,j}^{-}\left( \frac{\|\lambda\|^2}{R^2}\right) \le \chi_{I_j}(r_{\lambda}) \le W_{R,j}^{+}\left( \frac{\|\lambda\|^2}{R^2}\right), \quad \widehat W_{j}^{\pm}(0)=2|I_j|+O(P^{-1}).
\end{equation}
Also take $g^{\pm}$ to be Schwartz functions with
\begin{equation} \label{eq:f1bd2}
0 \le g^{-}(x) \le \chi_{[0,1]}(x) \le g^+(x), \quad \widehat g^{\pm}(0) =1+O(1/P), \quad (g^{\pm})^{(l)} (x)\ll \frac{P^l}{(1+|x|)^A},
\end{equation}
for all nonnegative integers $A,l$.
Given $J=[c,d) \subseteq \mathbb R$ with $|J|\le 1$
we take
$f_2^{\pm}$ to be Schwartz functions with
\begin{equation} \label{eq:f2bd}
0 \le f_2^{-}(x) \le \chi_{J}(x) \le f_2^+(x), \quad \widehat f_2^{\pm}(0) =|J|+O(1/P), \quad (f_2^{\pm})^{(l)} (x)\ll \frac{P^l}{(1+|x|)^A},
\end{equation}
for all nonnegative integers $A,l$.
Let $N=|d-c|^{-1}$. Also, let
\begin{equation} \label{eq:jpdef}
c'=cN, \quad d'=dN \quad \text{ and } \quad  J'=[c',d').
\end{equation}
We take
$f_3^{\pm}$ to be Schwartz functions with
\begin{equation} \label{eq:f3bd}
0 \le f_3^{-}(x) \le \chi_{J'}(x) \le f_3^+(x), \quad \widehat f_3^{\pm}(0) =1+O(1/P), \quad (f_3^{\pm})^{(l)} (x)\ll \frac{P^l}{(1+|x|)^A},
\end{equation}
for all nonnegative integers $A,l$.
Also, let $f_1^{\pm}(x)=g^{\pm}(x-k)$ and
\begin{equation}\label{eq:Fdef}
\begin{split}
F_{\pazocal K,1}^{\pm}(\theta) =&\sum_{j \in \mathbb Z} g^{\pm}\bigg(\frac{\pazocal K}{\pi/2}\bigg(\theta-\frac{\pi}{2} \cdot \frac{k}{\pazocal K}+j\cdot \frac{\pi}{2} \bigg) \bigg) \\
=&\sum_{j \in \mathbb Z} f_1^{\pm}\bigg(\frac{\pazocal K}{\pi/2}\bigg(\theta+j\cdot \frac{\pi}{2} \bigg) \bigg).
\end{split}
\end{equation}
We have by \eqref{eq:f1bd2} that
\begin{equation}\label{eq:Fkineq}
F_{\pazocal K,1}^{-}(\theta) \le \chi_{[\frac{k}{\pazocal K}, \frac{k+1}{\pazocal K}]}( \tfrac{2}{\pi} \theta \Mod 1) \le F_{\pazocal K,1}^{+}(\theta).
\end{equation}
 Periodicizing $f_2^{\pm}$ we define
\begin{equation}\label{eq:F2ineq}
F_{1,2}^{\pm}(\theta)=\sum_{j \in \mathbb Z} f_{2}^{\pm}(\tfrac{2}{\pi} \theta+j) \quad \text{ so that } \quad
F_{1,2}^{-}(\theta) \le \chi_{J}( \tfrac{2}{\pi}  \theta \Mod 1) \le F_{1,2}^{+}(\theta),
\end{equation}
by \eqref{eq:f2bd}. By \eqref{eq:jpdef} we similarly have
\begin{equation}\label{eq:F3ineq}
F_{N,3}^{\pm}(\theta)=\sum_{j \in \mathbb Z} f_3^{\pm}\bigg(\frac{N}{\pi/2}\bigg(\theta+j\cdot \frac{\pi}{2} \bigg) \bigg) \quad \text{ so that } \quad
F_{N,3}^{-}(\theta) \le \chi_{J}( \tfrac{2}{\pi}  \theta \Mod 1) \le F_{N,3}^{+}(\theta).
\end{equation}

\subsection{Proof of Theorem \ref{thm:varbd}}

\begin{proof}[Proof of Theorem \ref{thm:varbd} assuming Proposition \ref{thm:varasymp}]
Write $M=R/N$ where $N=|d-c|^{-1}$. Also, recall $c'=cN$, $d'=dN$ and $J'=[c',d')$ as in \eqref{eq:jpdef}.
Write $\alpha=9/10+\delta_2$, since $N \ge R^{\alpha}$,  by the assumptions of Theorem \ref{thm:varbd}, and $MN=R$ we have that
\begin{equation} \label{eq:Nbd}
N \ge M^{\frac{\alpha-\delta_2/10}{1-\alpha}}  R^{\delta_2/10} \ge M^{9+9\delta_2}   R^{\delta_2/10}.
\end{equation}
In Proposition \ref{thm:varasymp}
we take $P=(N/M)^{1/8} R^{-\delta_2/100} M^{-9\delta_2/40}$
and note that we have $P^{5/2} R^{1/2+\delta_1} \le R^{13/16+\delta_1}\le N,$ as needed. Let us also note that by construction and recalling \eqref{eq:Nbd}
we have that $P \ge M^{1+9 \delta_2/10} R^{\delta_2/400}$.
Hence, with these choices and recalling that $N \ge R^{9/10+\delta_2}$ the error term in
Proposition \ref{thm:varasymp} is
\[
\begin{split}
\ll& \frac{R^2}{N^2 P^{1-\varepsilon}}+\frac{R^{2+\varepsilon}}{N^2} \bigg( P^4 \bigg( \frac{RP}{N^2}+\frac{1}{R^{1/2}}\bigg)\bigg) \ll \frac{R^2}{N^2} \bigg( \frac{1}{P^{1-\varepsilon}}+\frac{P^4}{R^{1/2-\varepsilon}} \bigg)\\
=& M^2\bigg(
\frac{1}{P^{1-\varepsilon}}+\frac{N^{1/2}}{M^{1/2+9\delta_2/10} R^{\delta_2/25}} \cdot \frac{R^{\varepsilon}}{ \sqrt{MN}} \bigg) \ll M^{1-9\delta_2/10}.
\end{split}
\]
To complete the proof, it suffices to show that
for $n=\|\lambda\|^2 \in \mathbb N$ and $M \ge 1$ we have that
\begin{equation} \label{eq:asymp}
 \int_{-1/2}^{1/2} | \widehat \chi_I( \sqrt{n}\sin(2\pi t) )
\widehat \chi_{J'}( \pi^2 M \sqrt{n} \cos(2\pi t) ) |^2 dt
=\frac{|\widehat \chi_I(\sqrt{n})|^2}{\pi^3 M \sqrt{n}}
+O\bigg( \frac{1}{M^{6/5} n (\log(2n))^2}\bigg),
\end{equation}
(note that $\sum_{n \ge 1} (n (\log(2n))^2)^{-1}$ is convergent).

Let $\Delta=M^{2/5} \sqrt{n} (\log(2n))^{-2}$.
We first record the following simple estimate
\begin{equation}
\label{eq:chi hat simple}
\widehat \chi_{[a,b]}(\xi) \ll \min( |\xi|^{-1},|b-a|).
\end{equation}
For $|t\pm 1/4| \ge 1/\Delta$ and $|t| \le 1/2$ we have $|M \sqrt{n} \cos(2\pi t) | \gg M \sqrt{n}/\Delta$, so that
using \eqref{eq:chi hat simple} for $\widehat \chi_{[a,b]}$ we get that
\begin{equation} \label{eq:firstbd}
\begin{split}
&\int_{\substack{ |t\pm \frac14| \ge \frac{1}{\Delta} \\ \& \, |t| \le \frac12}} |\widehat \chi_{I}(\sqrt{n} \sin(2\pi t)) \widehat \chi_{J'}(\pi^2 M \sqrt{n} \cos(2\pi t))|^2 \, dt  \\
&\ll \int_{|t| \le 1/(\log(2n))^2}
\frac{1}{nM^2} \, dt+\int_{1/(\log(2n))^2}^1 \frac{1}{nt^2} \cdot \frac{\Delta^2}{M^2 n} \, dt \ll \frac{1}{M^2 n (\log(2n))^2}+\frac{\Delta^2 (\log(2n))^2}{M^2 n^2}.
\end{split}
\end{equation}

For $|h| \le |\xi|/2$ we have that
\begin{equation} \label{eq:funbd}
 \qquad \widehat \chi_{[a,b]}(\xi+h)=e(-\tfrac{(a+b)}{2}h)\widehat \chi_{[a,b]}(\xi)+O(|h| \min(|\xi|^{-1},|\xi|))
\end{equation}
where the implied constant depends on $|b-a|$. Using this estimate along with Taylor expansions
for sine and cosine, we have that
\begin{equation} \label{eq:taylorasymp}
\begin{split}
&\int_{|t\pm \frac14| \le \frac{1}{\Delta}}
\bigg| \widehat \chi_I( \sqrt{n}\sin(2\pi t) )
\widehat \chi_{J'}\bigg( \pi^2 M \sqrt{n} \cos(2\pi t) \bigg) \bigg|^2 dt
\\
&\qquad \qquad =2 \int_{|t|\le \frac{1}{\Delta}} |\widehat \chi_{I}( \sqrt{n})+O(t^2)|^2|\widehat \chi_{J'}(2\pi^3 M \sqrt{n} t)+O(t^2)|^2 dt \\
&\qquad \qquad =2 |\widehat \chi_{I}( \sqrt{n})|^2 \int_{|t|\le \frac{1}{\Delta}} |\widehat \chi_{J'}(2\pi^3 M \sqrt{n} t)|^2 dt +O\bigg( \frac{1}{\sqrt{n} \Delta^3}+\frac{1}{M\sqrt{n}\Delta^2}+\frac{1}{\Delta^5}\bigg).
\end{split}
\end{equation}
Also, making a linear transformation of variables then extending the integral to all of $\mathbb R$ we have that
\begin{equation}
\label{eq:int chi Plancherel}
2\int_{|t|\le \frac{1}{\Delta}} |\widehat \chi_{J'}(2\pi^3 M \sqrt{n} t)|^2 dt =
\frac{1}{\pi^3 M\sqrt{n}} \bigg( \int_{\mathbb R} |\widehat \chi_{J'}( t)|^2 \, dt
+O\bigg(\frac{\Delta}{M \sqrt{n}} \bigg) \bigg).
\end{equation}
By Plancherel's theorem
the integral on the r.h.s. of \eqref{eq:int chi Plancherel} equals $\int_{\mathbb R} \chi_{J'}(t)^2 \, dt=|J'|=1$. Applying this estimate in \eqref{eq:taylorasymp} and combining the resulting formula along with \eqref{eq:firstbd} yields \eqref{eq:asymp}.
\end{proof}

\subsection{Proof of Theorem \ref{thm:jointdistribution}}
Let us state the following result, which follows from the method of
stationary phase, whose proof will be given in
Section \ref{sec:lemproof}. Recall $h(t,\lambda,b)=2\pi(t+\tfrac{b}{4})+\theta_{\lambda}$, as in \eqref{eq:hdef}.

\begin{lemma} \label{lem:stationaryphaseapplied2}
Let $\varepsilon>0$.
Let $\delta_1>0$ be a fixed number, and suppose that $k\ge P^{3+3\delta_1}$.
For $\lambda \in \mathbb Z^2\setminus\{0\}$ denote
\[
\begin{split}
&I:=\int_{\mathbb R} \widehat W_1\bigg(-\frac{\|\lambda\|}{2} \sin(h(t,\lambda,b))\bigg)\widehat W_2\bigg( \frac{\|\lambda\|}{2} \sin(h(t,\lambda,b))\bigg)   \\
& \qquad \qquad \times \widehat f_1\bigg(\frac{\pi}{2\sqrt{2}} \|\lambda\| \cos(h(t,\lambda,b)) \bigg) f_2(-t) \, dt.
\end{split}
\]
Then, uniformly w.r.t. $1 \le \|\lambda\| \le k^{\varepsilon}P$, one has that
\begin{equation*}
I = O\bigg( \frac{1}{k^{1/2-\varepsilon} \sqrt{\|\lambda\|}}\bigg).
\end{equation*}
\end{lemma}

\begin{proof}[Proof of Theorem \ref{thm:jointdistribution} assuming
Proposition \ref{thm:paircorrelation} and Lemma \ref{lem:stationaryphaseapplied2}]
Recall \eqref{eq:Wbd}, \eqref{eq:Fdef}, \eqref{eq:Fkineq}, and \eqref{eq:F2ineq}. Also, recall Notation \ref{not:smoothfns}.
Applying Proposition \ref{thm:paircorrelation} with $W_{R,j}=W_{R,j}^{+}$ for $j=1,2$, $F_{1,2}=F_{1,2}^+$ and $F_{\pazocal K,1}=F_{\pazocal K,1}^{+}$
we have for $P^{5/2}k \le R^{1/2-\delta_0}$ that
\begin{equation} \label{eq:countbd}
\begin{split}
&\frac{1}{\pazocal K} \# \bigg\{ \lambda,\mu \in \Gamma_R :
r_{\lambda} \in I_1, r_{\mu} \in I_2,
\tfrac{2}{\pi}  \theta_{\lambda} \Mod 1 \in J, \, \& \,
\tfrac{2}{\pi} ( \theta_{\lambda}-\theta_{\mu}) \Mod 1 \in \bigg[ \frac{k}{\pazocal K}, \frac{k+1}{\pazocal K}\bigg]\bigg\} \\
&\le \frac{1}{\pazocal K}
\sum_{\lambda,\mu \in \mathbb Z^2}  F_{\pazocal K,1}^{+}(\theta_{\lambda}-\theta_{\mu}) F_{1,2}^{+}(\theta_{\lambda}) W_{R,1}^{+}\bigg( \frac{\|\lambda\|^2}{R^2}\bigg)W_{R,2}^+\bigg( \frac{\|\mu\|^2}{R^2}\bigg)\\
& =\frac{1}{8} \widehat W_1^{+}(0) \widehat W_2^{+}(0) \widehat f_1^{+}(0) \widehat f_2^{+}(0)+\mathcal M(W_1^{+},W_2^{+},f_1^{+},f_2^{+})+O\bigg(\frac{P^{7/2}}{ R^{1/2-\varepsilon}}+\frac{k^2 P^5}{R^{1-\varepsilon}} \bigg).
\end{split}
\end{equation}
Since $\widehat f_1^{+}(0)=\widehat g^{+}(0)$, we can apply \eqref{eq:Wbd}, \eqref{eq:f1bd2}, and \eqref{eq:f2bd} to get that
\begin{equation} \label{eq:MT}
 \frac{1}{8} \widehat W_1^{+}(0) \widehat W_2^{+}(0) \widehat f_1^{+}(0) \widehat f_2^{+}(0)=\frac{|I_1|}{\sqrt{2}}\frac{|I_2|}{\sqrt{2}}|J|+O(P^{-1})
\end{equation}
uniformly in $I_1,I_2,J$ (recall these are bounded intervals).
Note that by repeatedly integrating by parts and noting $|\widehat f_1^{\pm}|=|\widehat g^{\pm}|$ we have for any integer $A \ge 0$ that
\[
\widehat W_{j}^{\pm}(x),\widehat f_1^{\pm}(x) \ll \bigg(\frac{P}{1+|x|} \bigg)^A.
\]
Applying Lemma \ref{lem:stationaryphaseapplied2} and using the above estimates
we have for $P^{3+3\delta_1} \le k$ that
\begin{equation} \label{eq:ET}
\mathcal M(W_1^{+},W_2^{+},f_1^{+},f_2^{+})\ll \frac{1}{k^{1/2-\varepsilon}} \sum_{1 \le \|\lambda\| \le k^{\varepsilon}P} \frac{1}{\sqrt{\|\lambda\|}} \ll \frac{P^{3/2}}{k^{1/2-\varepsilon}}.
\end{equation}

Write $P=k^{\eta}$
and assume that $\eta<1/(3+3\delta_1)$. Then, with this choice for $P$, substituting the estimates \eqref{eq:ET} and \eqref{eq:MT}
into \eqref{eq:countbd}, one has the inequality
\begin{equation} \label{eq:ebd}
\begin{split}
\frac{1}{\pazocal K} &\# \bigg\{ \lambda,\mu \in \Gamma_R :
r_{\lambda} \in I_1, r_{\mu} \in I_2,
\tfrac{2}{\pi}  \theta_{\lambda} \Mod 1 \in J, \, \& \,
\tfrac{2}{\pi} ( \theta_{\lambda}-\theta_{\mu}) \Mod 1 \in \bigg[ \frac{k}{\pazocal K}, \frac{k+1}{\pazocal K}\bigg]\bigg\}
\\&\le \frac{|I_1|}{\sqrt{2}}\frac{|I_2|}{\sqrt{2}}|J|+ O\left(\frac{1}{k^{\eta-\varepsilon}}+\frac{1}{k^{\frac12-\frac{3\eta}{2} -\varepsilon}}+\frac{k^{2+5\eta}}{R^{1-\varepsilon}}+\frac{k^{7\eta/2}}{R^{1/2-\varepsilon}} \right).
\end{split}
\end{equation}
Arguing along the same lines we may obtain an analogous lower bound
\begin{equation}
\label{eq:ebdlow}
\ge \frac{|I_1|}{\sqrt{2}}\frac{|I_2|}{\sqrt{2}}|J|- O\left(\frac{1}{k^{\eta-\varepsilon}}+\frac{1}{k^{\frac12-\frac{3\eta}{2} -\varepsilon}}+\frac{k^{2+5\eta}}{R^{1-\varepsilon}}+\frac{k^{7\eta/2}}{R^{1/2-\varepsilon}} \right)
\end{equation}
for the l.h.s. of \eqref{eq:ebd}. Combine the lower and upper inequalities \eqref{eq:ebdlow} and \eqref{eq:ebd} for the l.h.s. of
\eqref{eq:ebd}. Upon noting in the l.h.s. of \eqref{eq:ebd} that by symmetry we may pass to counting angles $\Mod{\tfrac{\pi}{2}}$ to counting angles $\Mod{2\pi}$, rescaling both intervals $J$, $[k/\pazocal K,(k+1)/\pazocal K]$ and the parameter $k$ each by a factor of $\pi/2$, we obtain the estimate \eqref{eq:count},
hence completing the proof of Theorem \ref{thm:jointdistribution}.
\end{proof}

\section{Preliminary estimates and the proof of Proposition \ref{prop:equid}} \label{sec:lattice-proof1}

\subsection{ A summation formula}
\begin{notation} \label{not:hankel}
Given a Schwartz function $g:\mathbb R_{>0}\rightarrow \mathbb C$ we define the Hankel type transform
\[
\mathcal B_{l}(g)(\xi)=\int_{0}^{\infty} g(y) J_{l}(2\pi \sqrt{\xi y}) \, dy.
\]
\end{notation}
For $n \in \mathbb N$ and $\ell \in \mathbb Z$ we also write
\[
\tau_{4\ell}(n)=\sum_{\substack{\lambda \in \mathbb Z^2 \\ \|\lambda\|^2=n} } e^{4i\ell \theta_{\lambda}}.
\]
We note that $\tau_{-4\ell}(n)=\tau_{4\ell}(n)$ and $|\tau_{4\ell}(n)| \le \tau_0(n) \ll n^{\varepsilon}$.

\begin{lemma} \label{lem:voronoi}
Let $\ell \in \mathbb Z$. Suppose $F$ is a Schwartz function supported on the positive real numbers. Then for $X>0$ we have that
\begin{equation} \label{eq:sum}
\sum_{n \ge 1} \tau_{4\ell} (n) F\bigg( \frac{n}{X} \bigg)=\delta_{\ell=0} \cdot \pi X \int_0^{\infty} F(r)\,dr+ \pi X \sum_{n \ge 1} \tau_{4\ell}(n) \mathcal B_{4\ell}( F)(nX),
\end{equation}
where $\delta_{\ell=0}$ is one if $\ell=0$ and is zero otherwise.
\end{lemma}
\begin{proof}
Let $\vartheta: \mathbb R^2\setminus\{0\} \rightarrow [-2\pi,2\pi]$ be given by
\[
\vartheta(u,v)=\begin{cases}
4 \arctan(v/u) \, & \text{ if } \, u \neq 0,\\
0 \, & \text{ if } \, u=0.
\end{cases}
\]
For any non-zero $\lambda=(a,b) \in \mathbb Z^2$ note that $e^{i4\theta_{\lambda}}=e^{i\vartheta(a,b)}$. The l.h.s. of \eqref{eq:sum} is
\[
=\sum_{(a,b) \in \mathbb Z^2} e^{i\ell\vartheta(a,b)} F\bigg(\frac{a^2+b^2}{X} \bigg).
\]
We now apply Poisson summation to see the preceding expression is
\begin{equation} \label{eq:poisson}
= \sum_{(a,b) \in \mathbb Z^2} \int_{\mathbb R^2} e^{i \ell \vartheta(u,v) }  F\bigg(\frac{u^2+v^2}{X} \bigg) e(-au-bv) \, dudv.
\end{equation}
We treat the term with $(a,b)=0$ separately and use polar coordinates to see that it equals
\begin{equation} \label{eq:mtpoisson}
\begin{split}
\int_{0}^{\infty} \int_0^{2\pi} e^{4i \theta \ell}  F\bigg(\frac{r^2}{X} \bigg) r d\theta dr&=\frac{X}{2} \int_{0}^{\infty} F(r) \bigg( \int_0^{2\pi}  e^{4i \theta \ell}  d\theta\bigg)  dr
=\delta_{\ell=0} \cdot \pi X \int_0^{\infty} F(r)\,dr.
\end{split}
\end{equation}
The terms with $\lambda=(a,b) \in \mathbb Z^2\setminus\{0\}$ in \eqref{eq:poisson} are
\begin{equation} \label{eq:frequencies1}
\begin{split}
&= \sum_{\lambda \in \mathbb Z^2\setminus\{0\}} \int_{0}^{\infty} \int_0^{2\pi} e^{4i \theta \ell}  F\bigg(\frac{r^2}{X} \bigg) e(-\|\lambda\| r (\cos(\theta)\cos(\theta_{\lambda})+\sin(\theta)\sin(\theta_{\lambda}))) \, r d\theta dr\\
&=\pi X  \sum_{\lambda \in \mathbb Z^2\setminus\{0\}} \int_{0}^{\infty}   F(r ) \bigg( \int_0^{2\pi} e^{4i \theta \ell}   e(-\|\lambda\| \sqrt{Xr}(\cos(\theta-\theta_{\lambda}))) \frac{d\theta}{2\pi} \bigg) dr.
\end{split}
\end{equation}
Also
\begin{equation} \label{eq:frequencies2}
\int_0^{2\pi} e^{4i \theta \ell}   e(-\|\lambda\| \sqrt{Xr}(\cos(\theta-\theta_{\lambda}))) \frac{d\theta}{2\pi}=e^{4i\ell \theta_{\lambda}}J_{4\ell}(2\pi \|\lambda\| \sqrt{Xr}),
\end{equation}
and note \cite[8.404.2,8.411.1]{GR} that $J_{4\ell}(\cdot)=J_{-4\ell}(\cdot)$.
Combining \eqref{eq:frequencies1} and \eqref{eq:frequencies2} then using the resulting formula with \eqref{eq:mtpoisson} in \eqref{eq:poisson} we complete the proof.
\end{proof}

To quantify how the summation formula \eqref{eq:sum} transforms our sum let us record the following estimate, that essentially follows from integrating by parts.

\begin{lemma} \label{lem:IBP}
Let $l$ be an integer, and let $\mathcal B_{l}( W_{R,j})$ be the Hankel type transform of $W_{R,j}$ as in Notation \ref{not:hankel}.
Then, for $j=1,2$ and any integer $A \ge 1$, and real number $\xi\ge 1$, we have
\[
\mathcal B_{l}( W_{R,j})(\xi) \ll  \bigg(\frac{|l|+RP}{\sqrt{\xi}} \bigg)^A \frac{1}{R}.
\]
Additionally, if we also have that $\xi \ge 10 l^2$ and $ \int_{\mathbb R} |W_{j}^{(m)}(x)| dx \ll P^{m-1}$ for each integer $m \ge 1$ then
\[
\mathcal B_{l}( W_{R,j})(\xi) \ll
\frac{|l|^A+R^A P^{A-1}}{R\xi^{A/2+1/4}} .
\]
\end{lemma}

Our range of interest will be $|l| \ll \pazocal K \asymp R$ and the lemma implies that the dual sum after applying Poisson summation is essentially supported on integers $n$ with $n \le P^2  R^{\varepsilon}$. Coincidentally, this range is at the transition regime for the Bessel function.

\begin{proof}
By \cite[8.472.3]{GR} we have
\[
\frac{d}{dx} ( x^{v+1}J_{v+1}(x))=x^{v+1}J_v(x).
\]
For brevity write $W_R=W_{R,j}$.
We integrate by parts to get
\begin{equation} \label{eq:IBP1}
\begin{split}
\int_{0}^{\infty} W_R(y) J_{l}( \sqrt{\xi y}) \, dy
=&\frac{2}{\xi} \int_0^{\infty} y^{-l}W_R\bigg( \frac{y^2}{\xi}\bigg) d(y^{l+1} J_{l+1}(y))\\
=&\frac{2}{\xi} \int_0^{\infty} \bigg( \frac{2y^2}{\xi}  W_R'\bigg( \frac{y^2}{\xi}\bigg)
-l W_R\bigg( \frac{y^2}{\xi} \bigg)\bigg) J_{l+1}(y) \, dy.
\end{split}
\end{equation}
Integrating by parts once more, we obtain
\begin{equation}  \notag
\begin{split}
\int_0^{\infty} \frac{y^2}{\xi}  W_R'\bigg( \frac{y^2}{\xi}\bigg)J_{l+1}(y) \, dy&=\int_0^{\infty} \frac{y^{-l}}{\xi}  W_R'\bigg( \frac{y^2}{\xi}\bigg)d(y^{l+2}J_{l+2}(y))  \\
&=\int_0^{\infty} \frac{1}{y} \bigg( \frac{2y^4}{\xi^2}  W_R^{''}\bigg( \frac{y^2}{\xi}\bigg)
-\frac{ly^2}{\xi} W_R'\bigg( \frac{y^2}{\xi} \bigg)\bigg) J_{l+2}(y) \, dy.
\end{split}
\end{equation}
Similarly,
\begin{equation} \notag
\begin{split}
\int_0^{\infty} \frac{y}{\xi}  W_R'\bigg( \frac{y^2}{\xi}\bigg)J_{l+2}(y) \, dy&=\int_0^{\infty} \frac{y^{-l-2}}{\xi}  W_R'\bigg( \frac{y^2}{\xi}\bigg)d(y^{l+3}J_{l+3}(y))  \\
&=\int_0^{\infty} \frac{1}{y^2} \bigg( \frac{2y^4}{\xi^2}  W_R^{''}\bigg( \frac{y^2}{\xi}\bigg)
-(l+2) \frac{y^2}{\xi} W_R'\bigg( \frac{y^2}{\xi} \bigg)\bigg) J_{l+3}(y) \, dy.
\end{split}
\end{equation}
Thus repeatedly integrating by parts, and arguing as above one sees that
\begin{equation*}
\begin{split}
\mathcal B_{l}( W_R)(\tfrac{\xi}{4\pi^2}) \ll& \frac{1}{\xi} \int_0^{\infty} y^{-A+1}|J_{l+A}(y)| \sum_{j=0}^A \bigg(\frac{y^2}{\xi} \bigg)^j (|l|+1)^{A-j} \bigg|W_R^{(j)}\bigg(\frac{y^2}{\xi} \bigg) \bigg| \, dy \\
\ll & \frac{1 }{\xi^{A/2} R} \sum_{j=0}^A R^j (|l|+1)^{A-j} \int_{\mathbb R} (1+y/R)^{j-A/2} |J_{l+A}(\sqrt{\xi(1+\tfrac{y}{R}})| |W^{(j)}(y)| \, dy,
\end{split}
\end{equation*}
where in the second step we made the change of variables $y \rightarrow \sqrt{\xi(1+\tfrac{y}{R})}$ and used that $W_R^{(l)}(y) =R^l W^{(l)}((y-1)R)$.
Noting that $J_l(x)=\int_{-1/2}^{1/2} e(lt) e^{-ix \sin 2\pi t} \, dt$, we see $|J_l(x)|\le 1$.
Additionally, we have that $J_{l}(x) \ll 1/\sqrt{x}$ for $x \ge 2 l$ \cite[8.453.1]{GR}. Using these estimates in the above equation and bounding $\int_{\mathbb R} |W^{(j)}(y)| \, dy$ completes the proof.
\end{proof}

\subsection{Proof of Proposition \ref{prop:equid}}
\begin{proof}[Proof of Proposition \ref{prop:equid}]
Let $W_{R,2}^{\pm}$ and $f_3^{\pm}$ be as in \eqref{eq:Wbd} and \eqref{eq:f3bd}, respectively. Recall from \eqref{eq:jpdef} that $N=|c-d|^{-1}$, $c'=cN, d'=dN$ and $J'=[c',d')$.
By construction,
\begin{equation} \label{eq:sandwich}
\begin{split}
\sum_{\lambda \in \mathbb Z^2} F_{N,3}^-(\theta_{\lambda}) W_{R,2}^-\left(
\frac{\|\lambda\|^2}{R^2}\right) &\le
 \# \bigg\{ \lambda \in \Gamma_R : \widetilde \theta_{\lambda}  \Mod 1 \in [ c,d ] \, \& \, r_{\lambda} \in I \bigg\} \\&\le \sum_{\lambda \in \mathbb Z^2} F_{N,3}^+(\theta_{\lambda}) W_{R,2}^+\bigg(
\frac{\|\lambda\|^2}{R^2}\bigg).
\end{split}
\end{equation}
In what follows, the r.h.s. of \eqref{eq:sandwich} will be bounded above, whereas the analogous steps give a lower bound for the l.h.s. of \eqref{eq:sandwich}.
Applying \eqref{eq:fourierseries} and Lemma \ref{lem:voronoi} we have that
\begin{equation}
\label{eq:transformed}
\begin{split}
&\sum_{\lambda \in \mathbb Z^2} F_{N,3}^+(\theta_{\lambda}) W_{R,2}^+\bigg(
\frac{\|\lambda\|^2}{R^2}\bigg)\\&=\frac{\pi R}{N} \widehat W_2^+(0) \widehat f_3^+(0)
+\frac{\pi R^2}{N} \sum_{\ell \in \mathbb Z} \widehat f_3^+\bigg( \frac{\ell}{N} \bigg)  \sum_{n \ge 1} \tau_{4\ell}(n) \mathcal B_{4\ell}(W_{R,2}^+)(R^2 n).
\end{split}
\end{equation}
Note that $\widehat f_3(x) \ll P^A/(1+|x|)^A$ for any nonnegative integer $A$.
By this and Lemma \ref{lem:IBP} the sum over $\ell$ in \eqref{eq:transformed} can be truncated at $|\ell| \le R^{\varepsilon} PN$ up to a negligible error term of size $\ll R^{-100}$, and by Lemma \ref{lem:IBP}
the sum over $n$ is effectively truncated at $n \le R^{\varepsilon}P^2$.

Applying the second estimate in
Lemma \ref{lem:IBP} with $A=1$ we get uniformly for $|\ell| \le R^{\varepsilon} PN$ with $P \le R^{1/2-\varepsilon}$ (so that $PN \le R^{1-\varepsilon}$)
\begin{equation} \label{eq:boundedsum11}
\sum_{n \ge 1} \tau_{4\ell}(n) \mathcal B_{4\ell}(W_{R,2}^+)(R^2 n) \ll \frac{1}{R^{3/2-\varepsilon}} \sum_{1 \le n \le R^{\varepsilon}P^2} \frac{1}{n^{3/4}} \ll \frac{P^{1/2}}{R^{3/2-\varepsilon}}.
\end{equation}
Using \eqref{eq:f3bd} it is not hard to see that
\begin{equation} \label{eq:boundedsum12}
\widehat f_3^+(\xi) \ll \min\bigg\{ 1, \frac{1}{|\xi|} \bigg\}+\frac{1}{P}.
\end{equation}
The estimates \eqref{eq:boundedsum11} and \eqref{eq:boundedsum12} yield the bound
\begin{equation*}
\begin{split}
&\frac{\pi R^2}{N} \sum_{\ell \in \mathbb Z} \widehat f_3^+\bigg( \frac{\ell}{N} \bigg)  \sum_{n \ge 1} \tau_{4\ell}(n) \mathcal B_{4\ell}(W_{R,2}^+)(R^2 n)
\\&\ll \frac{R^{1/2+\varepsilon}}{N}P^{1/2}  \bigg(
\sum_{|\ell|\le N} 1+N \sum_{N < |\ell| \le R^{\varepsilon} PN} \frac{1}{|\ell|} +\frac1P \sum_{|\ell| \le R^{\varepsilon}PN}  1 \bigg) \ll  R^{1/2+\varepsilon}P^{1/2}
\end{split}
\end{equation*}
for the second term on the r.h.s. of \eqref{eq:transformed}, using which, and upon recalling
\eqref{eq:Wbd} and \eqref{eq:f2bd}, we conclude that
\begin{equation*}
\begin{split}
\sum_{\lambda \in \mathbb Z^2} F_{N,3}^+(\theta_{\lambda}) W_{R,3}^+\bigg(
\frac{\|\lambda\|^2}{R^2}\bigg)= 2 \sqrt{2} \pi R |c-d|\cdot \frac{|I|}{\sqrt{2}}
+O \bigg(  \frac{R}{NP}+R^{1/2+\varepsilon}P^{1/2}\bigg).
\end{split}
\end{equation*}
Take $P=R^{1/3}/N^{2/3}$ to balance error terms. Combining this along with a completely analogous estimate for $ \sum_{\lambda \in \mathbb Z^2} F_{N,3}^-(\theta_{\lambda}) W_{R,3}^-\left(
\frac{\|\lambda\|^2}{R^2}\right)$, and recalling \eqref{eq:sandwich} completes the proof of Proposition \ref{prop:equid}. \end{proof}

\subsection{Preliminary estimates for the proof of Proposition \ref{thm:paircorrelation}}
We begin by
applying \eqref{eq:fourierseries} to see that
\begin{equation} \label{eq:fourier}
\begin{split}
&\sum_{\lambda,\mu \in \mathbb Z^2}  F_{\pazocal K,1}(\theta_{\lambda}-\theta_{\mu}) F_{1,2}(\theta_{\lambda}) W_{R,1}\bigg( \frac{\|\lambda\|^2}{R^2}\bigg)W_{R,2}\bigg( \frac{\|\mu\|^2}{R^2}\bigg)\\
&\qquad =\frac{1}{\pazocal K} \sum_{\ell,m \in \mathbb Z} \widehat f_1\bigg(\frac{\ell}{\pazocal K} \bigg) \widehat f_2(m )  \sum_{\lambda \in \mathbb Z^2} e^{4i(\ell+m)\theta_{\lambda}} W_{R,1}\bigg( \frac{\|\lambda\|^2}{R^2}\bigg) \sum_{\mu \in \mathbb Z^2}e^{4i\ell \theta_{\mu}}W_{R,2}\bigg( \frac{\|\mu\|^2}{R^2}\bigg).
\end{split}
\end{equation}
To evaluate the sums over $\lambda,\mu$ our first step is to apply Lemma \ref{lem:voronoi} and prove:
\begin{lemma} \label{lem:summationapplied}
Let $\varepsilon>0$. We have that
\[
\begin{split}
&\sum_{\lambda,\mu \in \mathbb Z^2}  F_{\pazocal K,1}(\theta_{\lambda}-\theta_{\mu}) F_{1,2}(\theta_{\lambda}) W_{R,1}\bigg( \frac{\|\lambda\|^2}{R^2}\bigg)W_{R,2}\bigg( \frac{\|\mu\|^2}{R^2}\bigg)
= \frac{\pi R}{2\sqrt{2} } \widehat W_1(0) \widehat W_2(0)
\widehat f_1(0)\widehat f_2(0)\\
&+\frac{\pi R^3}{2\sqrt{2} } \bigg( \sum_{\substack{\ell \in \mathbb Z \\ }} \sum_{\substack{m \in \mathbb Z \\ }} \bigg( \sum_{1 \le n_1 \le P^2 R^{\varepsilon} } \tau_{4\ell+4m}(n_1) \mathcal B_{4\ell+4m}(W_{R,1})(R^2n_1)\bigg) \\
& \qquad \qquad  \times \bigg( \sum_{1 \le n_2 \le P^2 R^{\varepsilon} } \tau_{4\ell}(n_2) \mathcal B_{4\ell}(W_{R,2})(R^2n_2)\bigg) \widehat f_1 \bigg( \frac{\ell}{\pazocal K}\bigg) \widehat f_2 ( m)+O\bigg(P^{5/2}R^{1/2+\varepsilon}\bigg).
\end{split}
\]

\end{lemma}

\begin{proof}
Applying Lemma \ref{lem:voronoi} we see that the r.h.s. of \eqref{eq:fourier} is
\begin{equation} \label{eq:vorapplied}
\begin{split}
=\frac{\pi R^3}{2\sqrt{2}}\sum_{\ell,m \in \mathbb Z}
\widehat f_1  \bigg( \frac{\ell}{\pazocal K}\bigg) \widehat f_2 ( m)
&
\bigg( \delta_{\ell=-m} \widehat W_{R,1}(0) +\sum_{n_1\ge 1} \tau_{4\ell+4m}(n_1)\mathcal B_{4\ell+4m}(W_{R,1})(R^2n_1)\bigg) \\
 \times  & \bigg(\delta_{\ell=0} \widehat W_{R,2}(0) +\sum_{n_2\ge 1} \tau_{4\ell}(n_2)\mathcal B_{4\ell}(W_{R,2})(R^2n_2) \bigg).
\end{split}
\end{equation}
Using \eqref{eq:f1bd} and Lemma \ref{lem:IBP} it is not hard to see that the sum over $\ell$ can be truncated to $|\ell| \le P R^{1+\varepsilon}$
up to an error term of size $O(R^{-100})$, which is negligible and similarly the sum over $m$ can be truncated to $|m| \le P R^{\varepsilon}$. After having truncated the sums over $m,\ell$, we now
apply Lemma \ref{lem:IBP} and get that the contribution from the terms in the sum over $n_1$ with $n_1 > P^2 R^{\varepsilon}$ is $\ll R^{-100}$, which is negligible. Similarly, we can truncate the sum over $n_2$ to $n_2 \le P^2 R^{\varepsilon}  $ at the cost of a negligible error term of size $\ll R^{-100}$.

Recall $\widehat W_{R,1}(0) \ll 1/R$. Also, for $x \ge 2|m|$ we have $J_{m}(x) \ll 1/\sqrt{x}$ so that for $\xi \ge 10 m^2$
\[
\begin{split}
\mathcal B_{4m}(W_{R,2})(\xi)= &\int_0^{\infty} W_{R,2}(y) J_{4m}(2\pi \sqrt{\xi y }) \, dy \\
\ll& \frac1R \int_0^{\infty}| W_2(y) J_{4m}(2\pi \sqrt{\xi(1+\tfrac{y}{R})})| \, dy \ll \frac{1}{\xi^{1/4} R}.
\end{split}
\]
Hence, recalling also $P \le R^{1-\varepsilon}$ we have that
\begin{equation} \label{eq:vorapplied11}
\begin{split}
&\sum_{\ell,m \in \mathbb Z} \delta_{\ell=-m} \cdot
\widehat f_1  \bigg( \frac{\ell}{\pazocal K}\bigg) \widehat f_2 ( m)  \widehat W_{R,1}(0) \sum_{n_2\ge 1} \tau_{4\ell}(n_2)\mathcal B_{4\ell}(W_{R,2})(R^2n_2)  \\
&\qquad \qquad \ll \frac1R \sum_{|m| \le PR^{\varepsilon}}
 \sum_{n_2\le P^2R^{\varepsilon}} |\tau_{4m}(n_2)\mathcal B_{4m}(W_{R,2})(R^2n_2)|+R^{-100} \\
&\qquad \qquad \ll \frac{P}{R^{2-\varepsilon}}     \sum_{n_2 \le P^2 R^{\varepsilon}} \frac{1}{ (R^2n_2)^{1/4}} \ll \frac{P^{5/2}}{R^{5/2-\varepsilon}}  .
\end{split}
\end{equation}
By a similar argument we have that
\begin{equation} \label{eq:vorapplied12}
\begin{split}
\sum_{\ell,m \in \mathbb Z} \delta_{\ell=0} \cdot
\widehat f_1  \bigg( \frac{\ell}{\pazocal K}\bigg) \widehat f_2 (m)  \widehat W_{R,2}(0) \sum_{n_1\ge 1} \tau_{4\ell+4m}(n_1)\mathcal B_{4\ell+4m}(W_{R,1})(R^2n_1) \ll \frac{P^{5/2}}{R^{5/2-\varepsilon}} .
\end{split}
\end{equation}
Recall that in \eqref{eq:vorapplied} the inner sums over $n_1,n_2$ are effectively truncated at $n_1,n_2 \le P^2 R^{\varepsilon}$; hence using \eqref{eq:vorapplied11} and \eqref{eq:vorapplied12} in \eqref{eq:vorapplied} completes the proof of Lemma \ref{lem:summationapplied}.
\end{proof}

\subsection{Summing over frequencies}
Let us introduce the following notation.
For $b_1, b_2 \Mod 4$, and nonzero $\lambda,\mu \in \mathbb Z^2 $ let
\begin{equation} \label{eq:phidef}
\begin{split}
\phi(t)&=\phi(t;b_1,b_2,\lambda,\mu) \\
&= \|\lambda\| \sin(2\pi ( \tfrac{t}{4}+\tfrac{b_2}{4})+\theta_{\lambda})
- \|\mu\| \sin(2\pi ( \tfrac{t}{4}+\tfrac{b_2-b_1}{4})-\theta_{\mu}).
\end{split}
\end{equation}
Also, for $v_1,v_2 \in \mathbb R$, $b_1,b_2 \Mod 4$, and nonzero $\lambda,\mu \in \mathbb Z^2$ let
\begin{equation} \label{eq:adef}
\begin{split}
a(t,v_1,v_2)&=a(t,v_1,v_2;b_1,b_2, \lambda,\mu)\\
&=e \bigg(  \frac{v_1}{2}\|\lambda\| \sin(2\pi ( \tfrac{t}{4}+\tfrac{b_2}{4})+\theta_{\lambda})-
  \frac{v_2}{2}\|\mu\|\sin(2\pi ( \tfrac{t}{4}+\tfrac{b_2-b_1}{4})-\theta_{\mu}) \bigg)\\
& \qquad \qquad \times \widehat f_1\bigg(\frac{\pi}{2 \sqrt{2}} \|\mu\| \cos( 2\pi ( \tfrac{t}{4}+\tfrac{b_2-b_1}{4})-\theta_{\mu}) \bigg)
f_2(-t).
\end{split}
\end{equation}

In this section we will establish the following result.
\begin{lemma} \label{lem:summed}
Let $\varepsilon>0$.
We have that
\begin{equation} \label{eq:summedup}
\begin{split}
&\sum_{\substack{\ell \in \mathbb Z \\ }} \sum_{\substack{m \in \mathbb Z \\ }} \bigg( \sum_{n_1 \le R^{\varepsilon} P^2} \tau_{4\ell+4m}(n_1) \mathcal B_{4\ell+4m}(W_{R,2})(R^2n_1)\bigg) \\
& \qquad \qquad \qquad \qquad \qquad \times \bigg( \sum_{n_2 \le R^{\varepsilon} P^2} \tau_{4\ell}(n_2) \mathcal B_{4\ell}(W_{R,2})(R^2n_2)\bigg) \widehat f_1 \bigg( \frac{\ell}{\pazocal K}\bigg) \widehat f_2 (m) \\
&= \frac{1}{16R^2} \sum_{0\le b_1,b_2  <4}
 \sum_{\substack{\lambda, \mu \in \mathbb Z^2 \setminus \{0\} \\ \|\lambda\|,\|\mu\| < R^{\varepsilon} P}}  \int_{\mathbb R^3} \bigg( W_{1}(v_1) W_2(v_2) a(t,v_1,v_2;b_1,b_2, \lambda,\mu) \\
& \qquad \qquad \qquad \qquad \qquad \times e(R \phi(t;b_1,b_2,\lambda,\mu)) \bigg) \, dt \, dv_1 dv_2+O\bigg( \frac{k^2 P^5}{R^{3-\varepsilon}} \bigg).
\end{split}
\end{equation}

\end{lemma}
The main step in the proof of the preceding lemma is the following result in which we apply Poisson summation over $\ell,m$.

\begin{lemma} \label{lem:frequencies}
Let $N>0$, and let $g_1,g_2: \mathbb R \rightarrow \mathbb C$ be Schwartz functions.
For $\theta_1,\theta_2 \, \Mod{\tfrac14}$ and $x,y \in \mathbb R_{ > 0}$ we have that
\begin{equation}\label{eq:frequencies}
\begin{split}
&\sum_{\substack{\ell \in \mathbb Z \\ 4|\ell}} \sum_{\substack{m \in \mathbb Z \\ 4|m}} e(\ell \theta_1+m\theta_2) J_{\ell+m}(2\pi x)J_{\ell}(2\pi y) g_1\bigg( \frac{\ell}{4N}\bigg) g_2\bigg( \frac{m}{4}\bigg) \\
&=\frac{1}{16} \sum_{0 \le b_1,b_2 < 4}  \int_{\mathbb R} \int_{\mathbb R}   e\bigg(x \sin\bigg(2\pi \bigg(\frac{t_2}{4}+\theta_2+\frac{b_2}{4}\bigg)\bigg)\bigg)\\
& \qquad \qquad \qquad \times
e\bigg(y \sin\bigg(2\pi\bigg(\frac{t_1}{4N}-\frac{t_2}{4}+\theta_1-\theta_2+\frac{b_1-b_2}{4}\bigg)\bigg)\bigg) \widehat g_1(t_1) \widehat g_2(t_2) \, dt_1 dt_2.
\end{split}
\end{equation}
\end{lemma}


\begin{proof}
Recall for each integer $j$ that
\begin{equation}
\label{eq:Bessel Jj int}
J_{j}(2\pi x)=\int_{-1/2}^{1/2} e(-x \sin(2\pi t)) \cdot e(j t) \, dt.
\end{equation}
Using the fact that $\frac14 \sum_{0 \le b < 4} e(b\ell/4)$ is one if $4|\ell$ and is zero otherwise, along
with the identity \eqref{eq:Bessel Jj int}, for each $m\in\Z$ we have that
\begin{equation}
\label{eq:sum BesselJ subs int}
\begin{split}
&\sum_{\substack{\ell \in \mathbb Z \\ 4|\ell}}  e(\ell \theta_1) J_{\ell+m}(2\pi x)J_{\ell}(2\pi y) g_1\bigg( \frac{\ell}{4N}\bigg)
\\&=\frac14 \sum_{0 \le b_1 < 4} \int_{-1/2}^{1/2} \int_{-1/2}^{1/2} e(-x \sin(2\pi t_2)-y\sin(2\pi t_1)) e(mt_2)
  \\
& \qquad \qquad \qquad \times
\sum_{\ell \in \mathbb Z} e(\ell(t_1+t_2+\theta_1+\tfrac{b_1}{4})) g_1\bigg(\frac{\ell}{4N} \bigg)dt_1dt_2 \\
  &=N \sum_{0 \le b_1 < 4} \int_{-1/2}^{1/2} \int_{-1/2}^{1/2} e(-x \sin(2\pi t_2)-y\sin(2\pi t_1))  e(mt_2)
  \\
& \qquad \qquad \qquad \times \sum_{j \in \mathbb Z} \widehat g_1(4N(j-(t_1+t_2+\theta_1+\tfrac{b_1}{4})) \, dt_1 dt_2,
\end{split}
\end{equation}
where in the last step we applied Poisson summation. We now make the linear transformation of variables $t_1 \rightarrow -t_1-t_2-\theta_1-\tfrac{b_1}{4}+j$ in the inner integral on the r.h.s. of \eqref{eq:sum BesselJ subs int} to get that
\begin{equation}
\label{eq:ellsum}
\begin{split}
&\sum_{\substack{\ell \in \mathbb Z \\ 4|\ell}}  e(\ell \theta_1) J_{\ell+m}(2\pi x)J_{\ell}(2\pi y) g_1\bigg( \frac{\ell}{4N}\bigg) =N \sum_{0 \le b_1 < 4} \int_{-1/2}^{1/2}   e(-x \sin(2\pi t_2)) e(mt_2) \\
& \qquad \qquad \qquad \times \sum_{j \in \mathbb Z} \int_{-\frac12+j-t_2-\theta_1-\frac{b_1}{4}}^{\frac12+j-t_2-\theta_1-\frac{b_1}{4}} e(y \sin(2\pi(t_1+t_2+\theta_1+\tfrac{b_1}{4})) \widehat g_1(4Nt_1) \, dt_1 dt_2\\
&=N \sum_{0 \le b_1 < 4} \int_{-1/2}^{1/2} \int_{\mathbb R}   e(-x \sin(2\pi t_2)+y \sin(2\pi(t_1+t_2+\theta_1+\tfrac{b_1}{4})) e(mt_2) \widehat g_1(4Nt_1) \, dt_1 dt_2.
\end{split}
\end{equation}

Summing up the identity \eqref{eq:ellsum}, multiplied by $e(m\theta_2) \cdot  g_2( \frac{m}{4}) $, w.r.t. $m$, yields the equality
\begin{equation}
\label{eq:P1}
\begin{split}
&\sum_{\substack{\ell \in \mathbb Z \\ 4|\ell}} \sum_{\substack{m \in \mathbb Z \\ 4|m}} e(\ell \theta_1+m\theta_2) J_{\ell+m}(2\pi x)J_{\ell}(2\pi y) g_1\bigg( \frac{\ell}{4N}\bigg) g_2\bigg( \frac{m}{4}\bigg)
\\=&N \sum_{0 \le b_1 < 4} \int_{-1/2}^{1/2} \int_{\mathbb R}   e(-x \sin(2\pi t_2)+y \sin(2\pi(t_1+t_2+\theta_1+\tfrac{b_1}{4})) \\
& \qquad \qquad \qquad \times \bigg( \sum_{\substack{m \in \mathbb Z \\ 4|m}} e(m(\theta_2+t_2)) g_2\bigg(\frac{m}{4} \bigg) \bigg)\widehat g_1(4Nt_1) \, dt_1 dt_2.
\end{split}
\end{equation}
To evaluate the inner sum over $m$ on the r.h.s. of \eqref{eq:P1} we argue similarly, to see that
\begin{equation}
\label{eq:sum w.r.t. m PSF}
\sum_{\substack{m \in \mathbb Z \\ 4|m}} e(m(\theta_2+t_2)) g_2\bigg(\frac{m}{4} \bigg)= \sum_{0 \le b_2 <4} \sum_{j \in \mathbb Z} \widehat g_2(4(j-(t_2+\theta_2+\tfrac{b_2}{4}))).
\end{equation}
Substituting \eqref{eq:sum w.r.t. m PSF} into \eqref{eq:P1}, and making the change of variables $t_2 \rightarrow -t_2-\theta_2-\tfrac{b}{4}+j$ in the integral over $t_2$ yields the equality
\[
\begin{split}
&\sum_{\substack{\ell \in \mathbb Z \\ 4|\ell}} \sum_{\substack{m \in \mathbb Z \\ 4|m}} e(\ell \theta_1+m\theta_2) J_{\ell+m}(2\pi x)J_{\ell}(2\pi y) g_1\bigg( \frac{\ell}{4N}\bigg) g_2\bigg( \frac{m}{4}\bigg)
\\&=N \sum_{0 \le b_1,b_2 < 4} \sum_{j \in \mathbb Z} \int_{-\frac12-\theta_2-\frac{b_2}{4}+j}^{\frac12-\theta_2-\frac{b_2}{4}+j} \int_{\mathbb R}   e(x \sin(2\pi (t_2+\theta_2+\tfrac{b_2}{4}))\\
& \qquad \qquad \qquad \times e(y \sin(2\pi(t_1-t_2+\theta_1-\theta_2+\tfrac{b_1-b_2}{4})) \widehat g_1(4Nt_1) \widehat g_2(4t_2) \, dt_1 dt_2\\
&=\frac{1}{16} \sum_{0 \le b_1,b_2 < 4}  \int_{\mathbb R} \int_{\mathbb R}   e\bigg(x \sin\bigg(2\pi \bigg(\frac{t_2}{4}+\theta_2+\frac{b_2}{4}\bigg)\bigg)\bigg)\\
& \qquad \qquad \qquad \times
e\bigg(y \sin\bigg(2\pi\bigg(\frac{t_1}{4N}-\frac{t_2}{4}+\theta_1-\theta_2+\frac{b_1-b_2}{4}\bigg)\bigg)\bigg) \widehat g_1(t_1) \widehat g_2(t_2) \, dt_1 dt_2,
\end{split}
\]
which establishes \eqref{eq:frequencies} and completes the proof of Lemma \ref{lem:frequencies}.
\end{proof}

\begin{proof}[Proof of Lemma \ref{lem:summed}]
Recall for $\ell \in \mathbb Z$ and $\xi>0$ that $\mathcal B_l(g)(\xi)=\int_0^{\infty} g(y) J_l(2\pi \sqrt{\xi y}) \, dy$ and $\tau_{4l}(n)=\sum_{\lambda \in \mathbb Z^2 : \| \lambda \|^2=n} e^{4il \theta_{\lambda}}$, which we substitute into the
l.h.s. of \eqref{eq:summedup} to yield
\begin{equation} \label{eq:initial}
\begin{split}
&\sum_{\substack{\ell \in \mathbb Z \\ }} \sum_{\substack{m \in \mathbb Z \\ }} \bigg( \sum_{n_1 \le R^{\varepsilon} P^2} \tau_{4\ell+4m}(n_1) \mathcal B_{4\ell+4m}(W_{R,2})(R^2n_1)\bigg) \\
& \qquad \qquad \qquad \qquad \qquad \times \bigg( \sum_{n_2 \le R^{\varepsilon} P^2} \tau_{4\ell}(n_2) \mathcal B_{4\ell}(W_{R,2})(R^2n_2)\bigg) \widehat f_1 \bigg( \frac{\ell}{\pazocal K}\bigg) \widehat f_2 (m)
\\&=\sum_{\substack{\ell \in \mathbb Z \\ }} \sum_{\substack{m \in \mathbb Z \\ }} \bigg( \sum_{1 \le \| \lambda \| \le R^{\varepsilon} P}  e^{4i(\ell+m)\theta_{\lambda}} \mathcal B_{4\ell+4m}(W_{R,2})(R^2 \| \lambda \|^2)\bigg) \\
& \qquad \qquad \qquad \qquad \qquad \times \bigg( \sum_{1 \le \| \mu \| \le R^{\varepsilon} P} e^{4i\ell\theta_{\mu}}  \mathcal B_{4\ell}(W_{R,2})(R^2\| \mu \|^2)\bigg) \widehat f_1 \bigg( \frac{\ell}{\pazocal K}\bigg) \widehat f_2 (m)\\
&=
\sum_{\substack{\lambda, \mu \in \mathbb Z^2 \\ 1 \le
\|\lambda\|,\|\mu\| \le R^{\varepsilon}P}} \int_0^{\infty}
\int_0^{\infty} W_{R,1}(v_1)W_{R,2}(v_2) \bigg( \sum_{\ell, m \in \mathbb Z} \widehat f_1\bigg( \frac{\ell}{\pazocal K}\bigg) \widehat f_2(m) \\
& \qquad \qquad \qquad \times
e^{4i\ell(\theta_{\lambda}+\theta_{\mu})} e^{4im \theta_{\lambda}} J_{4\ell+4m}(2\pi \|\lambda\| R \sqrt{ v_1}) J_{4\ell}(2\pi \|\mu\| R \sqrt{ v_2}) \bigg) dv_1 dv_2.
\end{split}
\end{equation}
Applying Lemma \ref{lem:frequencies}, with  $N=\pazocal K$, $\theta_1=\tfrac{1}{2\pi}(\theta_{\lambda}+\theta_{\mu})$, $\theta_2=\tfrac{\theta_{\lambda}}{2\pi}$, $x=\|\lambda\|R \sqrt{v_1}$, and $y=\|\mu\|R\sqrt{v_2}$, on the inner sum on the r.h.s. of \eqref{eq:initial}, we obtain the equality
\begin{equation}
\label{eq:innersum}
\begin{split}
&\sum_{\ell, m \in \mathbb Z} \widehat f_1\bigg( \frac{\ell}{\pazocal K}\bigg) \widehat f_2(m)
e^{4i\ell(\theta_{\lambda}+\theta_{\mu})} e^{4im \theta_{\lambda}} J_{4\ell+4m}(2\pi \|\lambda\| R \sqrt{ v_1}) J_{4\ell}(2\pi \|\mu\| R \sqrt{ v_2})
\\&=\frac{1}{16} \sum_{0 \le b_1,b_2 < 4} \int_{\mathbb R} \int_{\mathbb R} e(\|\lambda\|R \sqrt{v_1} \sin(2\pi(\tfrac{t_2}{4}+\tfrac{b_2}{4})+\theta_{\lambda})) \\
& \qquad \qquad \times e(\|\mu\|R \sqrt{v_2} \sin(2\pi(\tfrac{t_1}{4\pazocal K}-\tfrac{t_2}{4}+\tfrac{b_1-b_2}{4})+\theta_{\mu})) f_1(-t_1)f_2(-t_2) \, dt_1 dt_2.
\end{split}
\end{equation}

For brevity, let $h=\tfrac{-t_2}{4}+\frac{\theta_{\mu}}{2\pi}+\frac{b_1-b_2}{4}$. Plainly,
\begin{equation}
\label{eq:taylor}
\sin(2\pi( h+\tfrac{t_1}{4\pazocal K}))=\sin(2\pi h)+\frac{\pi^2 t_1}{\pazocal K} \cos(2\pi h) +O(t_1^2 R^{-2}).
\end{equation}
We also recall the estimates $e^{iu}=1+O(|u|)$, $u \in \mathbb R$ and $\sqrt{1+u}=1+u/2+O(u^2)$, $u>-1/2$.
Substituting \eqref{eq:taylor} into \eqref{eq:innersum}, and using the second estimate in \eqref{eq:f1bd} with $l=2$ shows that
\begin{equation}
\label{eq:id12}
\begin{split}
&\sum_{\ell, m \in \mathbb Z} \widehat f_1\bigg( \frac{\ell}{\pazocal K}\bigg) \widehat f_2(m)
e^{4i\ell(\theta_{\lambda}+\theta_{\mu})} e^{4im \theta_{\lambda}} J_{4\ell+4m}(2\pi \|\lambda\| R \sqrt{ v_1}) J_{4\ell}(2\pi \|\mu\| R \sqrt{ v_2})
\\&=\frac{1}{16} \sum_{0 \le b_1,b_2 < 4} \int_{\mathbb R} e(\|\lambda\|R \sqrt{v_1} \sin(2\pi(\tfrac{t_2}{4}+\tfrac{b_2}{4})+\theta_{\lambda})+\|\mu\|R \sqrt{v_2} \sin(2\pi h))  \\
& \qquad \qquad \times \bigg(
\int_{\mathbb R} e(\tfrac{\pi^2 t_1}{\pazocal K} \|\mu\| R\sqrt{v_2} \cos(2\pi h))  f_1(-t_1) dt_1 \bigg) f_2(-t_2) \, dt_2+O\bigg(\frac{k^2 \sqrt{v_2}\| \mu \|}{R} \bigg).
\end{split}
\end{equation}

Recall $\pazocal K=2 \sqrt{2} \pi R$. Observe that the inner integral of the r.h.s. of \eqref{eq:id12} is
\begin{equation}
\label{eq:id11}
\int_{\mathbb R} e\left(\tfrac{\pi^2 t_1}{\pazocal K} \|\mu\| R\sqrt{v_2} \cos(2\pi h)\right)  f_1(-t_1) dt_1
=\widehat f_1\left(\tfrac{\pi\sqrt{v_2}}{2\sqrt{2}} \|\mu\| \cos(2\pi h)\right).
\end{equation}
Substituting \eqref{eq:id11} into \eqref{eq:id12}, and then into \eqref{eq:initial}, and transforming the variables
$(v_j-1)R \rightarrow v_j$ for $j=1,2$ we have that
\begin{equation}
\label{eq:nearlythere}
\begin{split}
&\sum_{\substack{\ell \in \mathbb Z \\ }} \sum_{\substack{m \in \mathbb Z \\ }} \bigg( \sum_{n_1 \le R^{\varepsilon} P^2} \tau_{4\ell+4m}(n_1) \mathcal B_{4\ell+4m}(W_{R,2})(R^2n_1)\bigg) \\
& \qquad \qquad \qquad \qquad \qquad \times \bigg( \sum_{n_2 \le R^{\varepsilon} P^2} \tau_{4\ell}(n_2) \mathcal B_{4\ell}(W_{R,2})(R^2n_2)\bigg) \widehat f_1 \bigg( \frac{\ell}{\pazocal K}\bigg) \widehat f_2 (m)
\\&=\frac{1}{16R^2} \sum_{0 \le b_1,b_2 < 4} \sum_{\substack{\lambda, \mu \in \mathbb Z^2 \\ 1 \le
\|\lambda\|,\|\mu\| \le R^{\varepsilon}P}} \int_{\mathbb R}
\int_{\mathbb R} W_{1}(v_1)W_{2}(v_2)
\int_{\mathbb R} \widetilde a(t,v_1,v_2)e(R \phi(t)) dt \, dv_1 dv_2+O\bigg(\frac{k^2 P^5}{R^{3-\varepsilon}} \bigg)
\end{split}
\end{equation}
where $\phi$ is as given in \eqref{eq:phidef}, and
\begin{equation} \label{eq:awidedef}
\begin{split}
\widetilde  a(t,v_1,v_2)
&=e \bigg( \bigg(\sqrt{1+\frac{v_1}{R}}-1\bigg) R\|\lambda\| \sin(2\pi ( \tfrac{t}{4}+\tfrac{b_2}{4})+\theta_{\lambda})+
  \bigg(\sqrt{1+\frac{v_2}{R}}-1\bigg) R\|\mu\|\sin(2\pi h) \bigg)\\
  & \qquad \qquad \qquad   \times \widehat f_1\bigg(\frac{\pi\sqrt{(1+\frac{v_2}{R})}}{2\sqrt{2}} \|\mu\| \cos(2\pi h) \bigg)  f_2(-t).
\end{split}
\end{equation}
Also, we have that
\begin{equation} \label{eq:mvtf1}
\widehat f_1\bigg(\frac{\pi\sqrt{(1+\frac{v_2}{R})}}{2\sqrt{2}} \|\mu\| \cos(2\pi h) \bigg)
=\widehat f_1\bigg(\frac{\pi}{2\sqrt{2}} \|\mu\| \cos(2\pi h) \bigg)
+O\bigg(\frac{|v_2| \|\mu\|}{R} \int_{\mathbb R} |t f(t)| \, dt \bigg)
\end{equation}
and by \eqref{eq:f1bd} the error term is $\ll R^{\varepsilon} P k/ R$ for $\|\mu\| \le R^{\varepsilon}P$.
Recall that $a$ is defined in \eqref{eq:adef}. Using \eqref{eq:awidedef} and \eqref{eq:mvtf1}, we have for $\|\lambda\|, \|\mu\| \le R^{\varepsilon}P$ that
\begin{equation} \label{eq:finalesta}
\widetilde  a(t,v_1,v_2)=  a(t,v_1,v_2)+
O\bigg(R^{\varepsilon}|f_2(-t)|\bigg( \frac{P(k+v_1^2+v_2^2)}{R} \bigg) \bigg).
\end{equation}
Applying this estimate in \eqref{eq:nearlythere} completes the proof of Lemma \ref{lem:summed}.
\end{proof}

\section{Proofs of propositions \ref{thm:paircorrelation}, \ref{thm:varasymp}, and Lemma \ref{lem:stationaryphaseapplied2}} \label{sec:lattice-proof2}
In this section we apply
the method of stationary phase to the oscillatory integral on the r.h.s. of \eqref{eq:summedup} as well as the integral appearing in Lemma \ref{lem:stationaryphaseapplied2}. We then modify this argument to prove Proposition \ref{thm:varasymp}.

\subsection{Stationary phase estimates}
We begin by quoting the following estimates for oscillatory integrals due to Blomer, Khan, and Young \cite{BKY}, which are uniform in multiple parameters. We first quote the following result which is special case of \cite[Lemma 8.1]{BKY}.

\begin{lemma} \label{lem:stationary1}
Let $Y  \ge 1$ and $Q,V,S>0$, and assume that $w$ is a smooth function supported on a closed interval $J \subseteq \mathbb R$ with
\[
w^{(j)}(t) \ll V^{-j},
\]
and $h$ is a smooth function on $J$ such that
\[
|h'(t)| \ge S, h^{(j)}(t) \ll_j Y Q^{-j}, \quad j=2,3,\ldots.
\]
Then, for any $A>0$,
\[
\int_{\mathbb R} w(t) e^{ih(t)} \, dt \ll_A |J|((QS/\sqrt{Y})^{-A}+(SV)^{-A}).
\]
\end{lemma}
We also require the following stationary phase estimate, which is uniform in multiple parameters \cite[Proposition 8.2]{BKY}.
\begin{lemma}\label{lem:stationary2}
Let $0<\delta<1/10$. Let $Y, V, U ,Q>0$. Set $Z=Q+U+Y+1$. Suppose that
\[
Y \ge Z^{3\delta}, U \ge V \ge \frac{Q Z^{\delta/2}}{\sqrt{Y}}.
\]
Additionally, suppose that $w$ is a smooth function
on $\mathbb R$
that is supported on an interval $J$ of length $U$, satisfying
\[
w^{(j)}(t) \ll_j V^{-j}
\]
for all nonnegative integers $j$. Suppose $h$ is a smooth function on $J$ such that there exists a unique point $t_0 \in J$ such that $h'(t_0)=0$, and furthermore
\[
h''(t) \gg Y Q^{-2}, \quad h^{(j)}(t) \ll_j Y Q^{-j}, \quad \text{for } j\in \mathbb N, t\in J.
\]
Then for any fixed $A>0$ we have that
\[
\int_{\mathbb R} w(t) e^{ih(t)} \,dt=\frac{e^{ih(t_0)+i\frac{\pi}{4}\tmop{sgn}(h''(t_0))}}{\sqrt{|h''(t_0)|}}\sum_{0 \le n \le 3 A/\delta} p_n(t_0)+O(Z^{-A})
\]
where for each nonnegative integer $n$ and $t_0 \in \mathbb R$
\[
p_n(t_0)=\frac{\sqrt{2\pi}}{n!} \bigg( \frac{i}{2h''(t_0)}\bigg)^n
G^{(2n)}(t_0),
\]
and
\[
G(t)=w(t)e^{iH(t)}, \quad H(t)=h(t)-h(t_0)-\frac12 h''(t_0) (t-t_0)^2.
\]
Furthermore, each $p_n$ is a rational function in $h'',h''',\ldots,$ satisfying
\begin{equation} \label{eq:pbd}
\frac{d^j}{dt_0^j}  p_n(t_0) \ll_{j,n} (V^{-j}+Q^{-j})((V^2Y/Q^2)^{-n}+Y^{-n/3}).
\end{equation}
\end{lemma}
The leading term $n=0$ in the asymptotic expansion \eqref{eq:pbd} is given by
\begin{equation} \label{eq:leadingterm}
\bigg(\frac{2\pi}{|h''(t_0)|}\bigg)^{1/2} e^{i \frac{\pi}{4} \tmop{sgn}(h''(t_0))}  e^{ih(t_0)} w(t_0) \ll \frac{Q}{Y^{1/2}}.
\end{equation}



\begin{lemma} \label{lem:stationaryphaseapplied}
Let $\varepsilon>0$.
Let $\mathcal M$ be as given in \eqref{eq:Mdef}, $a$ be as in \eqref{eq:adef}, and $\phi$ be as in \eqref{eq:phidef}.
Suppose that $P^2k \le R^{1/2-\delta_0}$ for some fixed $\delta_0>0$. Then
\begin{equation} \label{eq:diagonal}
\begin{split}
 &\frac{1}{16} \sum_{0\le b_1,b_2  <4}
 \sum_{\substack{\lambda, \mu \in \mathbb Z^2 \setminus \{0\} \\ \|\lambda\|,\|\mu\| < R^{\varepsilon}P}}  \int_{\mathbb R^3}  W_{1}(u_1) W_2(u_2)   a(t,u_1,u_2;b_1,b_2, \lambda,\mu) e(R \phi(t;b_1,b_2,\lambda,\mu)) \, dt \, du_1 du_2\\
&\qquad \qquad \qquad \qquad \qquad \qquad \qquad =8 \mathcal M(W_1,W_2,f_1,f_2)+O\bigg(\frac{P^{7/2}R^{\varepsilon}}{\sqrt{R}}\bigg).
 \end{split}
\end{equation}
\end{lemma}

The proof of Lemma \ref{lem:stationaryphaseapplied} gives a full asymptotic expansion of the error term (see \eqref{eq:stationary-expansion}).

\begin{proof}
To clarify the following computation we will work over the complex numbers. We have for nonzero $z,w \in \mathbb C$ with $z=|z|e^{i\theta_z}, w=|w| e^{i\theta_w}$ that
\[
\begin{split}
&|z| \sin(2\pi ( \tfrac{t}{4}+\tfrac{b_2}{4})+\theta_{z})
- |w| \sin(2\pi ( \tfrac{t}{4}+\tfrac{b_2-b_1}{4})-\theta_{w})\\
&\qquad \qquad \qquad =\frac{z e(\tfrac{t}{4}) i^{b_2}-\overline z e(\tfrac{-t}{4}) i^{-b_2}}{2i}-\bigg( \frac{ \overline w e(\tfrac{t}{4}) i^{b_2-b_1}- w e(\tfrac{-t}{4}) i^{b_1-b_2}}{2i} \bigg) \\
& \qquad \qquad \qquad =\frac{1}{2i} \bigg((z-i^{-b_1} \overline w)i^{b_2} e(\tfrac{t}{4})-\overline{(z-i^{-b_1} \overline w)i^{b_2} e(\tfrac{t}{4})} \bigg) \\
&\qquad \qquad \qquad = \tmop{Im}((z-i^{-b_1} \overline w) i^{b_2}e(\tfrac{t}{4})).
\end{split}
\]
Given $\lambda=(a,b),\mu=(c,d) \in \mathbb R^2$ we let
\begin{equation} \label{eq:gamma-def}
\gamma=\gamma(\lambda,\mu,b_1,b_2)=(\tmop{Re}(((a+ib)-i^{-b_1}  (c-id) )i^{b_2}),\tmop{Im} (((a+ib)-i^{-b_1}  (c-id) )i^{b_2})).
\end{equation}
Writing $\theta_{\gamma}$ for the angle between the positive real axis and the vector corresponding to $\gamma$ we see that
\begin{equation} \label{eq:phase}
\phi(t)=\phi(t;b_1,b_2,\lambda,\mu)=\| \gamma\| \sin(\tfrac{\pi t}{2}+\theta_{\gamma})
\end{equation}
and in particular
\begin{equation} \label{eq:stationarypts}
\phi'(t)=\frac{\pi\| \gamma\|}{2} \cos(\tfrac{\pi t}{2}+\theta_{\gamma}), \qquad \phi''(t)=\frac{-\pi^2\| \gamma\|}{4} \sin(\tfrac{\pi t}{2}+\theta_{\gamma})
\end{equation}
so that $\phi'(t_j)=0$ if and only if
\begin{equation}\label{eq:tj-def}
t_j=1+2j-\tfrac{2\theta_{\gamma}}{\pi}
\end{equation}
for some $j \in\mathbb Z$. Also,
\begin{equation} \label{eq:phipp}
\phi''(t_j)=(-1)^{j+1}\pi^2\| \gamma\|/4.
\end{equation}
We now separately analyze the contribution from the terms on the l.h.s. of \eqref{eq:diagonal} with $\gamma=0$ (diagonal) and $\gamma \neq 0$ (off-diagonal) cases.
Recall $h(t,\lambda,b)=2\pi(\tfrac{t}{4}+\tfrac{b}{4})+\theta_{\lambda}$.
The contribution to the sum on the l.h.s. of \eqref{eq:diagonal} over nonzero $\lambda,\mu \in \mathbb Z^2$ with $\gamma=0$ is
\begin{equation} \label{eq:secondary}
\begin{split}
&=\frac{1}{16} \sum_{0 \le b_1,b_2 < 4} \sum_{\substack{\lambda,\mu \in \mathbb Z^2 \\ 1 \le \|\lambda\|,\|\mu\| < R^{\varepsilon}P \\ \mu=\overline{\lambda} i^{-b_1}}
} \int_{\mathbb R} \int_{\mathbb R} W_1(v_1)W_2(v_2) \int_{\mathbb R} a(t,u_1,u_2)\, dt dv_1 dv_2\\
&=\frac{1}{16} \sum_{0 \le b_1,b_2 < 4} \sum_{\substack{\lambda,\mu \in \mathbb Z^2 \\ 1 \le \|\lambda\|,\|\mu\| < R^{\varepsilon}P \\ \mu=\overline{\lambda} i^{-b_1}}} \int_{\mathbb R^3} W_1(v_1) W_2(v_2)
e ( \tfrac12(v_1-v_2) \|\lambda\| \sin(2\pi ( \tfrac{t}{4}+\tfrac{b_2}{4})+\theta_{\lambda})) \\
&\qquad \qquad \qquad \qquad \qquad \qquad \times \widehat f_1\bigg(\frac{\pi}{2\sqrt{2}} \|\mu\| \cos( 2\pi ( \tfrac{t}{4}+\tfrac{b_2-b_1}{4})-\theta_{\mu}) \bigg) f_2(-t) \, dt dv_1 dv_2 \\
&= \frac{1}{4} \sum_{0 \le b < 4} \sum_{\substack{\lambda \in \mathbb Z^2 \\ 1 \le \|\lambda\| < R^{\varepsilon} P }} \int_{\mathbb R} \widehat W_1\bigg(- \tfrac12 \|\lambda\| \sin(h(t,\lambda,b))\bigg)\widehat W_2\bigg( \tfrac12 \|\lambda\| \sin(h(t,\lambda,b))\bigg) \\
& \qquad \qquad \qquad \times \widehat f_1\bigg(\frac{\pi}{2\sqrt{2}} \|\lambda\| \cos(h(t,\lambda,b)) \bigg) f_2(-t) \, dt.
\end{split}
\end{equation}
Using the rapid decay of $\widehat V_1,\widehat f_1$ we can extend the inner sum on the r.h.s. of \eqref{eq:secondary} to a sum over all $\lambda \in \mathbb Z^2$ at the cost of a negligible error term of size $O(R^{-100})$.

It remains to treat the terms on the l.h.s. of \eqref{eq:diagonal} with $\gamma \neq 0$.
Using a smooth partition of unity $\sum_j \varrho_j(t)=1$ we rewrite the integral
\[
\int_{\mathbb R} a(t,u_1,u_2) e(R \phi(t)) \, dt
=\sum_{l} \int_{\mathbb R} a(t,u_1,u_2) e(R \phi(t))  \varrho_l(t)\, dt.
\]
We choose the smooth functions $\varrho_l$ so that each is supported on an interval $[u_l,v_l]$ with $|u_l-v_l|\le 1/2$ so that $[u_l,v_l]$ contains at most one stationary point (cf. \eqref{eq:tj-def}) $ t_j=1+2j-\tfrac{2\theta_{\gamma}}{\pi}$
and the distance between each of $u_l,v_l$ and the corresponding nearest stationary point is  $\ge 1/100$
 and $\varrho_l^{(\ell)}(t) \ll_l 1$.
 Also, due to the rapid decay of $f_2$ we can assume that each stationary point $t_j$ satisfies $|t_j| \le R^{\varepsilon}$.
 Splitting
$a$ into its real and imaginary parts, we will next apply Lemma \ref{lem:stationary1} or Lemma \ref{lem:stationary2}
with $w=\varrho_l \cdot \tmop{Re}(a)$, $w=\varrho_l \cdot \tmop{Im}(a)$, and $h= 2\pi R \phi$, respectively depending on whether the support of $\varrho_l$
contains a stationary point, and recalling \eqref{eq:f1bd}, \eqref{eq:stationarypts} we may take
\begin{equation} \label{eq:paramchoice}
V=\frac{1}{(\|\lambda\|+\|\mu\|)Pk}, \quad Y=\|\gamma\| R, \quad  Q=1, \quad S=\frac{\|\gamma\| R}{200},
\end{equation}
say.
By construction we may also take $U=1/2$.
Recalling \eqref{eq:phase} and \eqref{eq:stationarypts},
and $\|\lambda\|,\|\mu\| \le R^{\varepsilon}P$,
we see that for each integer $l$, the hypothesis of either Lemma \ref{lem:stationary1} or Lemma \ref{lem:stationary2}
is satisfied by taking $\delta$ sufficiently small in terms of $\delta_0$ since $V \gg (P^2k R^{\varepsilon})^{-1}$ and $P^2 k \le R^{1/2-\delta_0}$. If we are in the case where $\varrho_l$ does not contain a stationary point we observe that $QS/\sqrt{Y} \gg R^{1/2}, SV \gg R^{1/2}$ so that by Lemma \ref{lem:stationary1} the contribution from this case is $\ll R^{-100}$.

It remains to treat the case where $\varrho_l$ contains a stationary point and we will apply Lemma \ref{lem:stationary2}. We note for $\|\lambda\|,\|\mu\| < R^{\varepsilon}P$ that $V^2Y/Q^2 \gg \| \gamma \| R^{1-\varepsilon}/(P^4k^2)$ and $Y^{1/3} \gg  (\| \gamma \| R)^{1/3}$. Hence,  using \eqref{eq:pbd} with $j=0$, \eqref{eq:leadingterm}, and \eqref{eq:phipp} we conclude that
uniformly for $|v_1|,|v_2| \le 2$, nonzero $\lambda,\mu \in \mathbb Z^2$ with $\|\lambda\|,\|\mu\| < R^{\varepsilon}P$ and $\gamma \neq 0$ that
\begin{equation} \label{eq:stationary-expansion}
\begin{split}
&\int_{\mathbb R} a(t,u_1,u_2) e(R \phi(t)) \, dt \\
&=\sum_{|j| \le R^{\varepsilon}}
\bigg(\frac{4}{\pi^2 R \|\gamma\|} \bigg)^{1/2} e^{-i (-1)^{j}\frac{\pi}{4}} e(R\phi(t_j))a(t_j)+
O\bigg(\frac{R^{\varepsilon}}{\sqrt{R\|\gamma\|}} \bigg( \frac{P^4k^2}{R\|\gamma\|}+\frac{1}{(R\|\gamma\|)^{1/3}} \bigg) \bigg).
\end{split}
\end{equation}
Hence the contribution of the $\lambda,\mu $ with $\gamma \neq 0$ to the l.h.s. of \eqref{eq:diagonal} is
\begin{equation} \label{eq:gammabd}
\ll \frac{1}{R^{1/2-\varepsilon}} \sum_{0 \le b_1,b_2 <4} \sum_{\substack{\lambda,\mu \in \mathbb Z^2 \\ 1 \le \|\lambda\|,\|\mu\| < R^{\varepsilon}P \\ \mu \neq \overline{\lambda} i^{-b_1}}} \frac{1}{\| \gamma\|^{1/2}} \ll \frac{1}{R^{1/2-\varepsilon}}
\sum_{1 \le \| \gamma\| \le R^{\varepsilon}P} \frac{1}{\sqrt{\| \gamma\|}} \sum_{\substack{\lambda,\mu \in \mathbb Z^2 \\ 1 \le \|\lambda\|,\|\mu\| < R^{\varepsilon}P \\ \gamma=i^{b_2}(\mu i^{b_1} + \overline{\lambda}) }} 1 \ll \frac{P^{7/2}R^{\varepsilon}}{\sqrt{R}},
\end{equation}
where the last estimate follows by noting that for each $\lambda$ there is one $\mu$ for which $\gamma=i^{b_2}(\mu i^{b_1} + \overline{\lambda})$.
Combining the estimate \eqref{eq:gammabd} with \eqref{eq:secondary} concludes the proof of Lemma \ref{lem:stationaryphaseapplied}.
\end{proof}

\subsection{Proof of Proposition \ref{thm:paircorrelation}}
\begin{proof}[Proof of Proposition \ref{thm:paircorrelation}]
Combining the results of Lemmas \ref{lem:summationapplied},\ref{lem:summed}, and \ref{lem:stationaryphaseapplied} we have that
\[
\begin{split}
&\sum_{\lambda,\mu \in \mathbb Z^2}  F_{\pazocal K,1}(\theta_{\lambda}-\theta_{\mu}) F_{1,2}(\theta_{\lambda}) W_{R,1}\bigg( \frac{\|\lambda\|^2}{R^2}\bigg)W_{R,2}\bigg( \frac{\|\mu\|^2}{R^2}\bigg)\\
&= \frac{\pi R}{2 \sqrt{2}} \widehat W_1(0) \widehat W_2(0) \widehat f_1(0)\widehat f_2(0)+ \frac{4\pi R}{\sqrt{2}} \mathcal M(W_1,W_2,f_1,f_2)+O(R^{\varepsilon}(P^{7/2} R^{1/2}+k^2 P^5))
\end{split}
\]
where $\mathcal M$ is as given in \eqref{eq:Mdef}. Recalling that $\pazocal K=2\pi\sqrt{2} R$ completes the proof.
\end{proof}

\subsection{Proof of Lemma \ref{lem:stationaryphaseapplied2}} \label{sec:lemproof}
\begin{proof}[Proof of Lemma \ref{lem:stationaryphaseapplied2}]
Recall $h(t,\lambda,b)=2\pi(\tfrac{t}{4}+\tfrac{b}{4})+\theta_{\lambda}$.
Since $f_1(x)=g(x-k)$ we have that
\[
\widehat f_1\bigg(\frac{\pi}{2\sqrt{2}} \|\lambda\| \cos(h(t,\lambda,b))\bigg)=e\bigg(-k\frac{\pi}{2\sqrt{2}} \|\lambda\| \cos(h(t,\lambda,b))\bigg) \widehat g\bigg( \frac{\pi}{2\sqrt{2}} \|\lambda\| \cos(h(t,\lambda,b)) \bigg).
\]
We now view $k$ as a large parameter so that $I$ is an oscillatory integral, which can be estimated using the method of stationary phase.
Taking
\[
\varphi(t)=-k \frac{2\pi^2}{2\sqrt{2}} \|\lambda\| \cos(2\pi(\tfrac{t}{4}+\tfrac{b}{4})+\theta_{\lambda})
\]
we see that
\[
\varphi'(t)=k \frac{\pi^3}{2\sqrt{2}} \|\lambda\| \sin(2\pi(\tfrac{t}{4}+\tfrac{b}{4})+\theta_{\lambda})
\]
and
\[
\varphi''(t)=k \frac{\pi^4}{4\sqrt{2}} \|\lambda\| \cos(2\pi(\tfrac{t}{4}+\tfrac{b}{4})+\theta_{\lambda}).
\]
Hence, $I$ will have stationary points
at $t_j=2j-b-\tfrac{2\theta_{\lambda}}{\pi}$, $j\in\mathbb Z$.
In particular,
\begin{equation} \label{eq:varphipp}
\varphi''(t_j)=(-1)^{j}k \frac{\pi^4}{4\sqrt{2}}.
\end{equation}
Additionally,
take
\[
b(t)= \widehat W_1\bigg(-\frac{\|\lambda\|}{2} \sin(h(t,\lambda,b))\bigg) \widehat W_2\bigg(\frac{\|\lambda\|}{2} \sin(h(t,\lambda,b))\bigg) \widehat g\bigg( \frac{\pi}{2\sqrt{2}} \|\lambda\| \cos(h(t,\lambda,b) \bigg) f_2(-t).
\]
Arguing as in the proof of Lemma \ref{lem:stationaryphaseapplied}, we use a smooth partition of unity $1=\sum_l \varrho_l$ and apply Lemma \ref{lem:stationary1}
or Lemma \ref{lem:stationary2} with $h=\varphi$, $w=\varrho_l\tmop{Re}(b)$ or $w=\varrho_j\tmop{Im}(b)$ depending on whether the support of $w$ contains a stationary point. By the rapid decay of $f_2$ we can assume the stationary points $t_j$ satisfy $|t_j| \le k^{\varepsilon}$. We take our parameters as follows:
\[
V=\frac{1}{\|\lambda\|+P}, \quad Y=k\|\lambda\|, \quad  Q=1, \quad S=\frac{k \|\lambda\| }{600},
\]
and $U=1/2$. In the case where $\varrho_l$ does not contain a stationary point we note $(QS/\sqrt{Y}) \gg k^{1/2}$ and $SV \gg k/P \gg k^{1/2}$, since $P^{3+3\delta_1} \le k$, so that by Lemma \ref{lem:stationary1} the contribution from this case is $\ll k^{-100}$.
In the case where $\varrho_l$ contains a stationary point, we apply Lemma \ref{lem:stationary2}, the hypotheses are satisfied since we may take $\delta$ sufficiently small in terms of $\delta_1$, $k>P^{3+3\delta_1}$ and $\|\lambda \|\le k^{\varepsilon} P$ so that $V \ge Q Z^{\delta/2}/\sqrt{Y}$ for $\delta$ sufficiently small in terms of $\delta_1$. Also, note that $V^2Y/Q \gg k/(P^2 \| \lambda \|)$ and $Y^{1/3} \gg (k\| \lambda | )^{1/3}$.
Hence, using \eqref{eq:pbd}, \eqref{eq:leadingterm}, and \eqref{eq:varphipp} we can conclude that

\begin{equation} \label{eq:stationary-expansion2}
\begin{split}
I&=\int_{\mathbb R} b(t) e^{i\varphi(t)} \, dt \\
&=\sum_{|j| \le k^{\varepsilon}} \bigg(\frac{8 \sqrt{2}}{\pi^3 k\|\lambda\|} \bigg)^{1/2} e^{\frac{\pi}{4}i(-1)^j} e^{i\varphi(t_j)} b(t_j)+O\bigg(\frac{1}{\sqrt{k\|\lambda\|}}\bigg(\frac{\|\lambda\|P^2}{k}+ \frac{1}{(k\|\lambda\|)^{1/3}}\bigg) \bigg) .
\end{split}
\end{equation}
The sum on the r.h.s. of \eqref{eq:stationary-expansion2}
is $\ll k^{\varepsilon} (\|\lambda\|k)^{-1/2}$, which establishes the claim.
\end{proof}

\subsection{Proof of Proposition \ref{thm:varasymp}}
The key ingredient in the proof of Proposition \ref{thm:varasymp} is the following result.
\begin{proposition} \label{prop:var} Let $c,d \in \mathbb R$ with $c<d$. Also, let $N=|d-c|^{-1}$.
Suppose that $ P^{5/2} R^{1/2+\delta_1} \le N \le 100R$ where $\delta_1>0$ is fixed. Then we have that
\[
\begin{split}
&\int_0^{\pi/2} \bigg| \sum_{\substack{ \lambda \in \mathbb Z^2}} F_{N,2}(\theta_{\lambda}-\theta) W_{R,2}\bigg(\frac{\|\lambda\|^2}{R^2} \bigg)-\frac{\pi R}{N} \widehat W_2(0) \widehat f_2(0) \bigg|^2 \frac{d\theta}{\pi/2}\\
&\qquad \qquad \qquad =\frac{\pi^2 R^2}{N^2} \sum_{1 \le \|\lambda \| \le R^{\varepsilon} P} \int_{-1/2}^{1/2} \bigg| \widehat W_2\bigg( \frac{\|\lambda \|}{2} \sin(2\pi t) \bigg)
\widehat f_2\bigg( \frac{\pi^2 R}{N} \|\lambda \| \cos(2\pi t) \bigg) \bigg|^2 \, dt\\
&\qquad \qquad\qquad \qquad \qquad+ O\bigg( \frac{R^{2+\varepsilon}P^4}{N^2} \bigg( \frac{RP}{N^2}+\frac{1}{R^{1/2}}\bigg)\bigg).
\end{split}
\]
\end{proposition}

\begin{proof}
The proof closely follows the arguments given in the proof of Theorem \ref{thm:paircorrelation} and we will give a detailed sketch of the argument. For brevity we write $f=f_2$, $F_N=F_{N,2}$, $W=W_2$ and $W_{R,2}=W_R$. Use \eqref{eq:fourierseries} to expand $F_N$, then apply Lemma \ref{lem:voronoi} and recall $\widehat W_R(0)=R^{-1} \widehat W(0)$ to get that
\begin{equation} \label{eq:expand1}
\begin{split}
&\sum_{\lambda \in \mathbb Z^2} F_N(\theta_{\lambda}-\theta) W_R\bigg(\frac{\|\lambda\|^2}{R^2} \bigg)
=\frac{\pi R}{N} \widehat W(0) \widehat f(0)\\
&\qquad \qquad \qquad \qquad \qquad \qquad+\frac{\pi R^2}{N} \sum_{\ell \in \mathbb Z}
 \widehat f\bigg( \frac{\ell}{N} \bigg) e^{-4i\ell \theta} \sum_{\lambda \in \mathbb Z^2 \setminus\{0\}} e^{4i\ell \theta_{\lambda}}
\mathcal B_{4\ell}(W_R)(R^2 \|\lambda\|^2).
\end{split}
\end{equation}
We now integrate over $\theta$ to obtain
\begin{equation} \label{eq:initialmanip}
\begin{split}
&\int_{0}^{\pi/2}\bigg|\sum_{\lambda \in \mathbb Z^2} F_N(\theta_{\lambda}-\theta) W_R\bigg(\frac{\|\lambda\|^2}{R^2} \bigg)
- \frac{\pi^2 R}{N} \widehat W(0) \widehat f(0)\bigg|^2 \frac{d\theta}{\pi/2} \\
&\qquad \qquad \qquad =\frac{\pi^2 R^4}{N^2} \sum_{\ell \in \mathbb Z}
 \bigg|\widehat f\bigg( \frac{\ell}{N} \bigg)\bigg|^2  \sum_{\mu,\lambda \in \mathbb Z^2 \setminus\{0\}} e^{4i\ell (\theta_{\lambda}-\theta_{\mu})}
\mathcal B_{4\ell}(W_R)(R^2 \|\lambda\|^2) \mathcal B_{4\ell}(W_R)(R^2 \|\mu\|^2).
\end{split}
\end{equation}

Recalling Lemma \ref{lem:IBP} the sums over $\lambda, \mu \in \mathbb Z^2$ are effectively truncated at $\|\lambda\|,\|\mu\| \le R^{\varepsilon} P$ and we will use this fact later on without further mention. Let $0 \le v_1,v_2 \le 4$. We now sum over $\ell$ by applying \eqref{eq:ellsum} with $g_1=|\widehat f|^2$, $x=R\|\lambda\|\sqrt{v_1}$, $y=R\|\mu\|\sqrt{v_1}$, $m=0$, $\theta_1=\tfrac{1}{2\pi}(\theta_{\lambda}-\theta_{\mu})$ and make a linear change of variables to see that
\begin{equation}
\label{eq:sum fhat BesselJ}
\begin{split}
&\sum_{\ell \in \mathbb Z}
 \bigg|\widehat f\bigg( \frac{\ell}{N} \bigg)\bigg|^2  e^{4i\ell (\theta_{\lambda}-\theta_{\mu})} J_{4\ell}(2\pi R\|\lambda\|\sqrt{v_1})
 J_{4\ell}(2\pi R\|\mu\|\sqrt{v_2})\\
 &\qquad \qquad = \frac{1}{4} \sum_{0 \le b < 4} \int_{-\frac{1}{2}+\frac{\theta_{\lambda}}{2\pi}}^{\frac{1}{2}+\frac{\theta_{\lambda}}{2\pi}} \int_{\mathbb R} e\bigg(-R \|\lambda\| \sqrt{v_1} \sin(2\pi t_2-\theta_{\lambda})\bigg)  \\
 & \qquad \qquad \qquad \times e\bigg(R\|\mu\| \sqrt{v_2} \sin(2\pi (\tfrac{t_1}{4N}+t_2+\tfrac{b}{4})-\theta_{\mu})  \bigg) (f \ast  f)(-t_1) \, dt_1 dt_2.
\end{split}
\end{equation}
Write $h_1=t_2+\tfrac{b}{4}-\tfrac{\theta_{\mu}}{2\pi}$. Taylor expanding $\sin(2\pi(\tfrac{t_1}{4N}+h_1))$ around $h_1$ (cf. \eqref{eq:taylor}) we see that up to an error term of size $\ll \sqrt{v_2} \| \mu \| R/N^2$, the l.h.s. of
\eqref{eq:sum fhat BesselJ} equals
\[
\frac{1}{4} \sum_{0 \le b < 4} \int_{-\frac{1}{2}+\frac{\theta_{\lambda}}{2\pi}}^{\frac{1}{2}+\frac{\theta_{\lambda}}{2\pi}} e\bigg(-R \|\lambda\| \sqrt{v_1} \sin(2\pi t_2-\theta_{\lambda})+R\|\mu\| \sqrt{v_2} \sin(2\pi h_1) \bigg)
\bigg|\widehat f\bigg( \frac{\pi^2 R \|\mu\| \sqrt{v_2}}{N} \cos(2\pi h_1) \bigg)\bigg|^2  dt_2.
\]
Using this expression to evaluate the sum over $\ell$ on the r.h.s. of \eqref{eq:initialmanip} we see that the l.h.s. of \eqref{eq:initialmanip} is
\begin{equation} \label{eq:transformedagain}
\begin{split}
&=\frac{\pi^2 R^4}{4N^2} \sum_{0 \le b < 4}
\sum_{\mu,\lambda \in \mathbb Z^2 \setminus\{0\}} \int_{0}^{\infty}\int_{0}^{\infty} W_R(v_1)W_R(v_2)
\int_{-\frac{1}{2}+\frac{\theta_{\lambda}}{2\pi}}^{\frac{1}{2}+\frac{\theta_{\lambda}}{2\pi}} e\bigg(-R \|\lambda\| \sqrt{v_1} \sin(2\pi t_2-\theta_{\lambda})\bigg) \\
&\times e\bigg(R\|\mu\| \sqrt{v_2} \sin(2\pi h_1) \bigg)
\bigg|\widehat f\bigg( \frac{\pi^2 R \|\mu\| \sqrt{v_2}}{N} \cos(2\pi h_1) \bigg)\bigg|^2  dt_2+O\bigg(\frac{R^{3+\varepsilon} P^5 }{N^4} \bigg).
\end{split}
\end{equation}

Make the linear change of variables $(v_j-1)R \rightarrow v_j$ for $j=1,2$, and $t_2 \rightarrow -t_2$ in the integrals above.
Let $\phi_1(t)=\|\lambda\| \sin(2\pi t_2+\theta_{\lambda})-\|\mu\|\sin(2\pi(t_2+\tfrac{b}{4})+\theta_{\mu})$ and
\[
a_1(t)= e\bigg( \frac{v_1}{2 } \|\lambda\| \sin(2\pi t_2+\theta_{\lambda})-\frac{v_2}{2} \|\mu\| \sin(2\pi(t_2+\tfrac{b}{4})
+\theta_{\mu})\bigg)  \bigg|\widehat f\bigg( \frac{\pi^2 R \|\mu\|}{N} \cos(2\pi h_1) \bigg)\bigg|^2.
\]
Note that $|\widehat f(x+y)|^2=|\widehat f(x)|^2+O(|y|)$ for $x,y\in \mathbb R$ with $|y|\le 1$ and recall $\sqrt{1+u}=1+u/2+O(u^2)$ for $u>-1/2$. Using the above estimates in \eqref{eq:transformedagain} (cf. \eqref{eq:awidedef}-\eqref{eq:finalesta}) we conclude that the l.h.s. of \eqref{eq:initialmanip} is
\begin{equation} \label{eq:oscform}
\begin{split}
=\frac{\pi^2 R^2}{4N^2} \sum_{0 \le b < 4} \sum_{\mu,\lambda \in \mathbb Z^2 \setminus \{0\}} \int_{\mathbb R^2}
W(v_1)W(v_2) \int_{-\frac{1}{2}-\frac{\theta_{\lambda}}{2\pi}}^{\frac{1}{2}-\frac{\theta_{\lambda}}{2\pi}} e(R\phi_1(t)) a_1(t) \, dt dv_1 dv_2
+O\bigg(\frac{R^{2+\varepsilon} P^5}{N^2} \bigg(\frac{R}{N^2}+\frac{1}{N} \bigg) \bigg).
\end{split}
\end{equation}
Note that the second error term is absorbed by the first.

We now analyze the integral over $t$ in \eqref{eq:oscform}. Recalling that $\phi_1(t)=\|\lambda\| \sin(2\pi t_2+\theta_{\lambda})-\|\mu\|\sin(2\pi(t_2+\tfrac{b}{4})+\theta_{\mu})$,
it is not difficult to see that $\phi_1(t)=\tmop{Im}(\gamma e(t))$ where $\gamma=\lambda-\mu i^b$. The cases $\gamma=0$ and $\gamma \neq 0$ are treated separately. After making a linear change of variables, switching the order of integration, and recalling the definition of $a_1$, the contribution of the terms with $\gamma=0$ to the main term in \eqref{eq:oscform} is
\begin{equation} \label{eq:gammazero}
= \frac{\pi^2 R^2}{N^2} \sum_{ \lambda \in \mathbb Z^2\setminus\{0\}}
 \int_{-1/2}^{1/2} \bigg| \widehat W\bigg( \frac{\|\lambda\|}{2} \sin(2\pi t) \bigg)
\widehat f\bigg( \frac{\pi^2 R}{N} \|\lambda\|\cos(2\pi t) \bigg) \bigg|^2 \, dt.
\end{equation}

Finally, we need to bound the contribution from the terms in \eqref{eq:oscform} with $\gamma\neq 0$. In this case, the integral
over $t$ is oscillatory and we use the method of stationary phase as in the proof of Theorem \ref{thm:paircorrelation}.
First we replace the indicator function
of $[-\tfrac12-\tfrac{\theta_{\lambda}}{2\pi}, \tfrac12-\tfrac{\theta_{\lambda}}{2\pi}]$ in the integrand over $t$ in \eqref{eq:oscform} by
a smooth function $v(t)$ which is equal to this indicator function outside of neighborhoods of length $\nu$ near the two boundary points, and we suppose $v^{(l)} \ll \nu^{-l}$. This step gives rise to an error of size $\ll R^{2+\varepsilon} P^4\nu/N^2$.
We then repeat the argument used to establish \eqref{eq:stationary-expansion} to bound the smoothed integral over $t$ in \eqref{eq:oscform} using Lemmas \ref{lem:stationary1} and \ref{lem:stationary2}, where one takes $Y,Q,S$ as in \eqref{eq:paramchoice}
and $V=N/(PR(\|\lambda\|+\|\mu\|))+\Delta$. Note that $P^2 R/N \le R^{1/2-\delta_1}$, so taking $\nu=R^{-1/2-\delta_3}$, where $\delta_3>0$ is sufficiently small, the hypotheses of Lemma \ref{lem:stationary2} are satisfied. Also, by \eqref{eq:leadingterm}
the leading order term in the stationary phase expansion is $ \ll \frac{1}{R^{1/2} \| \gamma\|^{1/2}}$. Hence arguing as in \eqref{eq:gammabd} we conclude that the contribution from the terms with $\gamma \neq 0$ to \eqref{eq:oscform} are
\[
\ll \frac{R^{3/2+\varepsilon}}{N^2} \sum_{0 \le b <4} \sum_{\substack{\lambda,\mu \in \mathbb Z^2 \\ 1 \le \|\lambda\|,\|\mu\| < R^{\varepsilon}P \\ \lambda \neq \mu i^{b}}} \frac{1}{\|\lambda-i^b \mu\|^{1/2}} \ll \frac{R^{3/2+\varepsilon} P^{7/2}}{N^2}.
\]
Using the estimate above along with \eqref{eq:gammazero} in \eqref{eq:oscform} completes the proof of Proposition \ref{prop:var}.
\end{proof}

\begin{proof}[Proof of Proposition \ref{thm:varasymp}]
Let $W_2^{\pm}$ and $f_3^{\pm}$ be as in \eqref{eq:Wbd} and \eqref{eq:f3bd}, respectively also recall that $N=|c-d|^{-1}$, $c'=cN$, $d'=dN$ and $J'=[c',d')$, as in \eqref{eq:jpdef}. By construction,
\begin{equation} \label{eq:sandwich1}
 \sum_{\substack{ \lambda \in \Gamma_R, r_{\lambda} \in I \\ \widetilde \theta_{\lambda}-\widetilde \theta \Mod 1 \in [ c,d ]  }} 1 \le \sum_{\lambda \in \mathbb Z^2} F_{N,3}^+(\theta_{\lambda}-\theta) W_{R,2}^+\bigg(
\frac{\|\lambda\|^2}{R^2}\bigg),
\end{equation}
and the l.h.s. of \eqref{eq:sandwich1} is
\begin{equation} \label{eq:sandwich2}
\ge \sum_{\lambda \in \mathbb Z^2} F_{N,3}^-(\theta_{\lambda}-\theta) W_{R,2}^-\bigg(
\frac{\|\lambda\|^2}{R^2}\bigg).
\end{equation}
Applying \eqref{eq:expand1}, we have that
\[
\int_0^{\pi/2}  \sum_{\lambda \in \mathbb Z^2} F_{N,3}^{\pm}(\theta_{\lambda}-\theta) W_{R,2}^{\pm}\bigg(
\frac{\|\lambda\|^2}{R^2}\bigg)\frac{d\theta}{\pi/2}= \frac{\pi R}{N} \widehat W_2^{\pm}(0) \widehat f_3^{\pm}(0)+\widehat f_3^{\pm}(0) \frac{\pi R^2}{N} \sum_{\lambda \in \mathbb Z^2 \setminus\{0\}}
\mathcal B_{0}(W_{R,2}^{\pm})(R^2 \|\lambda\|^2).
\]
Using Lemma \ref{lem:IBP}, including the second estimate in the lemma with $A=1$,
we see that $$
\frac{ R^2}{N} \sum_{\lambda \in \mathbb Z^2 \setminus\{0\}} \mathcal B_{0}(W_{R,2}^{\pm})(R^2 \|\lambda\|^2) \ll \frac{R^{2}}{N} \sum_{n \le P^2 R^{\varepsilon}} \frac{1}{R^{3/2}n^{3/4}} \ll \frac{R^{1/2+\varepsilon} P^{1/2}}{N}.$$
Additionally, recalling \eqref{eq:Wbd} gives $\widehat W_2^{\pm}(0)=2|I|+O(1/P)$ and also by \eqref{eq:f2bd} we have
$\widehat f_3^{\pm}(0)=1+O(1/P)$, so that
\begin{equation} \label{eq:expectation}
\int_0^{\pi/2}  \sum_{\lambda \in \mathbb Z^2} F_{N,3}^{\pm}(\theta_{\lambda}-\theta) W_{R,2}^{\pm}\bigg(
\frac{\|\lambda\|^2}{R^2}\bigg)\frac{d\theta}{\pi/2}=\frac{2 \sqrt{2}\pi R}{N} \frac{|I|}{\sqrt{2}}+O\bigg( \frac{R}{NP}+ \frac{R^{1/2+\varepsilon} P^{1/2}}{N}\bigg).
\end{equation}
Hence \eqref{eq:sandwich1},\eqref{eq:sandwich2}, and \eqref{eq:expectation} gives
\begin{equation} \label{eq:vareq}
\begin{split}
\int_0^{\pi/2} \bigg( \sum_{\substack{ \lambda \in \Gamma_R, r_{\lambda} \in I \\ \widetilde \theta_{\lambda}-\widetilde \theta \Mod 1 \in [ c,d ]  }} 1 \bigg) \frac{d\theta}{\pi/2} = \frac{2 \sqrt{2}\pi R}{N} \frac{|I|}{\sqrt{2}} +O\bigg( \frac{R}{N}\bigg(\frac1P+ \frac{ P^{1/2}}{R^{1/2-\varepsilon}}\bigg) \bigg).
\end{split}
\end{equation}
Using Proposition \ref{prop:var} along with \eqref{eq:sandwich1}, \eqref{eq:expectation}, and \eqref{eq:vareq} we get that
\begin{equation} \label{eq:varupperbd}
 \begin{split}
&\int_0^{\pi/2} \bigg|  \sum_{\substack{ \lambda \in \Gamma_R, r_{\lambda} \in I \\ \widetilde \theta_{\lambda}-\widetilde \theta \Mod 1 \in [ c,d ]  }} 1- \frac{2 \sqrt{2} \pi R }{N} \frac{|I|}{\sqrt{2}}\bigg|^2 \frac{d\theta}{\pi/2} \\
&\qquad \qquad \qquad \le \frac{\pi^2 R^2}{N^2} \sum_{1 \le n \le  P^{2+\varepsilon}} \int_{-1/2}^{1/2} \bigg| \widehat W_2^+\bigg( \frac{\sqrt{n}}{2} \sin(2\pi t) \bigg)
\widehat f_3^+\bigg( \frac{\pi^2 R}{N} \sqrt{n} \cos(2\pi t) \bigg) \bigg|^2 \, dt\\
&\qquad \qquad\qquad \qquad \qquad+ O\bigg( \frac{R^{2}P^4}{N^2} \bigg(  \frac{R^{1+\varepsilon}P}{N^2}+\frac{1}{P^5} +\frac{1}{R^{1/2-\varepsilon}}\bigg)\bigg),
\end{split}
\end{equation}
where we also used the rapid decay of $\widehat f_3,\widehat W_2$ to truncate the sum at $n \le P^{2+\varepsilon}$ up to an error term of $\ll R^2/(N^2P^{100})$.

Since $\widehat \chi_{J'}(\xi) \ll \min(1/|\xi|,1)$ we have
 for $\Delta_1, \Delta_2 \ge 1/100$, $a \ge 1/10$, and $|I|,|J| \ll 1$ that
\begin{equation} \label{eq:elementary}
\begin{split}
\int_{\mathbb R} | \widehat \chi_{I}(\Delta_1 \xi)| |\widehat \chi_{J}(\Delta_2 (\xi-a)) |^2 \, d\xi & \ll  \frac{1}{\Delta_2^2} \int_{{|\xi| \le \frac{1}{10^4\Delta_1}}} 1 \, d\xi
+ \frac{1}{\Delta_1} \int_{{|\xi-a| \le \frac{1}{10^4\Delta_2}}} 1 \, d\xi \\
& \qquad \qquad +\frac{1}{\Delta_1\Delta_2^2}
\int_{\substack{|\xi|>\frac{1}{10^4\Delta_1} \\ |\xi-a| > \frac{1}{10^4\Delta_2}}} \frac{1}{|\xi|(\xi-a)^2} \, d\xi \ll
\frac{\log (\Delta_1+1)}{\Delta_1 \Delta _2}.
\end{split}
\end{equation}
By \eqref{eq:Wbd} we have that $\widehat W_2^+(\xi)=2\widehat \chi_{I}(2\xi)+O(1/P)$ and by \eqref{eq:f3bd} we have that
$\widehat f_3^+(\xi)=\widehat \chi_{J'}(\xi)+O(1/P)$, for $\xi \in \mathbb R$. Hence, using \eqref{eq:elementary} together with the simple estimate $\int_{\mathbb R} |\chi_I(\Delta \xi)|^2 \, d \xi = |I|/\Delta$ for any $\Delta>0$, $I \subseteq \mathbb R$ we get that
\begin{equation} \label{eq:transfer}
\begin{split}
 &\int_{-1/2}^{1/2} \bigg| \widehat W_2^+\bigg( \frac{\sqrt{n}}{2} \sin(2\pi t) \bigg)
\widehat f_3^+\bigg( \frac{\pi^2 R}{N} \sqrt{n} \cos(2\pi t) \bigg) \bigg|^2 dt\\
& \qquad \qquad =4
 \int_{-1/2}^{1/2} \bigg| \widehat \chi_I( \sqrt{n}\sin(2\pi t) )
\widehat \chi_{J'}\bigg( \frac{\pi^2 R}{N} \sqrt{n} \cos(2\pi t) \bigg) \bigg|^2 dt+
O\bigg(\frac{1}{P^{1-\varepsilon} n} \bigg)
\end{split}
\end{equation}
for $n \le P^{2+\varepsilon}$. Use this in \eqref{eq:varupperbd} and sum over $n$ to obtain the claimed upper bound for the l.h.s. of \eqref{eq:varupperbd}. To get the claimed lower bound,
we argue similarly only now we replace $W_2^+$ and $f_3^+$ with $W_2^-$ and $f_3^-$. Combining these upper and lower bounds yields \eqref{eq:varform}.
\end{proof}

\bibliographystyle{abbrv}
\bibliography{references}

\end{document}